\documentclass[twoside,a4paper,backrefs,msc-links
]{amsart}
\usepackage{graphicx}
\PassOptionsToPackage{pdfauthor={Vladimir V. Kisil},%
    pdftitle={Erlangen Program at Large: Outline},%
    pdfsubject={general mathematics},%
    backref=page,%
    breaklinks=true,
    pdfkeywords={symmetries, Erlangen program}
}{hyperref}
\usepackage{hyperref}
\usepackage{amsrefs}
\usepackage{amssymb}
\input{cyr.fd}
\usepackage[all]{xy}
\IfFileExists{eulervm.sty}{\usepackage{eulervm}
\newcommand{\Aprime}{A\!'}}{\newcommand{\Aprime}{A'}}
\def\ifundefined#1{\expandafter\ifx\csname#1\endcsname\relax}
\makeindex
\newcommand{\amscite}[3]{\cite{#1}#2{#3}}

\numberwithin{equation}{section}
\newtheorem{thm}{Theorem}[section]
\newtheorem{prop}[thm]{Proposition}
\newtheorem{principle}[thm]{Principle}
\newtheorem{cor}[thm]{Corollary}
\theoremstyle{definition}
\newtheorem{defn}[thm]{Definition}
\newtheorem{example}[thm]{Example}
\theoremstyle{remark}
\newtheorem{rem}[thm]{Remark}

\PII{\relax}
\copyrightinfo{}{\relax}

\makeatletter
\providecommand{\diag}{\mathop {\operator@font diag}\nolimits}
\makeatother

\newcommand{\epigraph}[3]{\par
\hfill\parbox{0.6\textwidth}{\footnotesize #1 \par \hfil #2 
\textit{#3}}\par\medskip}%
\providecommand{\obj}[2][\,]{\ensuremath{\mathrm{#2}#1}}
\providecommand{\such}{\,\mid\,}
\providecommand{\norm}[2][\relax]{\left\|#2\right\|\ifx#1\relax\else_{#1}\fi}
\providecommand{\modulus}[2][\relax]{\left| #2 \right|\ifx#1\relax\else_{#1}\fi}
\providecommand{\lvec}[1]{\overrightarrow{#1}}
\providecommand{\uir}[3][0]{\ifcase #1{\rho^{#2}_{#3}}%
\or {\breve{\rho}^{#2}_{#3}}%
\or {\tilde{\rho}^{#2}_{#3}}\fi}
\providecommand{\cycle}[3][]{{#1 C^{#2}_{#3}}}
\newcommand{\zcycle}[3][]{#1 Z^{#2}_{#3}}
\newcommand{\realline}[3][]{#1 R^{#2}_{#3}}
\providecommand{\matr}[4]{{\ensuremath{ \left(\!\! \begin{array}{cc}
#1 & #2 \\ #3 & #4
\end{array}\!\!\right) }}}
\providecommand{\linv}[2][\relax]{\mathfrak{L}^{#2}_{#1}}

\providecommand{\spec}[1][]{\ensuremath{\mathbf{sp}}\,}

\providecommand{\bs}{\breve{\sigma}}
\providecommand{\Sp}[1][n]{\ensuremath{\mathrm{Sp}(#1)}}
\providecommand{\SL}[1][2]{\ensuremath{\FSpace{SL}{#1}(\Space{R}{})}}
\providecommand{\scalar}[3][\relax]{\left\langle #2,#3 
        \right\rangle\ifx#1\relax\else_{#1}\fi}
\providecommand{\Space}[3][]{\ensuremath{\mathbb{#2}^{#3}_{#1}{}}}
  \providecommand{\FSpace}[3][]{\ensuremath{\ifx#2l \ell_{#3}^{#1}{}\else
  #2_{#3}^{#1}{}\fi}} 
\providecommand{\rmi}{\mathrm{i}}
\providecommand{\rmh}{\mathrm{j}}
\providecommand{\rme}{\mathrm{e}}
\providecommand{\rmp}{\varepsilon}
\providecommand{\alli}{\iota}
\providecommand{\rmc}{\mathrm{\breve\i}}
\providecommand{\myh}{h}
\providecommand{\myhbar}{\hslash}
\providecommand{\tr}{\mathop{tr}}
\providecommand{\ladder}[2][]{L_{#1}^{\!#2}}
\providecommand{\algebra}[1]{\ensuremath{\mathfrak{#1}}}
\makeatletter
  \providecommand{\limlike}[1]{\mathop {\operator@font #1}}
  \providecommand{\loglike}[1]{\mathop {\operator@font #1}\nolimits}
\makeatother
\providecommand{\cosp}{\loglike{cosp}}
\providecommand{\anti}{\mathcal{A}}
\providecommand{\ub}[3][]{\left\{\!#1\left[#2,#3\right]\!#1\right\}}
\providecommand{\sinp}{\loglike{sinp}}

\providecommand{\modulus}[2][\relax]{\left| #2 \right|\ifx#1\relax\else_{#1}\fi}
\hypersetup{colorlinks=true,bookmarks=true}
\newcommand{\wiki}[3][\relax]{#3\ifx#1\relax\index{#3}\else\index{#1}\fi}

\providecommand{\oper}[1]{\mathcal{#1}}

\begin{document}

\title
{Erlangen Programme at Large: An Overview}

\author[Vladimir V. Kisil]%
{\href{http://maths.leeds.ac.uk/~kisilv/}{Vladimir V. Kisil}}
\address{
School of Mathematics,
University of Leeds,
Leeds, LS2\,9JT,
UK
}
\thanks{On  leave from the Odessa University.}
\email{
\href{mailto:kisilv@maths.leeds.ac.uk}{kisilv@maths.leeds.ac.uk}
}
\urladdr{
\url{http://www.maths.leeds.ac.uk/~kisilv/}}


\dedicatory{Dedicated to Prof. Hans G. Feichtinger on the occasion of his 60th birthday}
\begin{abstract}
  This is an overview of \emph{Erlangen Programme at
    Large}. Study of objects and properties, which are invariant under
  a group action, is very fruitful far beyond the traditional
  geometry. In this paper we demonstrate this on the example of the
  group \(\SL\). Starting from the conformal geometry we develop
  analytic functions and apply these to functional calculus. Finally we
  link this to quantum mechanics and conclude by a list of open problems.
\end{abstract}
\maketitle
\tableofcontents

\ifundefined{chapter}
\else
\chapter[Erlangen Programme at Large: An Overview]{Erlangen Programme at Large:\\ An Overview}

\label{cha:erlangen-program-at}
\author{Vladimir V. Kisil}
\address{School of Mathematics,
    University of Leeds,
    Leeds LS2\,9JT, England\\
    E-mail: \href{mailto:kisilv@maths.leeds.ac.uk}{\texttt{kisilv@maths.leeds.ac.uk}}}

\centerline{\emph{Dedicated to Prof. Hans G. Feichtinger on the
    occasion of his 60th birthday}\vspace{5mm}}

{\footnotesize {\scshape Abstract.}
  This is an overview of \emph{Erlangen Programme at
    Large}. Study of objects and properties, which are invariant under
  a group action, is very fruitful far beyond the traditional
  geometry. In this paper we demonstrate this on the example of the
  group \(\SL\). Starting from the conformal geometry we develop
  analytic functions and apply these to functional calculus. Finally we
  link this to quantum mechanics and conclude by a list of open problems.

\vspace*{3mm} {\scshape 2010 Mathematics Subject Classification:}
Primary 30G35; Secondary 22E46, 30F45, 32F45, 43A85, 30G30, 42C40, 46H30, 47A13, 81R30, 81R60.

\vspace*{1mm} {\scshape Key words and phrases:}
 Special linear group, Hardy space, Clifford algebra, elliptic,
  parabolic, hyperbolic, complex numbers, dual
  numbers, double numbers, split-complex numbers,
  Cauchy-Riemann-Dirac operator, M\"obius transformations, functional
  calculus, spectrum, quantum mechanics, non-commutative geometry.}\vspace{5mm}
\fi

\epigraph{A mathematical idea should not be petrified in a formalised
axiomatic settings, but should be considered instead as flowing as a
river.}{}{Sylvester (1878)} 

\section{Introduction}
\label{sec:introduction}

The simplest objects with non-commutative (but still associative)
multiplication may be \(2\times 2\) matrices with real entries.  The
subset of matrices \emph{of determinant one} has the following
properties: 
\begin{itemize}
\item form a closed set under multiplication since
  \(\det (AB)=\det A\cdot \det B\);
\item the identity matrix is the set; and
\item any such matrix has an inverse (since \(\det A\neq 0\)).
\end{itemize}
In
other words those matrices form a \emph{group}\index{group}, the 
\wiki[]{SL2(R)}{$\SL$ group}%
\index{$\SL$ group}%
\index{group!$\SL$}~\cite{Lang85}---one of the two most important Lie groups
in analysis. The other group is 
the \wiki{Heisenberg_group}{Heisenberg group}%
\index{Heisenberg!group}%
\index{group!Heisenberg}~\cite{Howe80a}. By contrast the
\wiki{Affine_transformation}{$ax+b$ group}%
\index{$ax+b$ group}%
\index{group!$ax+b$}%
\index{group!affine|see{$ax+b$ group}}%
\index{affine group|see{$ax+b$ group}}, which is often
used to build wavelets, is only a subgroup of \(\SL\), see the
numerator in~\eqref{eq:moebius}.

The simplest non-linear transforms of the real
line---\emph{linear-fractional}%
\index{M\"obius map}%
\index{map!M\"obius} or \emph{M\"obius maps}
---may
also be associated with \(2\times 2\)
matrices~\amscite{Beardon05a}*{Ch.~13}:
\begin{equation}
  \label{eq:moebius}
  g: x\mapsto g\cdot x=\frac{ax+b}{cx+d}, \quad\text{ where } 
  g=  \begin{pmatrix}
    a&b\\c&d
  \end{pmatrix},\  x\in\Space{R}{}.
\end{equation}
An enjoyable calculation shows that the composition of two
transforms~\eqref{eq:moebius} with different matrices \(g_1\) and
\(g_2\) is again a M\"obius transform with matrix the product
\(g_1 g_2\). In other words~\eqref{eq:moebius} it is a (left) action
of \(\SL\).

According to F.~Klein's \emph{\wiki{Erlangen_program}{Erlangen
    programme}}%
\index{Erlangen programme|indef} (which was influenced by S.~Lie) any
geometry is dealing with invariant properties under a certain
transitive%
\index{transitive}%
\index{action!transitive} group action. For example, we may ask:
\emph{What kinds of geometry are related to the \(\SL\)
  action~\eqref{eq:moebius}}?
  
The Erlangen programme has probably the highest rate of
\(\frac{\text{praised}}{\text{actually used}}\) among mathematical
theories not only due to the big numerator but also due to undeserving
small denominator. As we shall see below Klein's approach provides
some surprising conclusions even for such over-studied objects as
circles.

\subsection{Make a Guess in Three Attempts}
\label{sec:make-guess-three}

It is easy to see that the \(\SL\) action~\eqref{eq:moebius} makes
sense also as a map of complex numbers \(z=x+\rmi y\), \(\rmi^2=-1\)%
\index{$\rmi$ (imaginary unit)|indef}%
\index{imaginary!unit ($\rmi$)|indef}%
\index{unit!imaginary ($\rmi$)|indef} assuming the denominator is
non-zero.  Moreover, if \(y>0\) then \(g\cdot z\) has a positive
imaginary part as well, i.e. \eqref{eq:moebius} defines a map from
the upper half-plane%
\index{half-plane} to itself. Those transformations are isometries
of the Lobachevsky half-plane.%
\index{Lobachevsky!geometry}%
\index{geometry!Lobachevsky}

However there is no need to be restricted to the traditional route of
complex numbers only. Moreover in
  Subsection~\ref{sec:hypercomplex-numbers} we will naturally come to
  a necessity to work with all three kinds of hypercomplex numbers%
\index{number!hypercomplex|indef}%
\index{hypercomplex!number|indef}. Less-known
\emph{\wiki{Split-complex_number}{double}}%
\index{number!double}%
\index{double!number} and
\emph{\wiki{Dual_number}{dual}}%
\index{dual!number}%
\index{number!dual} numbers,
see~\amscite{Yaglom79}*{Suppl.~C}, have also the form \(z=x+\alli y\)
but different assumptions on the hypercomplex unit \(\alli\)%
\index{$\rma$@$\alli$ (hypercomplex unit)}: \(\alli^2=0\) or
\(\alli^2=1\) correspondingly. We will write \(\rmp\)%
\index{$\rmp$ (parabolic unit)|indef}%
\index{parabolic!unit ($\rmp$)|indef}%
\index{unit!parabolic ($\rmp$)|indef}%
\index{unit!nilpotent|see{parabolic unit}}%
\index{nilpotent unit|see{parabolic unit}} and \(\rmh\)%
\index{$\rmh$ (hyperbolic unit)|indef}%
\index{hyperbolic!unit ($\rmh$)|indef}%
\index{unit!hyperbolic ($\rmh$)|indef}
instead of \(\alli\) within dual and double numbers respectively.
Although the arithmetic of dual and double numbers is different from
the complex ones, e.g. they have divisors of zero%
\index{divisor!zero}%
\index{zero!divisor}, we are still able
to define their transforms by~\eqref{eq:moebius} in most cases.

Three possible values \(-1\), \(0\) and \(1\) of \(\sigma:=\alli^2\)%
\index{$\sigma$ ($\sigma:=\alli^2$)} will be refereed to here as
\emph{elliptic}, \emph{parabolic} and \emph{hyperbolic}%
\index{elliptic!case}%
\index{hyperbolic!case}%
\index{parabolic!case}%
\index{case!elliptic}%
\index{case!parabolic}%
\index{case!hyperbolic} cases respectively.  We repeatedly meet such a
division of various mathematical objects into three classes.  They are
named by the historically first example---the classification of conic
sections%
\index{conic!section}%
\index{sections!conic}---however the pattern persistently reproduces itself in many
different areas: equations, quadratic forms, metrics, manifolds,
operators, etc.  We will abbreviate this separation as
\emph{EPH-classification}.%
\index{EPH-classification}  The {common origin} of this fundamental
division of any family with one-parameter can be seen from the simple
picture of a coordinate line split by zero into negative and positive
half-axes:
\begin{equation}
 \label{eq:eph-class}
  \raisebox{-15pt}{\includegraphics[scale=1]{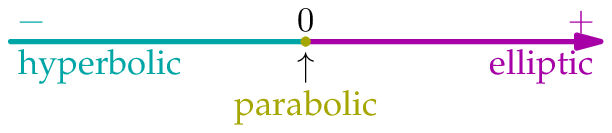}}
\end{equation}

Connections between different objects admitting EPH-classification are
not limited to this common source. There are many deep results
linking, for example, the {ellipticity} of quadratic forms, metrics
and operators, e.g. the \wiki {Atiyah-Singer_index_theorem}%
{Atiyah-Singer index theorem}%
\index{Atiyah-Singer index theorem}%
\index{theorem!Atiyah-Singer (index)}.  On the other hand there are still a
lot of white spots, empty cells, obscure gaps and missing connections
between some subjects as well.

To understand the action~\eqref{eq:moebius} in all EPH cases we
use the Iwasawa decomposition%
\index{Iwasawa decomposition}%
\index{decomposition!Iwasawa}~\amscite{Lang85}*{\S~III.1} of
\(\SL=ANK\) into three one-dimensional subgroups \(A\), \(N\)
and \(K\):%
  \index{$A$ subgroup}%
  \index{subgroup!$A$}%
  \index{$N$ subgroup}%
  \index{subgroup!$N$}%
  \index{$K$ subgroup}%
  \index{subgroup!$K$}
\begin{equation}
  \label{eq:iwasawa-decomp}
  \begin{pmatrix}
    a&b \\c &d
  \end{pmatrix}= {\begin{pmatrix} \alpha & 0\\0&\alpha^{-1}
    \end{pmatrix}} {\begin{pmatrix} 1&\nu \\0&1
    \end{pmatrix}} {\begin{pmatrix}
      \cos\phi &  -\sin\phi\\
      \sin\phi & \cos\phi
    \end{pmatrix}}.
\end{equation}

Subgroups \(A\) and \(N\) act in~\eqref{eq:moebius} irrespectively
to value of \(\sigma\): \(A\)%
\index{$A$-orbit}%
\index{orbit!subgroup $A$, of}%
\index{subgroup!$A$!orbit} makes a dilation by \(\alpha^2\), i.e.
\(z\mapsto \alpha^2z\), and \(N\)%
\index{$N$-orbit}%
\index{orbit!subgroup $N$, of}%
\index{subgroup!$N$!orbit} shifts points to left by \(\nu\), i.e.
\(z\mapsto z+\nu\).

\begin{figure}[htbp]
  \centering
  \includegraphics[scale=.8]{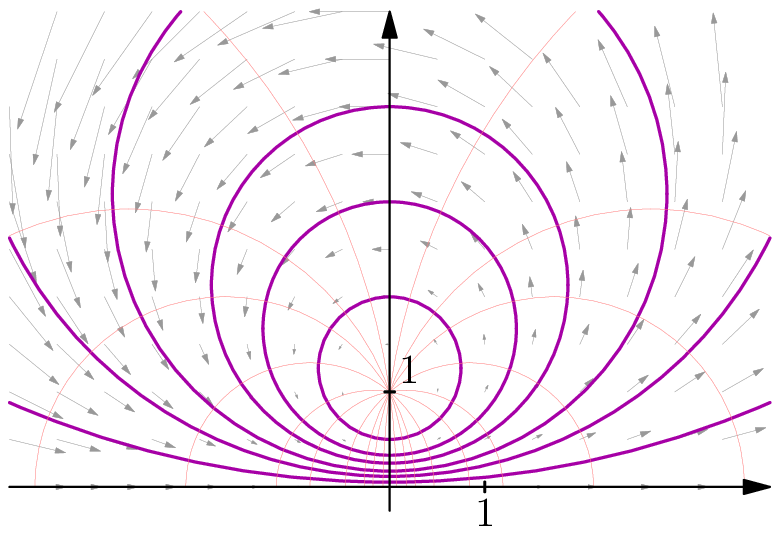} 
  \hfill
  \includegraphics[scale=.8]{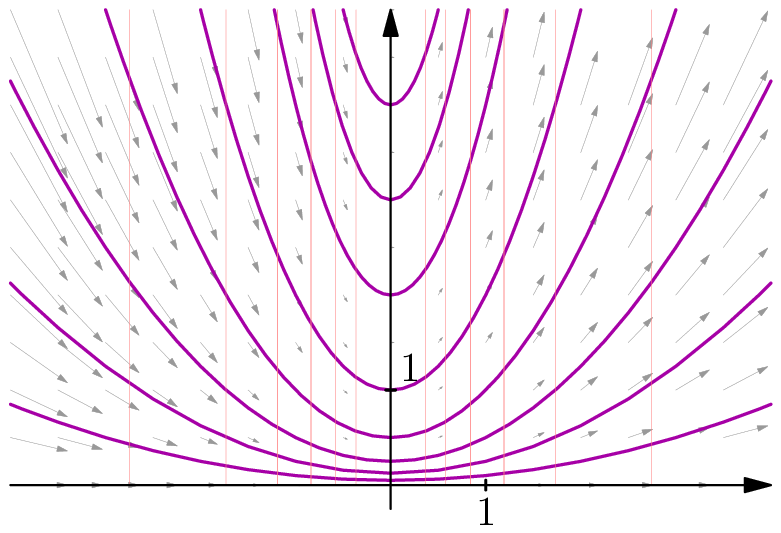}\\  
  \includegraphics[scale=.8]{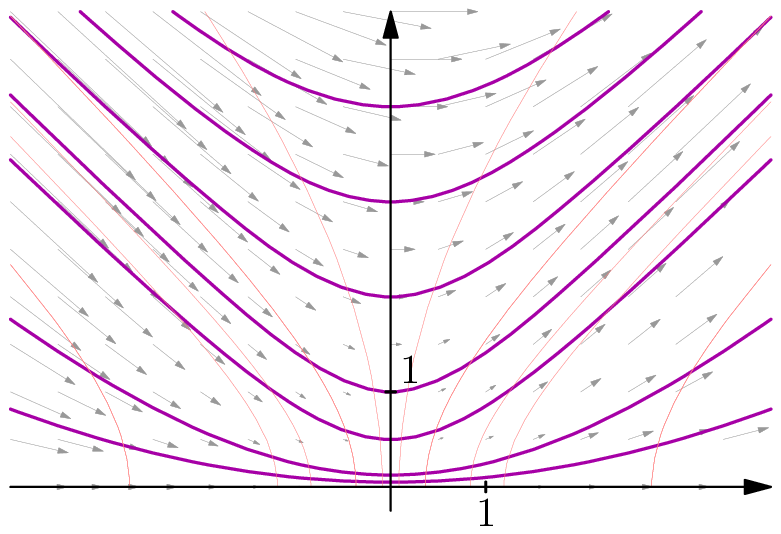}\hfill
  \parbox[b]{.45\textwidth}{ The corresponding orbits are circles,
    parabolas and hyperbolas shown by thick lines.  Transverse thin
    lines are images of the vertical axis under the action of the
    subgroup \(K\). Grey arrows show the associated
    derived action.\\}
  \caption[Action of the subgroup $K$]{Action of the subgroup \(K\).
}
  \label{fig:k-subgroup}
\end{figure}

By contrast, the action of the third matrix from the subgroup \(K\)%
\index{$K$-orbit}%
\index{orbit!subgroup $K$, of}%
\index{subgroup!$K$!orbit} 
sharply depends on \(\sigma\), see
Fig.~\ref{fig:k-subgroup}. In elliptic, parabolic and hyperbolic
cases \(K\)-orbits are circles, parabolas and (equilateral) hyperbolas
correspondingly.  Thin traversal lines in
Fig.~\ref{fig:k-subgroup} join points of orbits for the same
values of \(\phi\) and grey arrows represent ``local
velocities''---vector fields of derived representations.%
\index{derived action}%
\index{action!derived} We will describe some highlights of this
geometry in Section~\ref{sec:geometry}.

\subsection{Erlangen Programme at Large}
\label{sec:erlangen-program-at}
As we already mentioned the division of mathematics into areas is only
apparent. Therefore it is unnatural to limit Erlangen programme only to
``geometry''. We may continue to look for \(\SL\) invariant objects in
other related fields. For example, transform~\eqref{eq:moebius}
generates unitary representations%
\index{group!representation!linear}%
\index{representation!linear} on certain \(\FSpace{L}{2}\) spaces,
cf.~\eqref{eq:moebius} and Section~\ref{sec:induc-repr}:
\begin{equation}
  \label{eq:hardy-repres}
  g: f(x)\mapsto \frac{1}{(cx+d)^m}f\left(\frac{ax+b}{cx+d}\right).
\end{equation}

For \(m=1\), \(2\), \ldots\ the invariant subspaces of
\(\FSpace{L}{2}\) are Hardy%
\index{space!Hardy}%
\index{Hardy!space} and (weighted) Bergman spaces%
\index{space!Bergman}%
\index{Bergman!space} of complex
analytic functions.  All main objects of \emph{complex analysis}
(Cauchy%
\index{Cauchy!integral}%
\index{integral!Cauchy} and Bergman integrals, Cauchy-Riemann%
\index{Cauchy-Riemann operator}%
\index{operator!Cauchy-Riemann}%
\index{equation!Cauchy-Riemann|see{Cauchy-Riemann operator}}
 and Laplace equations%
\index{Laplacian},
Taylor series etc.) may be obtaining in terms of invariants of the
\emph{discrete series}%
\index{representation!discrete series} representations of
\(\SL\)~\amscite{Kisil02c}*{\S~3}.  Moreover two other series
(\emph{principal}%
\index{representation!principal series} and \emph{complimentary}%
\index{representation!complementary series}~\cite{Lang85}) play the
similar r\^oles for hyperbolic and parabolic
cases~\cite{Kisil02c,Kisil05a}. This will be discussed in
Sections~\ref{sec:covariant-transform} and~\ref{sec:analytic-functions}. 

Moving further we may observe that transform~\eqref{eq:moebius} is
defined also for an element \(x\) in any algebra \(\algebra{A}\) with
a unit \(\mathbf{1}\) as soon as \((cx+d\mathbf{1})\in\algebra{A}\)
has an inverse. If \(\algebra{A}\) is equipped with a topology, e.g.
is a Banach algebra, then we may study a \emph{functional calculus}%
\index{functional!calculus}%
\index{calculus!functional}
for element \(x\)~\cite{Kisil02a} in this way. It is defined as an
\emph{intertwining operator}%
\index{intertwining operator}%
\index{operator!intertwining}  between the
representation~\eqref{eq:hardy-repres} in a space of analytic
functions and a similar representation in a left
\(\algebra{A}\)-module. We will consider the
Section~\ref{sec:functional-calculus}. 

In the spirit of Erlangen programme such functional calculus is still a
geometry, since it is dealing with invariant properties under a group
action. However even for a simplest non-normal operator, e.g.  a
Jordan block of the length \(k\), the obtained space is not like a
space of point but is rather a space of \(k\)-th
\emph{jets}\index{jet}~\cite{Kisil02a}. Such non-point behaviour is oftenly
attributed to \emph{non-commutative geometry}%
\index{non-commutative geometry}%
\index{geometry!non-commutative} and Erlangen programme
provides an important input on this fashionable topic~\cite{Kisil02c}.

It is noteworthy that ideas of F.~Klein ans S.~Lie are spread more
in physics than in mathematics: it is a common viewpoint that laws of
nature shall be invariant under certain transformations. Yet
systematic use of Erlangen approach can bring new results even in this
domain as we demonstrate in
Section~\ref{sec:quantum-mechanics-1}. There are still many directions
to extend the present work thus we will conclude by a list of some
open problems in Section~\ref{sec:open-problems}.

Of course, there is no reasons to limit Erlangen programme to \(\SL\)
group only, other groups may be more suitable in different situations.
However \(\SL\) still possesses a big unexplored potential and is a
good object to start with.

\section{Geometry}
\label{sec:geometry}
We start from the natural domain of the Erlangen
Programme---geometry. Systematic use of this ideology allows to obtain
new results even for very classical objects like circles.

\subsection{Hypercomplex Numbers}
\label{sec:hypercomplex-numbers}
Firstly we wish to demonstrate that hypercomplex numbers appear very
naturally from a study of \(\SL\) action on the homogeneous
spaces~\cite{Kisil09c}. We begin from the standard definitions.

Let \(H\) be a subgroup of a group \(G\).  Let \( X=G / H\) be
the corresponding homogeneous space and \(s: X \rightarrow G\) be
a smooth section~\amscite{Kirillov76}*{\S~13.2}, which is a left inverse
to the natural projection \(p: G\rightarrow X \). The choice of
\(s\) is inessential in the sense that by a smooth map
\(X\rightarrow X\) we can always reduce one to another.  We
define a map \(r: G\rightarrow H\) associated to \(p\) and \(s\) from
the identities:
\begin{equation}
  \label{eq:r-map}
  r(g)={(s(x))}^{-1}g, \qquad \text{where }
  x=p(g)\in X .
\end{equation}
Note that \(X \) is a left homogeneous space with the
\(G\)-action defined in terms of \(p\) and \(s\) as follows:
\begin{equation}
  \label{eq:g-action}
  g: x  \mapsto g\cdot x=p(g* s(x)),
\end{equation}
\begin{example}[\cite{Kisil09c}]
  \label{ex:sl2-hypernumber}
  For \(G=\SL\), as well as for other semisimple groups, it is common
  to consider only the case of \(H\) being the maximal compact
  subgroup \(K\). However in this paper we admit \(H\) to be any
  one-dimensional subgroup.  Then \(X\) is a two-dimensional
  manifold and for any choice of \(H\) we
  define~\amscite{Kisil97c}*{Ex.~3.7(a)}:
  \begin{equation}
    \label{eq:s-map}
    s: (u,v) \mapsto
    \frac{1}{\sqrt{v}}
    \begin{pmatrix}
      v & u \\ 0 & 1
    \end{pmatrix}, \qquad (u,v)\in\Space{R}{2},\  v>0.
  \end{equation}
  Any continuous one-dimensional subgroup \(H\in\SL\) is conjugated to one
  of the following:
  \begin{eqnarray}
    \label{eq:k-subgroup}
    K&=&\left\{ {\begin{pmatrix}
          \cos t &  \sin t\\
          -\sin t & \cos t
        \end{pmatrix}=   \exp \begin{pmatrix} 0& t\\-t&0
        \end{pmatrix}},\ t\in(-\pi,\pi]\right\},\\
    \label{eq:np-subgroup}
    N'&=&\left\{   {\begin{pmatrix} 1&0\\t&1
        \end{pmatrix}=\exp \begin{pmatrix} 0 & 0\\t&0
        \end{pmatrix},}\  t\in\Space{R}{}\right\},\\
    \label{eq:ap-subgroup}
    \Aprime&=&\left\{  
      \begin{pmatrix} \cosh t & \sinh t\\ \sinh t& \cosh t
      \end{pmatrix}=\exp \begin{pmatrix} 0 & t\\t&0
      \end{pmatrix},\  t\in\Space{R}{}\right\}.
  \end{eqnarray}
  \index{subgroup!$A$}
  \index{subgroup!$N$}
  \index{subgroup!$K$}
  Then~\cite{Kisil09c} the action~\eqref{eq:g-action} of \(\SL\) on
  \(X=\SL/H\) coincides with M\"obius transformations~\eqref{eq:moebius}
  on complex, dual and double numbers respectively.
\end{example}

\subsection{Cycles as Invariant Families}
\label{sec:cycles-as-invariant}
We wish to consider all three hypercomplex systems at the same time,
the following definition is very helpful for this.

\begin{defn}
  \label{de:cycle}
  The common name \emph{cycle}\index{cycle}~\cite{Yaglom79} is used to
  denote circles, parabolas and hyperbolas (as well as straight lines
  as their limits) in the respective EPH case.
\end{defn}

\begin{figure}[htbp]
  \centering
  (a)\includegraphics[scale=1]{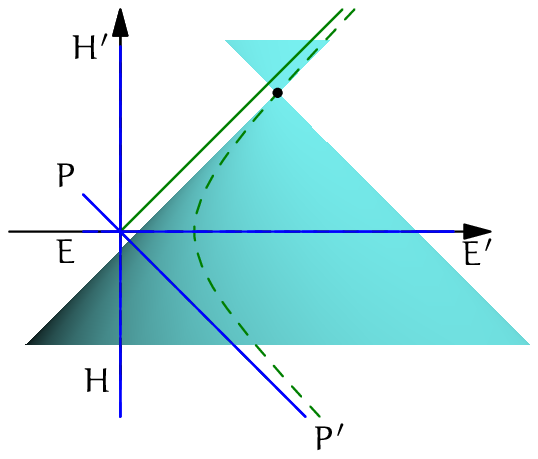}\hfill
  (b)\includegraphics[scale=1]{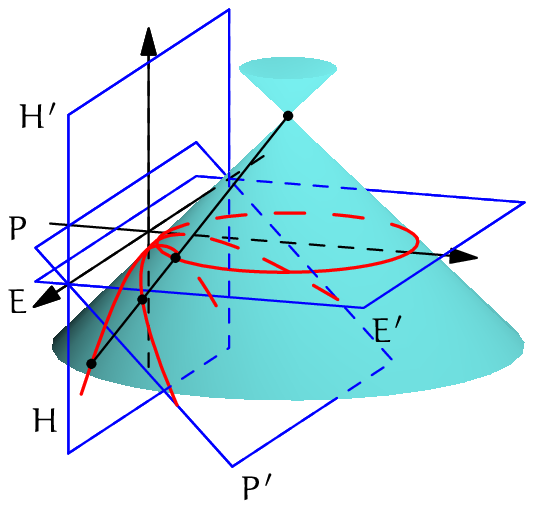}
  \caption[$K$-orbits as conic sections]{\(K\)-orbits as conic
    sections: circles are sections by the plane \(EE'\); parabolas are
    sections by \(PP'\); hyperbolas are sections by \(HH'\). Points on
    the same generator of the cone correspond to the same value of
    \(\phi\).}
  \label{fig:k-orbit-sect}
\end{figure}
It is well known that any cycle is a \emph{conic sections}%
\index{conic!section}%
\index{sections!conic} and an
interesting observation is that corresponding \(K\)-orbits are in fact
sections of the same two-sided right-angle cone, see
Fig.~\ref{fig:k-orbit-sect}.  Moreover, each straight line generating
the cone, see Fig.~\ref{fig:k-orbit-sect}(b), is crossing
corresponding EPH \(K\)-orbits at points with the same value of
parameter \(\phi\) from~\eqref{eq:iwasawa-decomp}. In other words,
all three types of orbits are generated by the rotations of this
generator%
\index{generator!of a cone}%
\index{cone!generator} along the cone.

\(K\)-orbits are \(K\)-invariant in a trivial way. Moreover since
actions of both \(A\) and \(N\) for any \(\sigma\) are extremely
``shape-preserving'' we find natural invariant objects of the M\"obius
map:

\begin{thm}
  \label{th:invariance-of-cycles}
  The family of all cycles from Definition~\ref{de:cycle} is invariant under the
  action~\eqref{eq:moebius}.
\end{thm}

According to Erlangen ideology we should now study invariant
properties of cycles. 

Fig.~\ref{fig:k-orbit-sect} suggests that we may get a unified
treatment of cycles in all EPH cases by consideration of a higher
dimension spaces. The standard mathematical method is to declare
objects under investigations (cycles in our case, functions in
functional analysis, etc.) to be simply points of some bigger space.
This space should be equipped with an appropriate structure to hold
externally information which were previously inner properties of our
objects.
  
A generic cycle is the set of points \((u,v)\in\Space{R}{2}\) defined
for all values of \(\sigma\) by the equation\index{cycle!equation}
\begin{equation}
  \label{eq:cycle-eq}
  k(u^2-\sigma v^2)-2lu-2nv+m=0.
\end{equation}
This equation (and the corresponding cycle) is defined by a point
\((k, l, n, m)\) from a \emph{projective space}%
\index{projective space}%
\index{space!projective} \(\Space{P}{3}\), since for a scaling factor
\(\lambda \neq 0\) the point \((\lambda k, \lambda l, \lambda n,
\lambda m)\) defines an equation equivalent
to~\eqref{eq:cycle-eq}. We call \(\Space{P}{3}\) the \emph{cycle
  space}%
\index{cycle!space}%
\index{space!cycles, of} and refer to the initial \(\Space{R}{2}\) as
the \emph{point space}%
\index{point space}%
\index{space!point}.

In order to get a connection with M\"obius action~\eqref{eq:moebius}
we arrange numbers \((k, l, n, m)\) into the matrix%
\index{cycle!matrix}%
\index{matrix!cycle, of a}  
\begin{equation}
  \label{eq:FSCc-matrix}
  C_{\bs}^s=\begin{pmatrix}
    l+\rmc s n&-m\\k&-l+\rmc s n
  \end{pmatrix}, 
\end{equation}
with a new hypercomplex unit \(\rmc\)%
\index{$\rmc$ (hypercomplex unit in cycle space)} and an additional
parameter \(s\) usually equal to \(\pm 1\). The values of
\(\bs:=\rmc^2\)%
\index{$\sigma1$@$\bs$ ($\bs:=\rmc^2$)} is \(-1\), \(0\) or \(1\)
independently from the value of \(\sigma\).  The
matrix~\eqref{eq:FSCc-matrix} is the cornerstone of an extended
\emph{Fillmore--Springer--Cnops construction}%
\index{Fillmore--Springer--Cnops construction}%
\index{construction!Fillmore--Springer--Cnops} (FSCc)~\cite{Cnops02a}.

The significance of FSCc in Erlangen framework is provided by the
following result:
\begin{thm}
  \label{th:FSCc-intertwine}
  The image  \(\tilde{C}_{\bs}^s\) of a cycle \(C_{\bs}^s\) under
  transformation~\eqref{eq:moebius} with \(g\in\SL\) is given by
  similarity of the matrix~\eqref{eq:FSCc-matrix}:
  \begin{equation}
    \label{eq:cycle-similarity}
    \tilde{C}_{\bs}^s= gC_{\bs}^sg^{-1}.
  \end{equation}
  In other words FSCc~\eqref{eq:FSCc-matrix} \emph{intertwines}
  M\"obius action~\eqref{eq:moebius}%
    \index{M\"obius map!on cycles}%
    \index{map!M\"obius!on cycles} on cycles with
  linear map~\eqref{eq:cycle-similarity}.
\end{thm}
  
There are several ways to prove~\eqref{eq:cycle-similarity}: either
by a brute force calculation (fortunately
\href{http://arxiv.org/abs/cs.MS/0512073}{performed by a
  CAS})~\cite{Kisil05a} or through the related orthogonality of
cycles~\cite{Cnops02a}%
  \index{orthogonality!cycles, of}%
  \index{cycles!orthogonal}, see the end of the next
Subsection~\ref{sec:invar-algebr-geom}. 

The important observation here is that our extended version of
FSCc~\eqref{eq:FSCc-matrix} uses an hypercomplex unit \(\rmc\), which
is not related to \(\alli\) defining the appearance of cycles on plane.
In other words any EPH type of geometry in the cycle space
\(\Space{P}{3}\) admits drawing of cycles in the point space
\(\Space{R}{2}\) as circles, parabolas or hyperbolas. We may think on
points of \(\Space{P}{3}\) as ideal cycles while their depictions on
\(\Space{R}{2}\) are only their shadows on the wall of
\wiki{Plato\#Metaphysics}%
{Plato's cave}%
\index{Plato's cave}.

\begin{figure}[htbp]
  \centering
  (a)\includegraphics[]{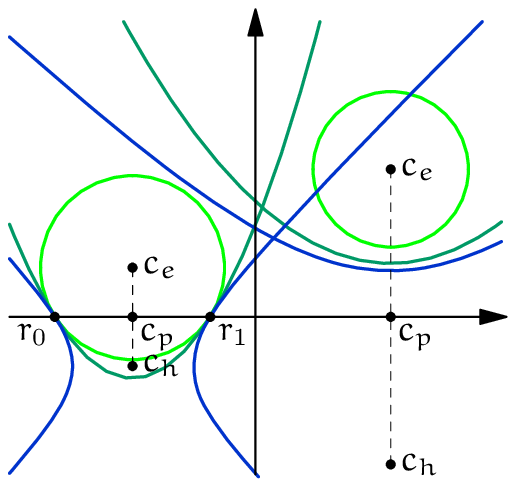}\hfill
  (b)\includegraphics[]{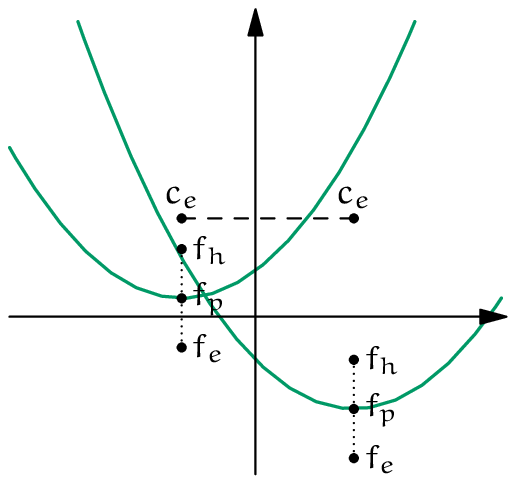}
  \caption[Cycle implementations, centres and foci]{
    (a) Different
    EPH implementations of the same cycles defined by quadruples of
    numbers.\\
    (b) Centres and foci of two parabolas with the same focal length.} 
  \label{fig:eph-cycle}
\end{figure}
Fig.~\ref{fig:eph-cycle}(a) shows the same cycles drawn in
different EPH styles. Points \(c_{e,p,h}=(\frac{l}{k}, -\bs
\frac{n}{k})\) are their respective e/p/h-centres%
\index{centre!of a cycle}%
\index{cycle!centre}. They are related to each other through several
identities:
\begin{equation}
  \label{eq:centres}
  c_e=\bar{c}_h, \quad c_p=\frac{1}{2}(c_e+c_h).
\end{equation}
Fig.~\ref{fig:eph-cycle}(b) presents two cycles drawn as
parabolas, they have the same focal length \(\frac{n}{2k}\) and thus
their e-centres are on the same level. In other words
\emph{concentric}%
\index{concentric} parabolas are obtained by a vertical shift, not
scaling as an analogy with circles or hyperbolas may suggest.

Fig.~\ref{fig:eph-cycle}(b) also presents points, called
e/p/h-foci%
\index{focus!of a cycle}%
\index{cycle!focus}:
\begin{equation}
  \label{eq:foci}
  f_{e,p,h}=\left(\frac{l}{k}, -\frac{\det C_{\bs}^s}{2nk}\right),
\end{equation}
which are independent of the sign of \(s\).  If a cycle is depicted as
a parabola then h-focus, p-focus, e-focus are correspondingly
geometrical focus%
\index{focus!of a parabola}%
\index{parabola!focus} of the parabola, its vertex%
\index{vertex!of a parabola}%
\index{parabola!vertex}, and the point on the directrix%
\index{directrix}%
\index{parabola!directrix} nearest to the vertex.

As we will see, cf. Theorems~\ref{th:ghost1} and~\ref{th:ghost2},
all three centres and three foci are useful attributes of a cycle even
if it is drawn as a circle.

\subsection{Invariants: Algebraic and Geometric}
\label{sec:invar-algebr-geom}

We use known algebraic invariants of matrices to build appropriate
geometric invariants of cycles. It is yet another demonstration that
any division of mathematics into subjects is only illusive.

For \(2\times 2\) matrices (and thus cycles) there are only two
essentially different invariants under
similarity~\eqref{eq:cycle-similarity} (and thus under M\"obius
action~\eqref{eq:moebius}): the \emph{trace}\index{trace} and the
\emph{determinant}\index{determinant}.  The latter was already used
in~\eqref{eq:foci} to define cycle's foci. However due to
projective nature of the cycle space \(\Space{P}{3}\) the absolute
values of trace or determinant are irrelevant, unless they are zero.

Alternatively we may have a special arrangement for normalisation%
\index{normalised!cycle}%
\index{cycle!normalised} of quadruples \((k,l,n,m)\). For example,
if \(k\neq0\) we may normalise the quadruple to
\((1,\frac{l}{k},\frac{n}{k},\frac{m}{k})\) with highlighted cycle's
centre. Moreover in this case \(\det \cycle{s}{\bs}\) is equal to the
square of cycle's radius, cf. Section~\ref{sec:dist-lenght-perp}.
Another normalisation \(\det \cycle{s}{\bs}=1\) is used
in~\cite{Kirillov06}%
\index{normalised!cycle!Kirillov}%
\index{cycle!normalised!Kirillov} to get a nice condition for touching
circles. Moreover, the Kirillov normalisation is preserved by the
conjugation~\eqref{eq:cycle-similarity}. 

We still get important characterisation even with non-normalised
cycles, e.g., invariant classes (for different \(\bs\)) of cycles are
defined by the condition \(\det C_{\bs}^s=0\). Such a class is
parametrises only by two real numbers and as such is easily attached to
certain point of \(\Space{R}{2}\). For example, the cycle
\(C_{\bs}^s\) with \(\det C_{\bs}^s=0\), \(\bs=-1\) drawn elliptically
represent just a point \((\frac{l}{k},\frac{n}{k})\), i.e. (elliptic)
zero-radius circle%
\index{zero-radius cycle}%
\index{cycle!zero-radius}.  The same condition with \(\bs=1\) in
hyperbolic drawing produces a null-cone originated at point
\((\frac{l}{k},\frac{n}{k})\):
\begin{displaymath}
  \textstyle (u-\frac{l}{k})^2-(v-\frac{n}{k})^2=0,
\end{displaymath}
i.e. a zero-radius cycle in hyperbolic metric. 

\begin{figure}[htbp]
  \centering
  \includegraphics[scale=.71]{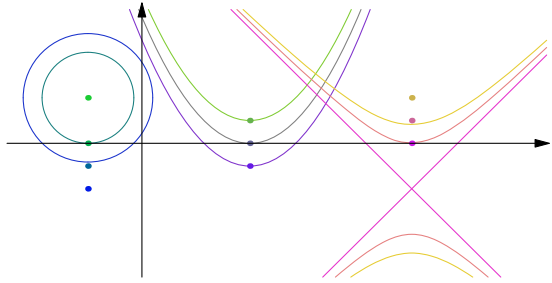}
  \caption[Different implementations of the same
  zero-radius cycles]{Different \(\sigma\)-implementations of the same
    \(\bs\)-zero-radius cycles and corresponding foci.}
  \label{fig:zero-radius}
\end{figure}
In general for every notion there are (at least) nine possibilities:
three EPH cases in the cycle space times three EPH realisations in the
point space. Such nine cases for ``zero radius'' cycles is shown in
Fig.~\ref{fig:zero-radius}. For example, p-zero-radius cycles in any
implementation touch the real axis.

This ``touching'' property is a manifestation of the \emph{boundary
  effect}%
\index{half-plane!boundary effect}%
\index{boundary effect on the upper half-plane} in the upper-half
plane geometry. The famous question
\wiki{Hearing_the_shape_of_a_drum}%
{on hearing drum's shape}%
\index{hearing drum's shape}%
\index{drum!hearing shape} has a sister:
\begin{quote}
  \emph{Can we see/feel the boundary from inside a domain?}
\end{quote}
Both orthogonality relations described below are ``boundary aware'' as
well. It is not surprising after all since \(\SL\) action on the
upper-half plane was obtained as an extension of its
action~\eqref{eq:moebius} on the boundary.

According to the \wiki{Category_theory}{categorical viewpoint}%
\index{category theory}
internal properties of objects are of minor importance in comparison
to their relations with other objects from the same class. As an
illustration we may put the proof of Theorem~\ref{th:FSCc-intertwine}
sketched at the end of of the next section. Thus from now on we will
look for invariant relations between two or more cycles.

\subsection{Joint Invariants: Orthogonality}
\label{sec:joint-invar-orth}

The most expected relation between cycles is based on the following
M\"obius invariant ``inner product'' build from a trace of product of
two cycles as matrices:
\begin{equation}
  \label{eq:inner-prod}
  \scalar{C_{\bs}^s}{\tilde{C}_{\bs}^s}= \tr (C_{\bs}^s\tilde{C}_{\bs}^s)
\end{equation}
By the way, an inner product of this type is used, for example, in
\wiki{Gelfand-Naimark-Segal_construction}{GNS construction}%
\index{GNS construction|see{Gelfand--Naimark--Segal construction}}%
\index{Gelfand--Naimark--Segal construction}%
\index{construction!Gelfand--Naimark--Segal} to make a Hilbert space
out of \(C^*\)-algebra.  The next standard move is given by the
following definition.
\begin{defn}
  \label{de:orthogonality}
  Two cycles are called \emph{\(\bs\)-orthogonal}%
  \index{orthogonality!cycles, of}%
  \index{cycles!orthogonal} if
  \(\scalar{C_{\bs}^s}{\tilde{C}_{\bs}^s}=0\).
\end{defn}
The orthogonality relation is preserved under the M\"obius
transformations, thus this is an example of a \emph{joint invariant}%
\index{joint!invariant}%
\index{invariant!joint} of two cycles.  For the case of \(\bs
\sigma=1\), i.e. when geometries of the cycle and point spaces are
both either elliptic or hyperbolic, such an orthogonality is the
standard one, defined in terms of angles between tangent lines in the
intersection points of two cycles.  However in the remaining seven
(\(=9-2\)) cases the innocent-looking
Definition~\ref{de:orthogonality} brings unexpected relations.

\begin{figure}[htbp]
  \includegraphics[scale=.95]{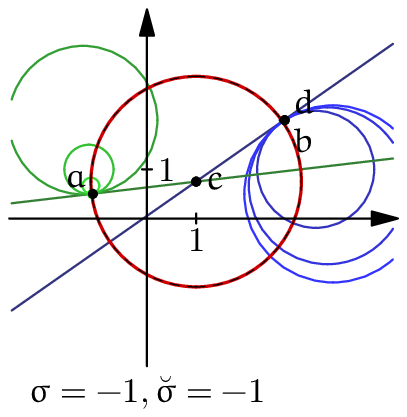}\hfill
  \includegraphics[scale=.95]{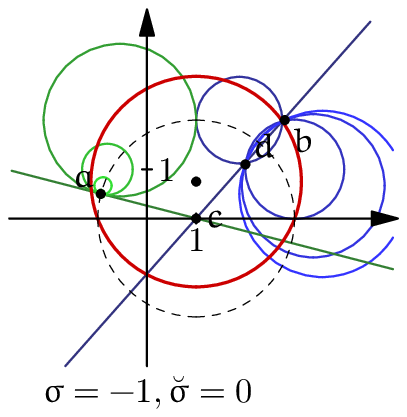}\hfill
  \includegraphics[scale=.95]{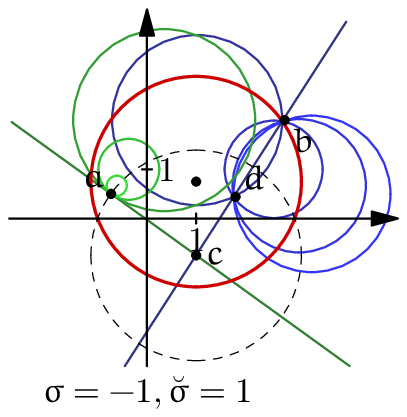}  \caption[Orthogonality of the first kind]{Orthogonality of the first
    kind in the elliptic point space.\\
    Each picture presents two groups (green and blue) of cycles which
    are orthogonal to the red cycle \(C^{s}_{\bs}\).  Point \(b\)
    belongs to \(C^{s}_{\bs}\) and the family of blue cycles
    passing through \(b\) is orthogonal to \(C^{s}_{\bs}\). They
    all also intersect in the point \(d\) which is the inverse of
    \(b\) in \(C^{s}_{\bs}\). Any orthogonality is reduced to the usual
    orthogonality with a new (``ghost'') cycle (shown by the dashed
    line), which may or may not coincide with \(C^{s}_{\bs}\). For
    any point \(a\) on the ``ghost'' cycle the orthogonality is
    reduced to the local notion in the terms of tangent lines at the
    intersection point. Consequently such a point \(a\) is always the
    inverse of itself.}
  \label{fig:orthogonality1}
\end{figure}
Elliptic (in the point space) realisations of
Definition~\ref{de:orthogonality}, i.e. \(\sigma=-1\) is shown in
Fig.~\ref{fig:orthogonality1}. The left picture corresponds to the
elliptic cycle space, e.g. \(\bs=-1\). The orthogonality between the
red circle and any circle from the blue or green families is given in
the usual Euclidean sense. The central (parabolic in the cycle space)
and the right (hyperbolic) pictures show non-local nature of the
orthogonality.  There are analogues pictures in parabolic and
hyperbolic point spaces as well, see~\cites{Kisil05a,Kisil12a}.

This orthogonality may still be expressed in the traditional sense if
we will associate to the red circle the corresponding ``ghost''
circle, which shown by the dashed line in
Fig.~\ref{fig:orthogonality1}.  To describe ghost cycle we need
the \emph{\wiki{Heaviside_step_function}{Heaviside function}}%
\index{function!Heaviside}%
\index{Heaviside!function}
\(\chi(\sigma)\):
\begin{equation}
  \label{eq:heaviside-function}
  \chi(t)=\left\{
    \begin{array}{ll}
      1,& t\geq 0;\\
      -1,& t<0.
    \end{array}\right.
\end{equation}

\begin{thm}
  \label{th:ghost1}
  A cycle is \(\bs\)-orthogonal to cycle \(C_{\bs}^s\) if it is
  orthogonal in the usual sense to the \(\sigma\)-realisation of
  ``ghost'' cycle \(\hat{C}_{\bs}^s\), which is defined by the
  following two conditions:
  \begin{enumerate}
  \item \label{item:centre-centre-rel}
    \(\chi(\sigma)\)-centre of \(\hat{C}_{\bs}^s\) coincides
    with  \(\bs\)-centre of \(C_{\bs}^s\).
  \item Cycles \(\hat{C}_{\bs}^s\) and \(C^{s}_{\bs}\) have the same
    roots, moreover \(\det \hat{C}_{\sigma}^1= \det C^{\chi(\bs)}_{\sigma}\).
  \end{enumerate}
\end{thm}
The above connection between various centres of cycles illustrates
their meaningfulness within our approach.

One can easy check the following orthogonality properties of the
zero-radius cycles defined in the previous section:
\begin{enumerate}
\item Due to the identity \(\scalar{C_{\bs}^s}{{C}_{\bs}^s}=\det
  {C}_{\bs}^s\) zero-radius cycles%
  \index{zero-radius cycle}%
  \index{cycle!zero-radius} are self-ortho\-go\-nal (isotropic)%
  \index{cycle!isotropic}%
  \index{isotropic cycles} ones.
\item \label{it:ortho-incidence} A cycle \(\cycle{s}{\bs}\) is
  \(\sigma\)-orthogonal to a zero-radius cycle \(\zcycle{s}{\bs}\) if
  and only if \(\cycle{s}{\bs}\) passes through the \(\sigma\)-centre
  of \(\zcycle{s}{\bs}\).
\end{enumerate}

As we will see, in parabolic case there is a more suitable notion of
an infinitesimal cycle which can be used instead of zero-radius ones.

\subsection{Higher Order Joint Invariants: f-Orthogonality}
\label{sec:higher-order-joint}

With appetite already wet one may wish to build more joint
invariants.%
\index{joint!invariant}%
\index{invariant!joint}
Indeed for any polynomial \(p(x_1,x_2,\ldots,x_n)\) of
several non-commuting variables one may define an invariant joint
disposition of \(n\) cycles \({}^j\!\cycle{s}{\bs}\) by the condition:
\begin{displaymath}
  \tr p({}^1\!\cycle{s}{\bs}, {}^2\!\cycle{s}{\bs}, \ldots,  {}^n\!\cycle{s}{\bs})=0.
\end{displaymath}
However it is preferable to keep some geometrical meaning of
constructed notions.

An interesting observation is that in the matrix similarity of
cycles~\eqref{eq:cycle-similarity} one may replace element
\(g\in\SL\) by an arbitrary matrix corresponding to another cycle.
More precisely the product%
\index{cycle!conjugation}%
\index{conjugation!cycles, of}
\(\cycle{s}{\bs}\cycle[\tilde]{s}{\bs}\cycle{s}{\bs}\) is again the
matrix of the form~\eqref{eq:FSCc-matrix} and thus may be
associated to a cycle. This cycle may be considered as the reflection%
\index{reflection!in a cycle}%
\index{cycle!reflection} of \(\cycle[\tilde]{s}{\bs}\) in
\(\cycle{s}{\bs}\).
\begin{defn}
  \label{de:f-ortho}
  A cycle \(\cycle{s}{\bs}\) is f-orthogonal%
    \index{f-orthogonality}%
    \index{cycles!f-orthogonal}%
  \index{orthogonality!cycles, of!focal|see{f-orthogonality}} (focal orthogonal%
    \index{focal!orthogonality|see{f-orthogonality}}%
    \index{cycles!f-orthogonal}%
    \index{cycles!focal orthogonal|see{f-orthogonal}}) \emph{to} a cycle
  \(\cycle[\tilde]{s}{\bs}\) if the reflection of
  \(\cycle[\tilde]{s}{\bs}\) in \(\cycle{s}{\bs}\) is orthogonal
  (in the sense of Definition~\ref{de:orthogonality}) to the real line.
  Analytically this is defined by:
  \begin{equation}
    \label{eq:f-orthog-def}
    \tr(\cycle{s}{\bs} \cycle[\tilde]{s}{\bs}\cycle{s}{\bs}\realline{s}{\bs})=0.
  \end{equation}
\end{defn}
Due to invariance of all components in the above definition
f-orthogonality is a M\"obius invariant condition. Clearly this is not
a symmetric relation: if \(\cycle{s}{\bs}\) is f-orthogonal to
\(\cycle[\tilde]{s}{\bs}\) then \(\cycle[\tilde]{s}{\bs}\) is not
necessarily f-orthogonal to \(\cycle{s}{\bs}\).
  
\begin{figure}[htbp]
  \includegraphics[scale=.95]{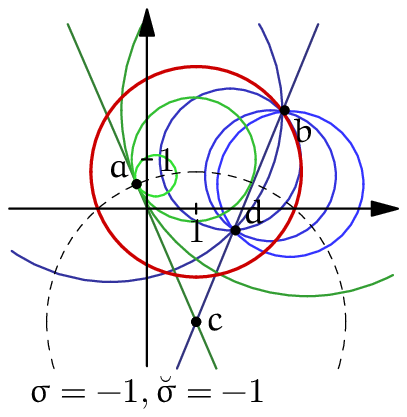}\hfill
  \includegraphics[scale=.95]{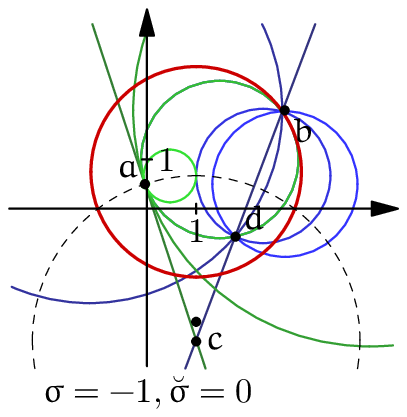}\hfill
  \includegraphics[scale=.95]{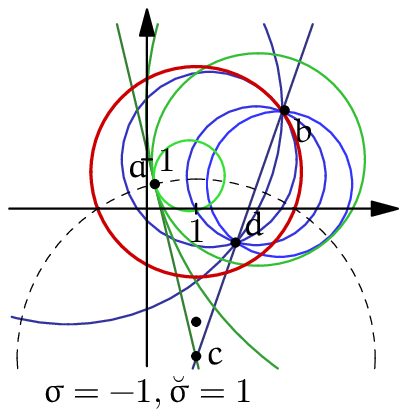}
  \caption[Focal Orthogonality]{Focal orthogonality for circles. To highlight both
    similarities and distinctions with the ordinary orthogonality we
    use the same notations as that in Fig.~\ref{fig:orthogonality1}.}
  \label{fig:orthogonality2}
\end{figure}
Fig.~\ref{fig:orthogonality2} illustrates f-orthogonality in the
elliptic point space. By contrast with
Fig.~\ref{fig:orthogonality1} it is not a local notion at the
intersection points of cycles for all \(\bs\). However it may be again
clarified in terms of the appropriate s-ghost cycle, cf.
Theorem~\ref{th:ghost1}.
\begin{thm}
  \label{th:ghost2}
  A cycle is f-orthogonal to a cycle \(C^{s}_{\bs}\) if its orthogonal
  in the traditional sense to its \emph{f-ghost cycle}
  \(\cycle[\tilde]{\bs}{\bs} = \cycle{\chi(\sigma)}{\bs}
  \Space[\bs]{R}{\bs} \cycle{\chi(\sigma)}{\bs}\), which is the
  reflection of the real line in \(\cycle{\chi(\sigma)}{\bs}\) and
  \(\chi\) is the \emph{Heaviside
    function}~\eqref{eq:heaviside-function}.  Moreover
  \begin{enumerate}
  \item \label{item:focal-centre-rel} \(\chi(\sigma)\)-Centre of
    \(\cycle[\tilde]{\bs}{\bs}\) coincides with the \(\bs\)-focus of
    \(\cycle{s}{\bs}\), consequently all lines f-orthogonal to
    \(\cycle{s}{\bs}\) are passing the respective focus.
  \item Cycles \(\cycle{s}{\bs}\) and \(\cycle[\tilde]{\bs}{\bs}\)
    have the same roots.
  \end{enumerate}
\end{thm}
Note the above intriguing interplay between cycle's centres and foci.
Although f-orthogonality may look exotic it will naturally appear in
the end of next Section again.

Of course, it is possible to define another interesting higher order
joint invariants of two or even more cycles.
  
\subsection{Distance, Length and Perpendicularity}
\label{sec:dist-lenght-perp}
Geo\emph{metry} in the plain meaning of this word deals with
\emph{distances} and \emph{lengths}. Can we obtain them from cycles?
  
\begin{figure}[htbp]
  \centering
  (a) \includegraphics[]{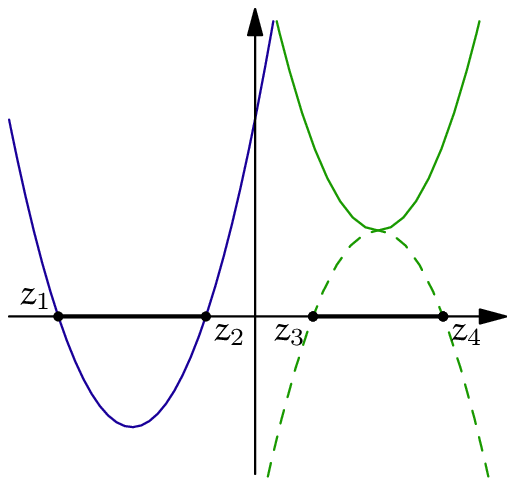}\hfill
  (b) \includegraphics[]{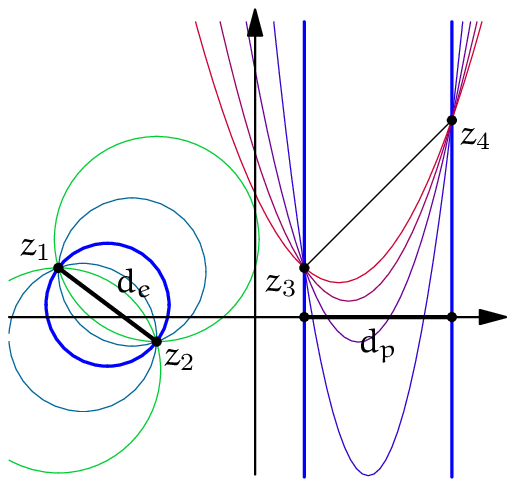}\hfill
  \caption[Radius and distance]{(a) The square of the parabolic
    diameter is the square of the distance between roots if they are
    real (\(z_1\) and \(z_2\)), otherwise the negative square of the
    distance between the adjoint roots (\(z_3\) and \(z_4\)).\\
    (b) Distance as extremum of diameters in elliptic (\(z_1\) and
    \(z_2\)) and parabolic (\(z_3\) and \(z_4\)) cases.}
  \label{fig:distances}
\end{figure}

We mentioned already that for circles normalised by the condition
\(k=1\)%
\index{$k$-normalised cycle}%
\index{normalised!cycle!$k$-}%
\index{cycle!normalised!$k$-} the value \(\det
\cycle{s}{\bs}=\scalar{\cycle{s}{\bs}}{\cycle{s}{\bs}}\) produces the
square of the traditional circle radius. Thus we may keep it as the
definition of the \(\bs\)-\emph{radius}%
\index{cycle!radius}%
\index{radius!cycle, of} for any cycle. But then we need to accept
that in the parabolic case the radius is the (Euclidean) distance
between (real) roots of the parabola, see
Fig.~\ref{fig:distances}(a).

Having radii of circles already defined we may use them for other
measurements in several different ways. For example, the following
variational definition may be used:

\begin{defn}
  \label{de:distance}
  The \emph{distance}\index{distance} between two points is the extremum of diameters
  of all cycles passing through both points, see
  Fig.~\ref{fig:distances}(b).
\end{defn}
  
If \(\bs=\sigma\) this definition gives in all EPH cases the following
expression for a distance \(d_{e,p,h}(u,v)\) between endpoints of any
vector \(w=u+\rmi v\):
\begin{equation}
  \label{eq:eph-distance}
  d_{e,p,h}(u,v)^2=(u+\rmi v)(u-\rmi v)=u^2-\sigma  v^2.
\end{equation}
The parabolic distance \(d_p^2=u^2\), see
Fig.~\ref{fig:distances}(b), algebraically sits between \(d_e\)
and \(d_h\) according to the general
principle~\eqref{eq:eph-class} and is widely
accepted~\cite{Yaglom79}. However one may be unsatisfied by its
degeneracy.

An alternative measurement is motivated by the fact that a circle is
the set of equidistant points from its centre. However the choice of
``centre'' is now rich: it may be either point from three
centres~\eqref{eq:centres} or three foci~\eqref{eq:foci}.
\begin{defn}
  \label{de:length}
  The \emph{length}\index{length} of a directed interval \(\lvec{AB}\)
  is the radius of the cycle with its \emph{centre}%
  \index{length!from centre}%
  \index{centre!length from} (denoted by \(l_c(\lvec{AB})\)) or
  \emph{focus}%
  \index{length!from focus}%
  \index{focus!length from} (denoted by \(l_f(\lvec{AB})\)) at the
  point \(A\) which passes through \(B\).
\end{defn}

This definition is less common and have some unusual properties like
non-symmetry: \(l_f(\lvec{AB})\neq l_f(\lvec{BA})\). However it
comfortably fits the Erlangen programme due to its
\(\SL\)-\emph{conformal invariance}%
\index{conformality}:

\begin{thm}[\cite{Kisil05a}]
  Let \(l\) denote either the EPH
  distances~\eqref{eq:eph-distance} or any length from
  Definition~\ref{de:length}. Then for fixed \(y\),
  \(y'\in\Space{R}{\sigma}\) the limit:
  \begin{displaymath}
    \lim_{t\rightarrow 0} \frac{l(g\cdot y, g\cdot(y+ty'))}{l(y,
      y+ty')}, \qquad
    \text{ where } g\in\SL, 
  \end{displaymath}
  exists and its value depends only from \(y\) and \(g\) and is
  independent from \(y'\).
\end{thm}
  
\begin{figure}[htbp]
  \centering
  \includegraphics[scale=.8]{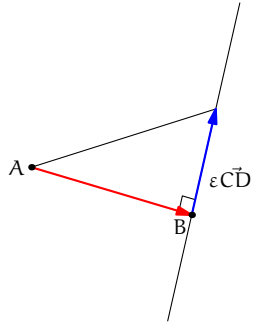}
  \caption{Perpendicular as the shortest route to a line.}
  \label{fig:perpendic}
\end{figure}
We may return from distances to angles recalling that in the Euclidean
space a perpendicular provides the shortest root from a point to a
line, see Fig.~\ref{fig:perpendic}.
\begin{defn}
  \label{de:perpendicular}
  Let \(l\) be a length or distance.  We say that a vector \(\lvec{AB}\) is
  \emph{\(l\)-perpendicular}%
  \index{perpendicular} to a vector \(\lvec{CD}\) if function
  \(l(\lvec{AB}+\varepsilon \lvec{CD})\) of a variable \(\varepsilon\) has a
  local extremum at \(\varepsilon=0\). 
\end{defn}
A pleasant surprise is that \(l_f\)-perpendicularity obtained thought
the length from focus (Definition~\ref{de:length}) coincides with
already defined in Section~\ref{sec:higher-order-joint}
f-orthogonality%
\index{f-orthogonality}%
\index{cycles!f-orthogonal} as follows from
Thm.~\ref{th:ghost2}(\ref{item:focal-centre-rel}). It is also
possible~\cite{Kisil08a} to make \(\SL\) action isometric in all three
cases.

Further details of the refreshing geometry of M\"obius transformation
can be found in the paper~\cite{Kisil05a}and the book~\cite{Kisil12a}.

All these study are waiting to be generalised to high dimensions,
\wiki{Quaternion}{quaternions}\index{quaternion} and
\wiki{Clifford_algebra}{Clifford algebras}%
\index{Clifford!algebra}%
\index{algebra!Clifford} provide a suitable language for
this~\cite{Kisil05a,JParker07a}.

\section{Linear Representations}
\label{sec:induc-repr}

A consideration of the symmetries in analysis is natural to start from
\emph{linear representations}%
\index{representations!linear}. The previous geometrical
actions~\eqref{eq:moebius} can be naturally extended to such
representations by
\wiki{Induced_representations}{induction}%
\index{representation!induced}%
\index{induced!representation}~\citelist{\amscite{Kirillov76}*{\S~13.2}
  \amscite{Kisil97c}*{\S~3.1}} from a representation of a subgroup \(H\).
If \(H\) is one-dimensional then its irreducible representation is a
character, which is always supposed to be a complex valued.  However
hypercomplex number naturally appeared in the \(\SL\)
action~\eqref{eq:moebius}, see
Subsection~\ref{sec:hypercomplex-numbers} and~\cite{Kisil09c}, why
shall we admit only \(\rmi^2=-1\) to deliver a character then?

\subsection{Hypercomplex Characters}
\label{sec:hyperc-char}

\begin{figure}[htbp]
  \centering
  \includegraphics[scale=.8]{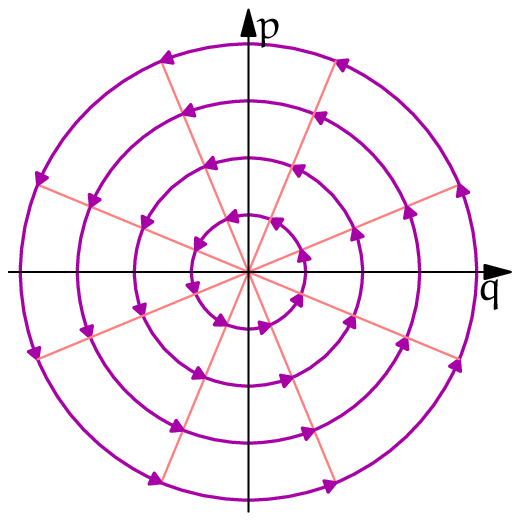}\hfill
  \includegraphics[scale=.8]{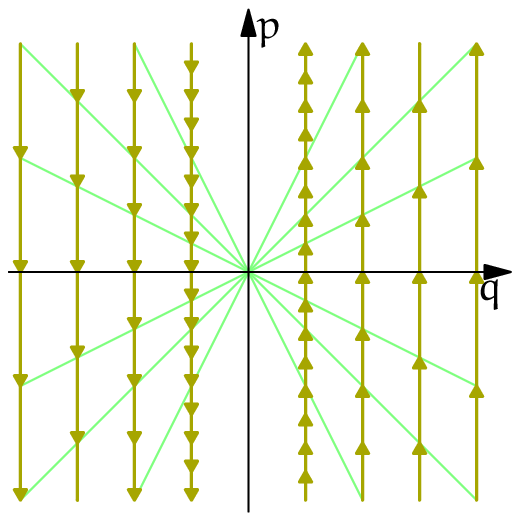}\hfill
  \includegraphics[scale=.8]{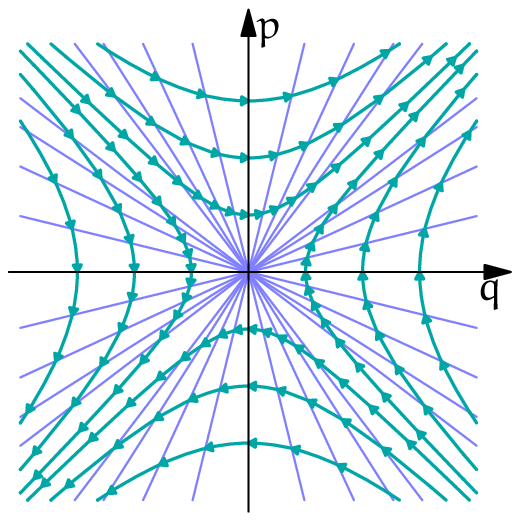}
  \caption[Rotations of wheels]{Rotations of algebraic wheels, i.e.
    the multiplication by \(e^{\alli t}\): elliptic (\(E\)), trivial
    parabolic (\(P_0\)) and hyperbolic (\(H\)). All blue orbits are
    defined by the identity \(x^2-\alli^2y^2=r^2\). Thin ``spokes''
    (straight lines from the origin to a point on the orbit) are
    ``rotated'' from the real axis. This is symplectic linear
    transformations of the classical phase space as well.}
  \label{fig:rotations}
\end{figure}

As we already mentioned the typical discussion of induced
representations of \(\SL\) is centred around the case \(H=K\) and a
complex valued character of \(K\).  A linear transformation defined by
a matrix~\eqref{eq:k-subgroup} in \(K\) is a rotation of
\(\Space{R}{2}\) by the angle \(t\).  After identification
\(\Space{R}{2}=\Space{C}{}\) this action is given by the
multiplication \(e^{\rmi t}\), with \(\rmi^2=-1\).  The rotation
preserve the (elliptic) metric given by:
\begin{equation}
  \label{eq:ell-metric}
  x^2+y^2=(x+\rmi y)(x-\rmi y).
\end{equation}
Therefore the orbits of rotations are circles, any line passing the
origin (a ``spoke'') is rotated by the angle \(t\), see
Fig.~\ref{fig:rotations}.

Dual%
\index{dual!number}%
\index{number!dual} and double numbers%
\index{number!double}%
\index{double!number} produces the most straightforward
adaptation of this result.

\index{symplectic!transformation}%
\index{transformation!symplectic}

\begin{prop}
  \label{pr:algebraic-char}
  The following table show correspondences between three types of
  algebraic characters:
  \begin{center}
    \begin{tabular}{||c|c|c||}
      \hline\hline
      Elliptic & Parabolic & Hyperbolic\\
      \hline\hline
      \(\strut^{\strut}\rmi^{2}=-1\)&    \(\rmp^2=0\)&    \(\rmh^2=1\)
      \\
      \hline
      \(w=x+\rmi y\) &\(w=x+\rmp y\) &\(w=x+\rmh y\) 
      \\
      \hline
      \(\bar{w}=x-\rmi y\) &\(\bar{w}=x-\rmp y\) &\(\bar{w}=x-\rmh y\) 
      \\
      \hline
      \(\strut^{\strut}e^{\rmi t} = \cos t +\rmi \sin t\)&\(e^{\rmp t} = 1 +\rmp t\)&\(e^{\rmh t} = \cosh t +\rmh \sinh t\)
      \\
      \hline
      \(\strut^{\strut}\modulus[e]{w}^{ 2}=w\bar{w}=x^2+y^2\) &\(\modulus[p]{w}^2=w\bar{w}=x^2\) &\(\modulus[h]{w}^2=w\bar{w}=x^2-y^2\) 
      \\
      \hline
      \(\arg w = \tan^{-1} \frac{y}{x} \frac{\strut}{\strut}\)& \(\arg w = \frac{y}{x}\)&\(\arg w = \tanh^{-1} \frac{y}{x}\)
      \\
      \hline
       unit circle \(\strut^{\strut} \modulus[e]{w}^{2}=1\) & ``unit'' strip \(x=\pm 1\) & unit hyperbola \(\modulus[h]{w}^2=1\)
      \\
      \hline\hline
    \end{tabular}
  \end{center}
  Geometrical action of multiplication by \(e^{\alli t}\) is drawn in
  Fig.~\ref{fig:rotations} for all three cases.
\end{prop}
Explicitly parabolic rotations associated with \(\rme^{\rmp t}\) acts on dual
numbers%
\index{dual!number|(}%
\index{number!dual|(} as follows:
\begin{equation}
  \label{eq:parab-rot-triv}
  \rme^{\rmp x}: a+\rmp b \mapsto a+\rmp (a x+b).
\end{equation}
This links the parabolic case with the Galilean
group~\cite{Yaglom79} of symmetries of the classic mechanics, with
the absolute time disconnected from space.

The obvious algebraic similarity and the connection to classical
kinematic is a wide spread justification for the following viewpoint on
the parabolic case, cf.~\cites{HerranzOrtegaSantander99a,Yaglom79}:
\begin{itemize}
\item the parabolic trigonometric functions are trivial:
  \begin{equation}
    \label{eq:par-trig-0}
    \cosp t =\pm 1, \qquad \sinp t=t;
  \end{equation}
\item the parabolic distance is independent from \(y\) if \(x\neq 0\):
  \begin{equation}
    \label{eq:par-metr-0}
    x^2=(x+\rmp y)(x-\rmp y);
  \end{equation}
\item the polar decomposition of a dual number is defined by~\amscite{Yaglom79}*{App.~C(30')}:
  \begin{equation}
    \label{eq:p-polar-yaglom}
    u+\rmp v = u(1+\rmp \frac{v}{u}), \quad \text{ thus }
    \quad \modulus{u+\rmp v}=u, \quad \arg(u+\rmp v)=\frac{v}{u};
  \end{equation}
\item the parabolic wheel looks rectangular, see  Fig.~\ref{fig:rotations}.
\end{itemize}

Those algebraic analogies are quite explicit and widely accepted as an
ultimate source for parabolic
trigonometry~\cites{LavrentShabat77,HerranzOrtegaSantander99a,Yaglom79}.
Moreover, those three rotations are all non-isomorphic symplectic
linear transformations%
\index{symplectic!transformation}%
\index{transformation!symplectic}  of the phase space%
\index{phase!space}%
\index{space!phase}, which makes them useful in
the context of classical and quantum
mechanics~\cites{Kisil10a,Kisil11a}, see
Section~\ref{sec:quantum-mechanics-1}.  There exist also
alternative characters~\cite{Kisil09a} based on M\"obius
transformations with geometric motivation and connections 
to equations of mathematical physics.%
\index{dual!number|)}%
\index{number!dual|)}
\subsection{Induced Representations}
\label{sec:concl-induc-repr}

Let \(G\) be a group, \(H\) be its closed subgroup with the
corresponding homogeneous space \(X=G/H\) with an invariant measure.
We are using notations and definitions of maps \(p: G\rightarrow X\),
\(s:X\rightarrow G\) and \(r: G\rightarrow H\) from
Subsection~\ref{sec:hypercomplex-numbers}.  Let \(\chi\) be an
irreducible representation of \(H\) in a vector space \(V\), then it
induces a representation of \(G\) in the sense of
Mackey~\amscite{Kirillov76}*{\S~13.2}.  This representation has the
realisation \(\uir{}{\chi}\) in the space \(\FSpace{L}{2}(X)\) of
\(V\)-valued functions by the
formula~\amscite{Kirillov76}*{\S~13.2.(7)--(9)}:
\begin{equation} 
  \label{eq:def-ind}
  [\uir{}{\chi}(g) f](x)= \chi(r(g^{-1} * s(x)))  f(g^{-1}\cdot x),
  .
\end{equation}
where \(g\in G\), \(x\in X\), \(h\in H\) and \(r: G
\rightarrow H\), \(s: X \rightarrow G\) are maps defined
above; \(*\)~denotes multiplication on \(G\) and \(\cdot\) denotes the
action~\eqref{eq:g-action} of \(G\) on \(X\).

Consider this scheme for representations of \(\SL\) induced from
characters of its one-dimensional subgroups. We can notice that only
the subgroup \(K\) requires a complex valued character due to the fact
of its compactness. For subgroups \(N'\) and \(\Aprime\) we can
consider characters of all three types---elliptic, parabolic and
hyperbolic.  Therefore we have seven essentially different induced
representations. We will write explicitly only three of them here.

\begin{example}%
  \index{representation!$\SL$ group|(}%
  \index{$\SL$ group!representation|(}
  Consider the subgroup \(H=K\), due to its compactness we are limited
  to complex valued characters of \(K\) only. All of them are of the
  form \(\chi_k\):
  \begin{equation}
    \label{eq:k-character}
    \chi_k\begin{pmatrix}
      \cos t &  \sin t\\
      -\sin t & \cos t
    \end{pmatrix}=e^{-\rmi k t}, \qquad \text{ where }
    k\in\Space{Z}{}.
  \end{equation}
  Using the explicit form~\eqref{eq:s-map} of the map \(s\) we find 
  the map \(r\) given in~\eqref{eq:r-map} as follows:
  \begin{displaymath}
    r
    \begin{pmatrix}
      a&b\\c&d
    \end{pmatrix}
    =\frac{1}{\sqrt{c^2+d^2}}
    \begin{pmatrix}
      d&-c\\c&d
    \end{pmatrix}\in K.
  \end{displaymath}
  Therefore:
  \begin{displaymath}
    r(g^{-1} * s(u,v))  =  
    \frac{1}{\sqrt{(c u+d)^2 +(cv)^2}}
    \begin{pmatrix}
      cu+d&-cv\\cv&cu+d
    \end{pmatrix}, \quad \text{where } g^{-1}=    \begin{pmatrix}
      a&b\\c&d
    \end{pmatrix}.
  \end{displaymath}
  Substituting this into~\eqref{eq:k-character} and combining with the
  M\"obius transformation of the domain~\eqref{eq:moebius} we get the
  explicit realisation \(\uir{}{k}{}\) of the induced representation~\eqref{eq:def-ind}:
  \begin{equation}
    \label{eq:discrete}
    \uir{}{k}{}(g) f(w)=\frac{\modulus{cw+d}^k}{(cw+d)^k}f\left(\frac{aw+b}{cw+d}\right),
    \quad \text{ where } g^{-1}=\begin{pmatrix}a&b\\c&d
    \end{pmatrix}, \ w=u+\rmi v.
  \end{equation}
  This representation acts on complex valued functions in the upper
  half-plane \(\Space[+]{R}{2}=\SL/K\) and belongs to the discrete
  series%
  \index{representation!discrete series}%
  ~\amscite{Lang85}*{\S~IX.2}.
  It is common to get rid of the factor \(\modulus{cw+d}^k\) from that
  expression in order to keep analyticity and we will follow this
  practise for a convenience as well.
\end{example}

\begin{example}
  \label{ex:n-induced}
  In the case of the subgroup \(N\) there is a wider choice of
  possible characters.
  \begin{enumerate}
  \item Traditionally only complex valued characters of the subgroup
    \(N\) are considered, they are:
    \begin{equation}
      \label{eq:np-character}
      \chi_{\tau}^{\Space{C}{}}\begin{pmatrix}
        1 &  0\\
        t & 1
      \end{pmatrix}=e^{\rmi \tau t}, \qquad \text{ where }
      \tau\in\Space{R}{}.
    \end{equation}
    A direct calculation shows that:
    \begin{displaymath}
      r
      \begin{pmatrix}
        a&b\\c&d
      \end{pmatrix}
      =
      \begin{pmatrix}
        1&0\\\frac{c}{d}&1
      \end{pmatrix}\in N'.
    \end{displaymath}
    Thus:
    \begin{equation}
      \label{eq:np-char-part}
      r(g^{-1}*s(u,v))=
      \begin{pmatrix}
        1&0\\\frac{cv}{d+cu}&1
      \end{pmatrix}, \quad\text{ where } g^{-1}=    \begin{pmatrix}
        a&b\\c&d
      \end{pmatrix}.
    \end{equation}
    A substitution of this value into the
    character~\eqref{eq:np-character} together with the M\"obius
    transformation~\eqref{eq:moebius} we obtain the next realisation of~\eqref{eq:def-ind}:
    \begin{displaymath}
      \uir{\Space{C}{}}{\tau}(g) f(w)= \exp\left(\rmi\frac{\tau c v}{cu+d} \right)
      f\left(\frac{aw+b}{cw+d}\right), 
      \quad \text{where } w=u+\rmp v, \  g^{-1}=\begin{pmatrix}a&b\\c&d
      \end{pmatrix}.
    \end{displaymath}
    The representation acts on the space of \emph{complex} valued
    functions on the upper half-plane \(\Space[+]{R}{2}\), which is
    a subset of \emph{dual} numbers%
    \index{dual!number|(}%
    \index{number!dual|(} as a homogeneous space \(\SL/N'\).
    The mixture of complex and dual numbers in the same expression is
    confusing.
  \item The parabolic character \(\chi_{\tau}\) with the algebraic flavour
    is provided by multiplication~\eqref{eq:parab-rot-triv} with the
    dual number:
    \begin{displaymath}
      \chi_{\tau}\begin{pmatrix}
        1 &  0\\
        t & 1
      \end{pmatrix}=e^{\rmp \tau t}=1+\rmp \tau t, \qquad \text{ where }
      \tau\in\Space{R}{}.
    \end{displaymath}
    If we substitute the value~\eqref{eq:np-char-part} into this
    character, then we receive the representation:
    \begin{displaymath}
      \uir{}{\tau}{}(g) f(w)= \left(1+\rmp\frac{\tau c v}{cu+d} \right)
      f\left(\frac{aw+b}{cw+d}\right),
    \end{displaymath}
    where \(w\), \(\tau\) and \(g\) are as above.  The representation
    is defined on the space of dual numbers valued functions on the
    upper half-plane of dual numbers.  Thus expression contains only
    dual numbers with their usual algebraic operations. Thus it is
    linear with respect to them.
  \end{enumerate}
\end{example}
All characters in the previous Example are unitary. Then the general
scheme of induced representations~\amscite{Kirillov76}*{\S~13.2} implies
their unitarity in proper senses.
\begin{thm}[\cite{Kisil09c}]
  \label{th:unitarity}
  Both representations of \(\SL\) from Example~\ref{ex:n-induced}
  are unitary on the space of function on the upper half-plane
  \(\Space[+]{R}{2}\) of dual numbers with the inner product:
  \begin{equation}
    \label{eq:inner-product}
    \scalar{f_1}{f_2}=\int_{\Space[+]{R}{2}} f_1(w)
    \bar{f}_2(w)\,\frac{du\,dv}{v^2}, \qquad \text{ where } w=u+\rmp v,
  \end{equation}
  and we use the conjugation and multiplication of functions' values
  in algebras of complex and dual numbers for representations
  \(\uir{\Space{C}{}}{\tau}\) and \(\uir{}{\tau}\) respectively.
\end{thm}
The inner product~\eqref{eq:inner-product} is positive defined for
the representation \(\uir{\Space{C}{}}{\tau}\) but is not for the
other. The respective spaces are parabolic cousins of the \emph{Krein
  spaces}%
\index{Krein!space}%
\index{space!Krein}~\cite{ArovDym08}, which are hyperbolic in our sense.
  \index{representation!$\SL$ group|)}%
  \index{$\SL$ group!representation|)}%
\index{dual!number|)}%
\index{number!dual|)}

\subsection{Similarity and Correspondence: Ladder Operators}
\label{sec:correspondence}

From the above observation we can deduce the following empirical
principle, which has a heuristic value.

\begin{principle}[Similarity and correspondence]
  \index{similarity|see{principle of similarity and correspondence}}
  \index{correspondence|see{principle of similarity and correspondence}}
  \index{principle!similarity and correspondence}
  \label{pr:simil-corr-principle}
  \begin{enumerate}
  \item Subgroups \(K\), \(N'\) and \(\Aprime\) play a similar r\^ole in the
    structure of the group \(\SL\) and its representations.
  \item The subgroups shall be swapped simultaneously with the
    respective replacement of hypercomplex%
    \index{number!hypercomplex}%
    \index{hypercomplex!number} unit \(\alli\).
  \end{enumerate}
\end{principle}
The first part of the Principle (similarity) does not look sound
alone. It is enough to mention that the subgroup \(K\) is compact (and
thus its spectrum is discrete) while two other subgroups are not.
However in a conjunction with the second part (correspondence) the
Principle have received the following confirmations so far,
see~\cite{Kisil09c} for details:
\begin{itemize}
\item The action of \(\SL\) on the homogeneous space \(\SL/H\) for
  \(H=K\), \(N'\) or \(\Aprime\) is given by linear-fractional
  transformations of complex, dual or double numbers respectively.  
\item Subgroups \(K\), \(N'\) or \(\Aprime\) are isomorphic to the groups of
  unitary rotations of respective unit cycles in complex, dual or
  double numbers.   
\item Representations induced from subgroups \(K\), \(N'\) or \(\Aprime\)
  are unitary if the inner product spaces of functions with values in
  complex, dual or double numbers.
\end{itemize}
\begin{rem}
  The principle of similarity and correspondence resembles
  \emph{supersymmetry}%
  \index{supersymmetry} between bosons and fermions in particle physics, but
  we have similarity between three different types of entities in our case.
\end{rem}
Let us give another illustration to the Principle. Consider the Lie
algebra \(\algebra{sl}_2\) of the group \(\SL\). Pick up the following
basis in \(\algebra{sl}_2\)~\amscite{MTaylor86}*{\S~8.1}:  
\begin{equation}
  \label{eq:sl2-basis}
  A= \frac{1}{2}
  \begin{pmatrix}
    -1&0\\0&1
  \end{pmatrix},\quad 
  B= \frac{1}{2} \
  \begin{pmatrix}
    0&1\\1&0
  \end{pmatrix}, \quad 
  Z=
  \begin{pmatrix}
    0&1\\-1&0
  \end{pmatrix}.
\end{equation}
The commutation relations between the elements are:
\begin{equation}
  \label{eq:sl2-commutator}
  [Z,A]=2B, \qquad [Z,B]=-2A, \qquad [A,B]=- \frac{1}{2} Z.
\end{equation} 
Let \(\uir{}{}\) be a representation of the group \(\SL\) in a space
\(V\). Consider the derived representation \(d\uir{}{}\) of the Lie
algebra \(\algebra{sl}_2\)~\amscite{Lang85}*{\S~VI.1} and denote
\(\tilde{X}=d\uir{}{}(X)\) for \(X\in\algebra{sl}_2\). To see the
structure of the representation \(\uir{}{}\) we can decompose the
space \(V\) into eigenspaces of the operator \(\tilde{X}\) for some
\(X\in \algebra{sl}_2\), cf. the Taylor series in
Section~\ref{sec:taylor-expansion}.

\subsubsection{Elliptic Ladder Operators}
\label{sec:ellipt-ladd-oper}
It would not be surprising that we are going to consider three cases:
Let \(X=Z\) be a generator of the subgroup%
\index{generator!of a subgroup}%
\index{subgroup!generator}
\(K\)~\eqref{eq:k-subgroup}. Since this is a compact subgroup the
corresponding eigenspaces%
\index{eigenvalue} \(\tilde{Z} v_k=\rmi k v_k\) are
parametrised by an integer \(k\in\Space{Z}{}\).  The
\emph{raising/lowering} or \emph{ladder operators}%
\index{ladder operator|(}%
\index{operator!ladder|(}%
\index{operator!raising|see{ladder operator}}%
\index{operator!lowering|see{ladder operator}}%
\index{operator!creation|see{ladder operator}}%
\index{operator!annihilation|see{ladder operator}}%
\index{raising operator|see{ladder operator}}%
\index{lowering operator|see{ladder operator}}%
\index{creation operator|see{ladder operator}}%
\index{annihilation operator|see{ladder operator}}
\(\ladder{\pm}\)~\citelist{\amscite{Lang85}*{\S~VI.2}
  \amscite{MTaylor86}*{\S~8.2}} are defined by the following
commutation relations:
\begin{equation}
  \label{eq:raising-lowering}
  [\tilde{Z},\ladder{\pm}]=\lambda_\pm \ladder{\pm}. 
\end{equation}
In other words \(\ladder{\pm}\) are eigenvectors for operators 
\(\loglike{ad}Z\) of adjoint representation of \(\algebra{sl}_2\)~\amscite{Lang85}*{\S~VI.2}.
\begin{rem}
  The existence of such ladder operators follows from the general
  properties of Lie algebras if the element \(X\in\algebra{sl}_2\)
  belongs to a \emph{Cartan subalgebra}%
  \index{Cartan!subalgebra}%
  \index{subalgebra!Cartan}. This is the case for vectors \(Z\)
  and \(B\), which are the only two non-isomorphic types of
  Cartan subalgebras in \(\algebra{sl}_2\). However the third case
  considered in this paper, the parabolic vector \(B+Z/2\), does
  not belong to a Cartan subalgebra, yet a sort of ladder
  operators is still possible with dual number coefficients.
  Moreover, for the hyperbolic vector \(B\), besides the standard
  ladder operators an additional pair with double number
  coefficients will also be described.
\end{rem}

From the commutators~\eqref{eq:raising-lowering} we deduce that
\(\ladder{+} v_k\) are eigenvectors of \(\tilde{Z}\) as well:
\begin{eqnarray}
  \tilde{Z}(\ladder{+} v_k)&=&(\ladder{+}\tilde{Z}+\lambda_+\ladder{+})v_k=\ladder{+}(\tilde{Z}v_k)+\lambda_+\ladder{+}v_k
  =\rmi k \ladder{+}v_k+\lambda_+\ladder{+}v_k\nonumber \\
  &=&(\rmi k+\lambda_+)\ladder{+}v_k.\label{eq:ladder-action}
\end{eqnarray}
Thus action of ladder operators on respective eigenspaces can be
visualised by the diagram:
\begin{equation}
  \label{eq:ladder-chain-1D}
  \xymatrix@1{
    \ldots\, \ar@<.4ex>[r]^{\ladder{+}} &
    \,V_{\rmi k-\lambda}\,  \ar@<.4ex>[l]^{\ladder{-}}\ar@<.4ex>[r]^{\ladder{+}} &
    \,V_{\rmi k}\, \ar@<.4ex>[l]^{\ladder{-}} \ar@<.4ex>[r]^{\ladder{+}} &
    \,V_{\rmi k+ \lambda}\,\ar@<.4ex>[l]^{\ladder{-}}  \ar@<.4ex>[r]^{\ladder{+}}
    &
    \,\ldots\ar@<.4ex>[l]^{\ladder{-}}}
\end{equation}
Assuming \(\ladder{+}=a\tilde{A}+b\tilde{B}+c\tilde{Z}\) from the
relations~\eqref{eq:sl2-commutator} and defining
condition~\eqref{eq:raising-lowering} we obtain linear equations
with unknown \(a\), \(b\) and \(c\): 
\begin{displaymath}
  c=0, \qquad 2a=\lambda_+ b, \qquad -2b=\lambda_+ a.
\end{displaymath}
The equations have a solution if and only if \(\lambda_+^2+4=0\),
and the raising/lowering operators are \(\ladder{\pm}=\pm\rmi \tilde{A}+\tilde{B}\).

\subsubsection{Hyperbolic Ladder Operators}
\label{sec:hyperb-ladd-oper}
Consider the case \(X=2B\) of a generator of the subgroup
\(\Aprime\)~\eqref{eq:ap-subgroup}. The subgroup is not compact
and eigenvalues of the operator \(\tilde{B}\) can be arbitrary,
however raising/lowering operators are still
important~\citelist{\amscite{HoweTan92}*{\S~II.1}
  \amscite{Mazorchuk09a}*{\S~1.1}}.  We again seek a solution in the
form \(\ladder[h]{+}=a\tilde{A}+b\tilde{B}+c\tilde{Z}\) for the commutator
\([2\tilde{B},\ladder[h]{+}]=\lambda \ladder[h]{+}\). We will get the system:
\begin{displaymath}
  4c=\lambda a,\qquad
  b=0,\qquad
  {a}=\lambda c.
\end{displaymath}
A solution exists if and only if \(\lambda^2=4\). There are
obvious values \(\lambda=\pm 2\) with the ladder operators
\(\ladder[h]{\pm}=\pm2\tilde{A}+\tilde{Z}\),
see~\citelist{\amscite{HoweTan92}*{\S~II.1}
  \amscite{Mazorchuk09a}*{\S~1.1}}. Each indecomposable
\(\algebra{sl}_2\)-module is formed by a one-dimensional chain of
eigenvalues\index{eigenvalue} with a transitive action of ladder
operators.

\begin{figure}[htbp]
  \centering
  \(  \xymatrix@R=2.5em@C=1.5em@M=.5em{
    & 
    \,\ldots\, \ar@<.4ex>[d]^{\ladder[\rmh]{+}} &  
    \,\ldots\, \ar@<.4ex>[d]^{\ladder[\rmh]{+}} & 
    \,\ldots\,  \ar@<.4ex>[d]^{\ladder[\rmh]{+}}  & 
    \\
    \ldots\, \ar@<.4ex>[r]^-{\ladder[h]{+}} & 
    \,V_{(n-2)+\rmh (k-2)}\,  \ar@<.4ex>[l]^-{\ladder[h]{-}}\ar@<.4ex>[r]^{\ladder[h]{+}}
    \ar@<.4ex>[u]^{\ladder[\rmh]{-}} \ar@<.4ex>[d]^{\ladder[\rmh]{+}} &  
    \,V_{n+\rmh (k-2)}\, \ar@<.4ex>[l]^{\ladder[h]{-}} \ar@<.4ex>[r]^{\ladder[h]{+}}
    \ar@<.4ex>[u]^{\ladder[\rmh]{-}} \ar@<.4ex>[d]^{\ladder[\rmh]{+}} & 
    \,V_{(n+2)+\rmh (k-2)}\,\ar@<.4ex>[l]^{\ladder[h]{-}}  \ar@<.4ex>[r]^-{\ladder[h]{+}}
    \ar@<.4ex>[u]^{\ladder[\rmh]{-}} \ar@<.4ex>[d]^{\ladder[\rmh]{+}}    & 
    \,\ldots\ar@<.4ex>[l]^-{\ladder[h]{-}}\\
    \ldots\, \ar@<.4ex>[r]^-{\ladder[h]{+}} & 
    \,V_{(n-2)+\rmh k}\,  \ar@<.4ex>[l]^-{\ladder[h]{-}}\ar@<.4ex>[r]^{\ladder[h]{+}}
    \ar@<.4ex>[u]^{\ladder[\rmh]{-}} \ar@<.4ex>[d]^{\ladder[\rmh]{+}} &  
    \,V_{n+\rmh k}\, \ar@<.4ex>[l]^{\ladder[h]{-}} \ar@<.4ex>[r]^{\ladder[h]{+}}
    \ar@<.4ex>[u]^{\ladder[\rmh]{-}} \ar@<.4ex>[d]^{\ladder[\rmh]{+}}& 
    \,V_{(n+2)+\rmh k}\,\ar@<.4ex>[l]^{\ladder[h]{-}}  \ar@<.4ex>[r]^-{\ladder[h]{+}}
    \ar@<.4ex>[u]^{\ladder[\rmh]{-}} \ar@<.4ex>[d]^{\ladder[\rmh]{+}}    & 
    \,\ldots\ar@<.4ex>[l]^-{\ladder[h]{-}}\\
    \ldots\, \ar@<.4ex>[r]^-{\ladder[h]{+}} & 
    \,V_{(n-2)+\rmh (k+2)}\,  \ar@<.4ex>[l]^-{\ladder[h]{-}}\ar@<.4ex>[r]^{\ladder[h]{+}}
    \ar@<.4ex>[u]^{\ladder[\rmh]{-}} \ar@<.4ex>[d]^{\ladder[\rmh]{+}} &  
    \,V_{n+\rmh (k+2)}\, \ar@<.4ex>[l]^{\ladder[h]{-}} \ar@<.4ex>[r]^{\ladder[h]{+}} 
    \ar@<.4ex>[u]^{\ladder[\rmh]{-}} \ar@<.4ex>[d]^{\ladder[\rmh]{+}}& 
    \,V_{(n+2)+\rmh (k+2)}\,\ar@<.4ex>[l]^{\ladder[h]{-}}  \ar@<.4ex>[r]^-{\ladder[h]{+}}
    \ar@<.4ex>[u]^{\ladder[\rmh]{-}} \ar@<.4ex>[d]^{\ladder[\rmh]{+}}    & 
    \,\ldots\ar@<.4ex>[l]^-{\ladder[h]{-}}\\
    & 
    \,\ldots\, \ar@<.4ex>[u]^{\ladder[\rmh]{-}} &  
    \,\ldots\, \ar@<.4ex>[u]^{\ladder[\rmh]{-}} & 
    \,\ldots\,  \ar@<.4ex>[u]^{\ladder[\rmh]{-}}  & }
  \)
  \caption[The action of hyperbolic ladder operators on a 2D
  lattice of eigenspaces]{The action of hyperbolic ladder operators on a 2D
    lattice of eigenspaces. Operators \(\ladder[h]{\pm}\) move the
    eigenvalues by \(2\), 
    making shifts in the horizontal direction. Operators
    \(\ladder[\rmh]{\pm}\) change the eigenvalues by \(2\rmh\), 
    shown as vertical shifts.}  
  \label{fig:2D-lattice}
\end{figure}
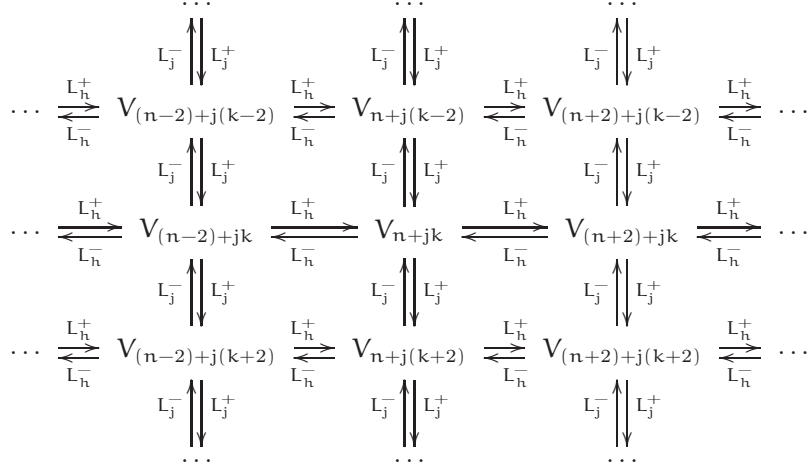

Admitting double numbers%
\index{number!double}%
\index{double!number} we have an extra possibility to satisfy
\(\lambda^2=4\) with values \(\lambda=\pm2\rmh\).  Then there is an
additional pair of hyperbolic ladder operators
\(\ladder[\rmh]{\pm}=\pm2\rmh\tilde{A}+\tilde{Z}\), which shift eigenvectors
in the ``orthogonal'' direction to the standard operators \(\ladder[h]{\pm}\).
Therefore an indecomposable \(\algebra{sl}_2\)-module can be
parametrised by a two-dimensional lattice of eigenvalues on the
double number plane,  see
Fig.~\ref{fig:2D-lattice}

\subsubsection{Parabolic Ladder Operators}
\label{sec:parab-ladd-oper}

Finally consider the case of a generator \(X=-B+Z/2\) of the
subgroup \(N'\)~\eqref{eq:np-subgroup}. According to the above
procedure we get the equations:
\begin{displaymath}
  b+2c=\lambda a,\qquad
  -a=\lambda b,\qquad
  \frac{a}{2}=\lambda c,
\end{displaymath}
which can be resolved if and only if \(\lambda^2=0\). If we restrict
ourselves with the only real (complex) root \(\lambda=0\), then the
corresponding operators \(\ladder[p]{\pm}=-\tilde{B}+\tilde{Z}/2\)
will not affect eigenvalues and thus are useless in the above context.
However the dual number roots \(\lambda =\pm \rmp t\),
\(t\in\Space{R}{}\) lead to the operators \(\ladder[\rmp]{\pm}=\pm
\rmp t\tilde{A}-\tilde{B}+\tilde{Z}/2\). These operators are suitable
to build an \(\algebra{sl}_2\)-modules with a one-dimensional chain of
eigenvalues.
\begin{rem}
  \label{re:hyper-number-necessity}
  The following r\^oles of hypercomplex numbers are noteworthy%
  \index{number!hypercomplex}%
  \index{hypercomplex!number}:
  \begin{itemize}
  \item the introduction of complex numbers is a necessity for the
    \emph{existence} of ladder operators in the elliptic
    case;
  \item in the parabolic case we need dual numbers%
    \index{dual!number}%
    \index{number!dual} to make
    ladder operators \emph{useful};
  \item in the hyperbolic case double numbers%
    \index{number!double}%
    \index{double!number} are not required 
    neither for the existence or for the usability of ladder operators, but
    they do provide an enhancement. 
  \end{itemize}
\end{rem}
We summarise the above consideration with a focus on the Principle of
similarity and correspondence:%
\index{principle!similarity and correspondence}
\begin{prop}
  \label{pr:ladder-sim-eq}
  Let a vector \(X\in\algebra{sl}_2\) generates the subgroup \(K\),
  \(N'\) or \(\Aprime\), that is \(X=Z\), \(B-Z/2\), or
  \(B\) respectively. Let \(\alli\) be the respective hypercomplex unit.   

  Then raising/lowering operators \(\ladder{\pm}\) satisfying to the
  commutation relation:
  \begin{displaymath}
    [X,\ladder{\pm}]=\pm\alli \ladder{\pm},\qquad [\ladder{-},\ladder{+}]=2\alli X.
  \end{displaymath}
  are:
  \begin{displaymath}
    \ladder{\pm}=\pm\alli \tilde{A} +\tilde{Y}.
  \end{displaymath}
  Here \(Y\in\algebra{sl}_2\) is a linear combination of  \(B\) and
  \(Z\) with the properties:
  \begin{itemize}
  \item \(Y=[A,X]\).
  \item \(X=[A,Y]\).
  \item Killings form \(K(X,Y)\)~\amscite{Kirillov76}*{\S~6.2} vanishes.
  \end{itemize}
  Any of the above properties defines the vector \(Y\in\loglike{span}\{B,Z\}\)
  up to a real constant factor.
\end{prop}

The usability of the Principle of similarity and correspondence will
be illustrated by more examples below.%
\index{ladder operator|)}%
\index{operator!ladder|)}

\section{Covariant Transform}
\label{sec:covariant-transform}
A general group-theoretical
construction~\cites{Perelomov86,FeichGroech89a,Kisil98a,AliAntGaz00,%
Fuhr05a,ChristensenOlafsson09a,KlaSkag85}
of \emph{wavelets}\index{wavelet} (or \emph{coherent state}%
\index{coherent state|see{wavelet}}%
\index{state!coherent|see{wavelet}}) starts from an
irreducible square integrable  representation%
\index{square integrable!representation}%
\index{representation!square integrable}---in the
proper sense or modulo a subgroup.  Then a mother wavelet is chosen to
be \emph{admissible}%
\index{admissible wavelet}%
\index{wavelet!admissible}. This leads to a wavelet transform which is an
isometry to \(\FSpace{L}{2}\) space with respect to the Haar measure%
\index{invariant!measure}%
\index{measure!invariant}
on the group or (quasi)invariant measure on a homogeneous space.

The importance of the above situation shall not be diminished, however
an exclusive restriction to such a setup is not necessary, in fact.
Here is a classical example from complex analysis: the Hardy space
\(\FSpace{H}{2}(\Space{T}{})\)%
\index{space!Hardy}%
\index{Hardy!space} on the unit circle and Bergman spaces%
\index{space!Bergman}%
\index{Bergman!space}
\(\FSpace[n]{B}{2}(\Space{D}{})\), \(n\geq 2\) in the unit disk produce wavelets
associated with representations \(\rho_1\) and \(\rho_n\) of the group
\(\SL\) respectively~\cite{Kisil97c}. While representations \(\rho_n\), \(n\geq 2\)
are from square integrable discrete series%
\index{representation!discrete series}, the mock discrete series
representation \(\rho_1\) is not square
integrable~\citelist{\amscite{Lang85}*{\S~VI.5}
  \amscite{MTaylor86}*{\S~8.4}}. However it would be natural to treat the
Hardy space in the same framework as Bergman ones. Some more examples
will be presented below.

\subsection{Extending Wavelet Transform}
\label{sec:wavelet-transform}

To make a sharp but still natural generalisation of wavelets\index{wavelet} we give the following
definition.
\begin{defn}\cite{Kisil09d}
  Let \(\uir{}{}\) be a representation of
  a group \(G\) in a space \(V\) and \(F\) be an operator from \(V\) to a space
  \(U\). We define a \emph{covariant transform}%
  \index{covariant!transform}%
  \index{transform!covariant}%
  \index{transform!wavelet|see{wavelet transform}}
  \(\oper{W}\) from \(V\) to the space \(\FSpace{L}{}(G,U)\) of
  \(U\)-valued functions on \(G\) by the formula:
  \begin{equation}
    \label{eq:coheret-transf-gen}
    \oper{W}: v\mapsto \hat{v}(g) = F(\uir{}{}(g^{-1}) v), \qquad
    v\in V,\ g\in G.
  \end{equation}
  Operator \(F\) will be called \emph{fiducial operator}%
  \index{fiducial operator}%
  \index{operator!fiducial} in this context.
\end{defn}
We borrow the name for operator \(F\) from fiducial vectors of
Klauder and Skagerstam~\cite{KlaSkag85}. 
\begin{rem}
  We do not require that fiducial operator \(F\) shall be linear.
  Sometimes the positive homogeneity, i.e. \(F(t v)=tF(v)\) for
  \(t>0\), alone can be already sufficient, see
  Example~\ref{ex:maximal-function}.
\end{rem}
\begin{rem}
  \label{re:range-dim} 
  Usefulness of the covariant transform is in the reverse proportion
  to the dimensionality of the space \(U\). The covariant transform
  encodes properties of \(v\) in a function \(\oper{W}v\) on \(G\).
  For a low dimensional \(U\) this function can be ultimately
  investigated by means of harmonic analysis. Thus \(\dim U=1\)
  (scalar-valued functions) is the ideal case, however, it is
  unattainable sometimes, see Example~\ref{it:direct-funct} below. We
  may have to use a higher dimensions of \(U\) if the given group
  \(G\) is not rich enough.
\end{rem}
As we will see below covariant transform is a close relative of
wavelet transform.  The name is chosen due to the following common
property of both transformations.
\begin{thm} 
  \label{pr:inter1} 
  The covariant transform~\eqref{eq:coheret-transf-gen}
  intertwines%
  \index{intertwining operator}%
  \index{operator!intertwining} \(\uir{}{}\) and the left regular representation
  \(\Lambda\)    on \(\FSpace{L}{}(G,U)\):
  \begin{displaymath}
    \oper{W} \uir{}{}(g) = \Lambda(g) \oper{W}.
  \end{displaymath}
  Here \(\Lambda\) is defined as usual by:
  \begin{equation}\label{eq:left-reg-repr}
    \Lambda(g): f(h) \mapsto f(g^{-1}h).
  \end{equation}
\end{thm}
\begin{proof}
  We have a calculation similar to wavelet
  transform~\amscite{Kisil98a}*{Prop.~2.6}. Take \(u=\uir{}{}(g) v\) and
  calculate its covariant transform:
   \begin{eqnarray*}{}
     [\oper{W}( \uir{}{}(g) v)] (h) & = &  [\oper{W}(\uir{}{}(g) v)] (h)=F(\uir{}{}(h^{-1}) \uir{}{}(g) v ) \\
     & = & F(\uir{}{}((g^{-1}h)^{-1}) v) \\
     & = & [\oper{W}v] (g^{-1}h)\\
     & = & \Lambda(g) [\oper{W}v] (h).
   \end{eqnarray*}
\end{proof}
The next result follows immediately:
\begin{cor}\label{co:pi}
  The image space \(\oper{W}(V)\) is invariant under the
  left shifts on \(G\).
\end{cor}

\begin{rem}
  A further generalisation of the covariant transform can be obtained
  if we relax the group structure. Consider, for example, a
  \emph{cancellative semigroup}%
  \index{cancellative semigroup}%
  \index{semigroup!cancellative} \(\Space[+]{Z}{}\) of non-negative
  integers. It has a linear presentation on the space of polynomials
  in a variable \(t\) defined by the action \(m: t^n \mapsto t^{m+n}\)
  on the monomials. Application of a linear functional \(l\), e.g.
  defined by an integration over a measure on the real line, produces
  \emph{umbral calculus}%
  \index{umbral calculus}%
  \index{calculus!umbral} \(l(t^n)=c_n\), which has a magic efficiency
  in many areas, notably in combinatorics~\cites{Rota95,Kisil97b}. In
  this direction we also find fruitful to expand the notion of an
  intertwining operator%
  \index{intertwining operator}%
  \index{operator!intertwining} to a
  \emph{token}\index{token}~\cite{Kisil01b}.
\end{rem}

\subsection{Examples of Covariant Transform}
\label{sec:exampl-covar-transf}
In this Subsection we will provide several examples of covariant
transforms. Some of them will be expanded in subsequent sections,
however a detailed study of all aspects will not fit into the present
work.  We start from the classical example of the group-theoretical
wavelet transform:
\begin{example}
  \label{ex:wavelet}
  Let \(V\) be a Hilbert space with an inner product
  \(\scalar{\cdot}{\cdot}\) and \(\uir{}{}\) be a unitary
  representation of a group \(G\) in the space \(V\). Let \(F: V
  \rightarrow \Space{C}{}\) be a functional \(v\mapsto
  \scalar{v}{v_0}\) defined by a vector \(v_0\in V\). The vector
  \(v_0\) is oftenly called the \emph{mother wavelet}%
  \index{mother wavelet}%
  \index{wavelet!mother} in areas related
  to signal processing or the \emph{vacuum state}%
  \index{vacuum state|see{mother wavelet}}%
  \index{state!vacuum|see{mother wavelet}} in quantum
  framework.

  Then the transformation~\eqref{eq:coheret-transf-gen} is the
  well-known expression for a \emph{wavelet
    transform}%
  \index{wavelet!transform}~\amscite{AliAntGaz00}*{(7.48)} (or \emph{representation
    coefficients}%
  \index{representation!coefficients|see{wavelet transform}}):
  \begin{equation}
    \label{eq:wavelet-transf}
    \oper{W}: v\mapsto \hat{v}(g) = \scalar{\uir{}{}(g^{-1})v}{v_0}  =
    \scalar{ v}{\uir{}{}(g)v_0}, \qquad
    v\in V,\ g\in G.
  \end{equation}
  The family of vectors \(v_g=\uir{}{}(g)v_0\) is called
  \emph{wavelets}\index{wavelet} or \emph{coherent states}. In this case we obtain
  scalar valued functions on \(G\), thus the fundamental r\^ole of
  this example is explained in Rem.~\ref{re:range-dim}.

  This scheme is typically carried out for a square integrable
  representation \(\uir{}{}\)%
  \index{square integrable!representation}%
  \index{representation!square integrable} and \(v_0\) being an admissible
  vector%
  \index{admissible wavelet}%
  \index{wavelet!admissible}~\cites{Perelomov86,FeichGroech89a,AliAntGaz00,Fuhr05a,ChristensenOlafsson09a}.
  In this case the wavelet (covariant) transform is a map into the
  square integrable functions~\cite{DufloMoore} with respect to the
  left Haar measure.%
   \index{invariant!measure}%
   \index{measure!invariant} The map becomes an isometry if \(v_0\) is properly
  scaled.
\end{example}
However square integrable representations and admissible vectors does not cover all
interesting cases.
\begin{example}
  \label{ex:ax+b}
  Let \(G=\mathrm{Aff}\) be the ``\(ax+b\)''%
  \index{$ax+b$ group}%
  \index{group!$ax+b$} (or \emph{affine})
  group~\amscite{AliAntGaz00}*{\S~8.2}: the set of points \((a,b)\),
  \(a\in \Space[+]{R}{}\), \(b\in \Space{R}{}\) in the upper
  half-plane with the group law:
   \begin{equation}
     (a, b) * (a', b') = (aa', ab'+b)
   \end{equation}
   and left invariant measure \(a^{-2}\,da\,db\)%
   \index{$ax+b$ group!invariant measure}%
   \index{group!$ax+b$!invariant measure}%
   \index{invariant!measure}%
   \index{measure!invariant}.  Its isometric
   representation on \(V=\FSpace{L}{p}(\Space{R}{})\) is given by the
   formula:%
   \index{$ax+b$ group!representation}%
   \index{representation!$ax+b$ group}
  \begin{equation}\label{eq:ax+b-repr-1}
     [\uir{}{p}(g)\, f](x)= a^{\frac{1}{p}}f\left(ax+b\right),
     \qquad\text{where } g^{-1}=(a,b).
  \end{equation}
  We consider the operators \(F_{\pm}:\FSpace{L}{2}(\Space{R}{})
  \rightarrow \Space{C}{}\) defined by:
  \begin{equation}
    \label{eq:cauchy-pm}
    F_{\pm}(f)=\frac{1}{2\pi i}\int_{\Space{R}{}} \frac{f(t)\,dt}{x\mp \rmi}.
  \end{equation}
  Then the covariant transform~\eqref{eq:coheret-transf-gen} is the
  Cauchy integral%
  \index{integral!Cauchy}%
  \index{Cauchy!integral} from \(\FSpace{L}{p}(\Space{R}{})\) to the space of
  functions \(\hat{f}(a,b)\) such that
  \(a^{-\frac{1}{p}}\hat{f}(a,b)\) is in the Hardy space%
  \index{space!Hardy}%
  \index{Hardy!space} in the
  upper/lower half-plane \(\FSpace{H}{p}(\Space[\pm]{R}{2})\).
  Although the representation~\eqref{eq:ax+b-repr-1} is square
  integrable for \(p=2\), the function \(\frac{1}{x\pm \rmi}\) used
  in~\eqref{eq:cauchy-pm} is not an admissible vacuum vector.  Thus the
  complex analysis become decoupled from the traditional wavelet
  theory. As a result the application of wavelet theory shall relay on
  an extraneous mother wavelets~\cite{Hutnik09a}.

  Many important objects in complex analysis are generated by
  inadmissible mother wavelets like~\eqref{eq:cauchy-pm}. For example, if
  \(F:\FSpace{L}{2}(\Space{R}{}) \rightarrow \Space{C}{}\) is defined
  by \(F: f \mapsto F_+ f + F_-f\) then the covariant
  transform~\eqref{eq:coheret-transf-gen} reduces to the \emph{Poisson
    integral}%
  \index{Poisson!kernel}%
  \index{kernel!Poisson}%
  \index{integral!Poisson|see{Poisson kernel}}.  If \(F:\FSpace{L}{2}(\Space{R}{}) \rightarrow
  \Space{C}{2}\) is defined by \(F: f \mapsto( F_+ f, F_-f)\) then the
  covariant transform~\eqref{eq:coheret-transf-gen} represents a
  function \(f\) on the real line as a jump:
  \begin{equation}
    \label{eq:jump-rl}
    f(z)=f_+(z)-f_-(z),\qquad f_\pm(z)\in \FSpace{H}{p}(\Space[\pm]{R}{2})
  \end{equation}
  between functions analytic in the upper and the lower half-planes.
  This makes a decomposition of \(\FSpace{L}{2}(\Space{R}{})\) into
  irreducible components of the representation~\eqref{eq:ax+b-repr-1}.
  Another interesting but non-admissible vector is the \emph{Gaussian}\index{Gaussian}
  \(e^{-x^2}\).
\end{example}
\begin{example}
  \label{ex:sl2}
  For the group \(G=\SL\)~\cite{Lang85} let us consider the unitary
  representation \(\uir{}{}\) on the space of square integrable function
  \(\FSpace{L}{2}(\Space[+]{R}{2})\) on the upper half-plane through
  the M\"obius transformations%
  \index{M\"obius map}%
  \index{map!M\"obius}~\eqref{eq:moebius}:
   \begin{equation}
     \label{eq:sl2-action}
     \uir{}{}(g): f(z) \mapsto \frac{1}{(c z + d)^2}\,
     f\left(\frac{a z+ b }{c z +d}\right), \qquad g^{-1}=\
     \begin{pmatrix}
       a & b \\ c & d 
     \end{pmatrix}.
   \end{equation}
   This is a representation from the discrete series%
   \index{representation!discrete series} and
   \(\FSpace{L}{2}(\Space{D}{})\) and irreducible invariant subspaces
   are parametrised by integers.  Let \(F_k\) be the functional
   \(\FSpace{L}{2}(\Space[+]{R}{2})\rightarrow \Space{C}{}\) of
   pairing with the lowest/highest \(k\)-weight vector in the
   corresponding irreducible component (Bergman space)%
   \index{space!Bergman}%
   \index{Bergman!space}
   \(\FSpace{B}{k}(\Space[\pm]{R}{2})\), \(k\geq 2\) of the discrete
   series%
   \index{representation!discrete series}~\amscite{Lang85}*{Ch.~VI}. Then we can build an operator \(F\)
   from various \(F_k\) similarly to the previous Example. In
   particular, the jump representation~\eqref{eq:jump-rl} on the real
   line generalises to the representation of a square integrable
   function \(f\) on the upper half-plane as a sum
   \begin{displaymath}
     f(z)=\sum_k a_k f_k(z), \qquad f_k\in\FSpace{B}{n}(\Space[\pm]{R}{2})
   \end{displaymath}
   for prescribed coefficients \(a_k\) and analytic functions \(f_k\) in
   question from different irreducible subspaces.

   Covariant transform is also meaningful for principal%
   \index{representation!principal series} and
   complementary%
   \index{representation!complementary series} series of representations of the group
   \(\SL\), which are not square integrable~\cite{Kisil97c}.
\end{example}
\begin{example}
  Let \(G=\mathrm{SU}(2)\times \mathrm{Aff}\) be the Cartesian product
  of the groups \(\mathrm{SU}(2)\) of unitary rotations of
  \(\Space{C}{2}\) and the \(ax+b\) group \(\mathrm{Aff}\). This group
  has a unitary linear representation on the space
  \(\FSpace{L}{2}(\Space{R}{},\Space{C}{2})\) of square-integrable (vector)
  \(\Space{C}{2}\)-valued functions by the formula:
  \begin{displaymath}
    \uir{}{}{}(g)
    \begin{pmatrix}
      f_1(t)\\f_2(t)
    \end{pmatrix}=    \begin{pmatrix}
      \alpha f_1(at+b)+ \beta f_2(at+b)\\ \mathsf{a}mma f_1(at+b)+\delta f_2(at+b)
    \end{pmatrix}, 
  \end{displaymath}
  where \(g=
  \begin{pmatrix}
    \alpha &\beta \\\mathsf{a}mma &\delta 
  \end{pmatrix}\times (a,b)\in\mathrm{SU}(2)\times \mathrm{Aff}\). It
  is obvious that the vector Hardy space, that is functions with both
  components being analytic, is invariant under such action of \(G\).

  As a fiducial operator \(F: \FSpace{L}{2}(\Space{R}{},\Space{C}{2})
  \rightarrow \Space{C}{}\) we can take, cf.~\eqref{eq:cauchy-pm}:
  \begin{equation}
    \label{eq:cauchy-pm-oper}
    F
    \begin{pmatrix}
      f_1(t)\\f_2(t)
    \end{pmatrix}
    =\frac{1}{2\pi i}\int_{\Space{R}{}} \frac{f_1(t)\,dt}{x- \rmi}.
  \end{equation}
  Thus the image of the associated covariant transform is a subspace
  of scalar valued bounded functions on \(G\). In this way we can
  transform (without a loss of information) vector-valued problems, e.g.
  matrix \emph{Wiener--Hopf factorisation}%
  \index{Wiener--Hopf factorisation}%
  \index{factorisation!Wiener--Hopf}~\cite{BoetcherKarlovichSpitkovsky02a},
  to scalar question of harmonic analysis on the group \(G\).
\end{example}
\begin{example}
  \label{it:direct-funct}
  A straightforward generalisation of Ex.~\ref{ex:wavelet} is
  obtained if \(V\) is a Banach space and \(F: V \rightarrow
  \Space{C}{}\) is an element of \(V^*\). Then the
  covariant transform coincides with the construction of wavelets in
  Banach spaces~\cite{Kisil98a}.
\end{example}
\begin{example}
  The next stage of generalisation is achieved if \(V\) is a
  Banach space and \(F: V \rightarrow \Space{C}{n}\) is a linear
  operator. Then the corresponding covariant transform is a map
  \(\oper{W}: V \rightarrow \FSpace{L}{}(G,\Space{C}{n})\). This is
  closely related to M.G.~Krein's works on \emph{directing
    functionals}%
  \index{Krein!directing functional}%
  \index{directing functional}%
  \index{functional!directing}~\cite{Krein48a}, see also \emph{multiresolution
    wavelet analysis}%
  \index{multiresolution analysis}%
  \index{analysis!multiresolution}~\cite{BratJorg97a},
  Clifford-valued%
  \index{Clifford!algebra}%
  \index{algebra!Clifford}
  Fock--Segal--Bargmann spaces%
  \index{Fock--Segal--Bargmann!space}%
  \index{space!Fock--Segal--Bargmann}~\cite{CnopsKisil97a}
  and~\amscite{AliAntGaz00}*{Thm.~7.3.1}. 
\end{example}
\begin{example}
  Let \(F\) be a projector \(\FSpace{L}{p}(\Space{R}{})\rightarrow
  \FSpace{L}{p}(\Space{R}{})\) defined by the relation \((Ff)\hat{\
  }(\lambda )=\chi(\lambda)\hat{f}(\lambda)\), where the hat denotes the Fourier
  transform and \(\chi(\lambda)\) is the characteristic function of
  the set \([-2,-1]\cup[1,2]\). Then the covariant transform
  \(\FSpace{L}{p}(\Space{R}{})\rightarrow \FSpace{C}{}(\mathrm{Aff},
  \FSpace{L}{p}(\Space{R}{}))\) generated by the
  representation~\eqref{eq:ax+b-repr-1} of the affine group from
  \(F\) contains all information provided by the \emph{Littlewood--Paley
    operator}%
  \index{Littlewood--Paley!operator}%
  \index{operator!Littlewood--Paley}~\amscite{Grafakos08}*{\S~5.1.1}.
\end{example}
\begin{example}
  \label{ex:maximal-function}
  A step in a different direction is a consideration of
  non-linear operators. Take again the ``\(ax+b\)'' group%
  \index{$ax+b$ group}%
  \index{group!$ax+b$} and its
  representation~\eqref{eq:ax+b-repr-1}. 
  We define \(F\) to be a homogeneous but non-linear functional
  \(V\rightarrow \Space[+]{R}{}\):
  \begin{displaymath}
    F (f) = \frac{1}{2}\int\limits_{-1}^1 \modulus{f(x)}\,dx.
  \end{displaymath}
  The covariant transform~\eqref{eq:coheret-transf-gen} becomes:
  \begin{equation}
    \label{eq:hardi-max}
    [\oper{W}_p f](a,b) =  F(\uir{}{p}(a,b) f) 
    = \frac{1}{2}\int\limits_{-1}^1
    \modulus{a^{\frac{1}{p}}f\left(ax+b\right)}\,dx
    = a^{\frac{1}{p}}\frac{1}{2a}\int\limits^{b+a}_{b-a}
    \modulus{f\left(x\right)}\,dx.
  \end{equation}
  Obviously \(M_f(b)=\max_{a}[\oper{W}_{\infty}f](a,b)\) coincides
  with the Hardy \emph{maximal function}%
  \index{Hardy!maximal function}%
  \index{maximal!function}%
  \index{function!maximal}, which contains important
  information on the original function \(f\). 
  From the Cor.~\ref{co:pi} we deduce that the operator \(M: f\mapsto
  M_f\) intertwines%
  \index{intertwining operator}%
  \index{operator!intertwining} \(\uir{}{p}\) with itself \(\uir{}{p}M=M
  \uir{}{p}\).

  Of course, the full covariant transform~\eqref{eq:hardi-max} is even more
  detailed than \(M\). For example,
  \(\norm{f}=\max_b[\oper{W}_{\infty}f](\frac{1}{2},b)\) is the shift
  invariant norm%
  \index{norm!shift invariant}~\cite{Johansson08a}.
\end{example}
\begin{example}
  Let \(V=\FSpace{L}{c}(\Space{R}{2})\) be the space of
  compactly supported bounded functions on the plane. We take \(F\)
  be the linear operator \(V\rightarrow \Space{C}{}\) of integration
  over the real line:
  \begin{displaymath}
    F: f(x,y)\mapsto F(f)=\int_{\Space{R}{}}f(x,0)\,dx.
  \end{displaymath}
  Let \(G\) be the group of Euclidean motions of the plane
  represented by \(\uir{}{}\) on \(V\) by a change of variables. Then
  the wavelet transform \(F(\uir{}{}(g)f)\) is the \emph{Radon
    transform}%
  \index{Radon transform}%
  \index{transform!Radon}~\cite{Helgason11a}.
\end{example}

\subsection{Symbolic Calculi}
\label{sec:symbolic-calculus}

There is a very important class of the covariant transforms which maps
operators to functions. Among numerous sources we wish to single out
works of Berezin~\cites{Berezin72,Berezin86}. We start from the
Berezin \emph{covariant symbol}.%
\index{covariant!symbol}%
\index{symbol!covariant}
\begin{example}
  Let a representation \(\uir{}{}\) of a group \(G\) act on a
  space \(X\). Then there is an associated representation
  \(\uir{}{B}\) of \(G\) on a space \(V=\FSpace{B}{}(X,Y)\) of
  linear operators \(X\rightarrow Y\) defined by the
  identity~\cites{Berezin86,Kisil98a}: 
  \begin{equation}
    \label{eq:oper-repres}
    (\uir{}{B}(g) A)x=A(\uir{}{}(g^{-1})x), \qquad x\in X,\ g\in G,\ A
    \in \FSpace{B}{}(X,Y). 
  \end{equation}
  Following the Remark~\ref{re:range-dim} we take \(F\) to be a
  functional \(V\rightarrow\Space{C}{}\), for example \(F\) can be
  defined from a pair \(x\in X\), \(l\in Y^*\)  by the expression
  \(F: A\mapsto \scalar{Ax}{l}\). Then the covariant
  transform is:
  \begin{displaymath}
    \oper{W}: A \mapsto \hat{A}(g)=F(\uir{}{B}(g) A).
  \end{displaymath}
  This is an example of \emph{covariant calculus}%
  \index{covariant!symbolic calculus}%
  \index{symbolic calculus!covariant}~\cite{Kisil98a,Berezin72}.
\end{example}
There are several variants of the last Example which are of a separate
interest.
\begin{example}
  A modification of the previous construction is obtained if we
  have two groups \(G_1\) and \(G_2\) represented by \(\uir{}{1}\)
  and \(\uir{}{2}\) on \(X\) and \(Y^*\) respectively. Then we have a covariant
  transform \(\FSpace{B}{}(X,Y)\rightarrow \FSpace{L}{}(G_1\times
  G_2, \Space{C}{})\) defined by the formula:
  \begin{displaymath}
    \oper{W}: A \mapsto \hat{A}(g_1,g_2)=\scalar{A\uir{}{1}(g_1)x}{\uir{}{2}(g_2)l}.
  \end{displaymath}
  This generalises the above \emph{Berezin covariant calculi}%
  \index{covariant!calculus!Berezin}%
  \index{calculus!covariant!Berezin}~\cite{Kisil98a}.
\end{example}
\begin{example}
  Let us restrict the previous example to the case when \(X=Y\) is a
  Hilbert space, \(\uir{}{1}{}=\uir{}{2}{}=\uir{}{}\) and \(x=l\) with
  \(\norm{x}=1\). Than the range of the covariant transform:
  \begin{displaymath}
    \oper{W}: A \mapsto \hat{A}(g)=\scalar{A\uir{}{}(g)x}{\uir{}{}(g)x}
  \end{displaymath}
  is a subset of the \emph{numerical range}%
  \index{numerical!range}%
  \index{range!numerical} of the operator \(A\). As
  a function on a group \(\hat{A}(g)\) provides a better description of
  \(A\) than the set of its values---numerical range. 
\end{example}

\begin{example}
  \label{ex:functional-model}
  The group \(\mathrm{SU}(1,1)\simeq \SL\)%
  \index{group!$\mathrm{SU}(1,1)$}%
  \index{group!$\mathrm{SU}(1,1)$|see{also $\SL$}} consists of
  \(2\times 2\) matrices of the form \(\begin{pmatrix}
    \alpha&\beta\\\bar{\beta}&\bar{\alpha}
  \end{pmatrix}\) with the unit
  determinant~\amscite{Lang85}*{\S~IX.1}.   Let \(T\) be an operator
  with the spectral radius less than \(1\). Then the associated
  M\"obius transformation%
  \index{M\"obius map}%
  \index{map!M\"obius}
  \begin{equation}
    \label{eq:meobius-T}
    g: T \mapsto g\cdot T =  \frac{\alpha T+\beta
      I}{\bar{\beta}T+\bar{\alpha}I}, \qquad \text{where} \quad
    g=
    \begin{pmatrix}
      \alpha&\beta\\\bar{\beta}&\bar{\alpha}
    \end{pmatrix}\in \SL,\ 
  \end{equation}
  produces a well-defined operator with the spectral radius less than
  \(1\) as well.  Thus we have a representation of \(SU(1,1)\).

  Let us introduce the \emph{defect operators}%
  \index{defect operator}%
  \index{operator!defect} \(D_T=(I-T^*T)^{1/2}\) and
  \(D_{T^*}=(I-TT^*)^{1/2}\). For the fiducial operator \(F=D_{T^*}\)
  the covariant transform is, cf.~\amscite{NagyFoias70}*{\S~VI.1,
    (1.2)}:
  \begin{displaymath}{}
    [\oper{W} T](g)=F(g\cdot T)=-e^{\rmi\phi}\,\Theta_T(z)\, D_T, \qquad
    \text{for } 
    g
    =     \begin{pmatrix}
      e^{\rmi\phi/2}&0\\0&e^{-\rmi\phi/2}
    \end{pmatrix}
    \begin{pmatrix}
      1&-z\\-\bar{z}&1
    \end{pmatrix},
  \end{displaymath}
  where 
  the  \emph{characteristic function}%
  \index{characteristic!function}%
  \index{function!characteristic}
  \(\Theta_T(z)\)~\amscite{NagyFoias70}*{\S~VI.1, (1.1)} is:
  \begin{displaymath}
    \Theta_T(z) = -T+D_{T^*}\,(I-zT^*)^{-1}\,z\,D_T.
  \end{displaymath}
  Thus we approached the \emph{functional model}%
  \index{functional!model}%
  \index{model!functional} of operators from the
  covariant transform. In accordance with Remark~\ref{re:range-dim}
  the model is most fruitful for the case of operator
  \(F=D_{T^*}\) being one-dimensional.
\end{example}
The intertwining property in the previous examples was obtained as a
consequence of the general Prop.~\ref{pr:inter1} about the covariant
transform. However it may be worth to select it as a separate definition:
\begin{defn}
  \label{de:covariant-calculus}
  A \emph{covariant calculus}%
  \index{covariant!calculus|indef}%
  \index{calculus!covariant|indef}, also known as
  \emph{symbolic calculus}%
  \index{symbolic calculus|see{covariant calculus}}%
  \index{calculus!symbolic|see{covariant calculus}}, is a map from operators to functions,
  which intertwines two representations of the same group in the
  respective spaces.
\end{defn}

There is a dual class of covariant transforms acting in the opposite
direction: from functions to operators. The prominent examples are the
Berezin \emph{contravariant symbol}%
\index{contravariant!symbol}%
\index{symbol!contravariant}~\cites{Berezin72,Kisil98a} and symbols 
of a \emph{pseudodifferential operators}%
\index{pseudodifferential operator}%
\index{operator!pseudodifferential}
(PDO)\index{PDO}~\cites{Howe80b,Kisil98a}. 
\begin{example}
  \label{ex:su-group}
  The classical \emph{Riesz--Dunford functional calculus}%
  \index{functional!calculus!Riesz--Dunford}%
  \index{calculus!functional!Riesz--Dunford}~\citelist{\amscite{DunfordSchwartzI}*{\S~VII.3}
    \amscite{Nikolskii86}*{\S~IV.2}} maps analytical functions on the
  unit disk to the linear operators, it is defined through the
  Cauchy-type formula with the resolvent\index{resolvent}.  The calculus is an
  intertwining operator~\cite{Kisil02a} between the M\"obius
  transformations of the unit disk, cf.~\eqref{eq:rho-1-1}, and the
  actions~\eqref{eq:meobius-T} on operators from the
  Example~\ref{ex:functional-model}.  This topic will be developed in
  Subsection~\ref{sec:anoth-appr-analyt}.
\end{example}
In line with the Defn.~\ref{de:covariant-calculus} we can
directly define the corresponding calculus through the intertwining
property~\cites{Kisil95i,Kisil02a}:
\begin{defn}
  \label{de:conravariant-calculus}
  A \emph{contravariant calculus}%
  \index{contravariant!calculus|indef}%
  \index{calculus!contravariant|indef}, also know as \emph{functional
    calculus}%
  \index{functional!calculus}%
  \index{calculus!functional}, is a map from functions to operators,
  which intertwines two representations of the same group in the
  respective spaces.
\end{defn}
The duality between co- and contravariant calculi is the particular
case of the duality between covariant transform and the inverse
covariant transform defined in the next Subsection. In many cases a
proper choice of spaces makes covariant and/or contravariant calculus
a bijection between functions and operators. Subsequently only one
form of calculus, either co- or contravariant, is defined explicitly,
although both of them are there in fact.

\subsection{Inverse Covariant Transform}
\label{sec:invar-funct-groups}
An object invariant under the left action
\(\Lambda\)~\eqref{eq:left-reg-repr} is called \emph{left invariant}.\index{invariant}
For example, let \(L\) and \(L'\) be two left invariant spaces of
functions on \(G\).  We say that a pairing \(\scalar{\cdot}{\cdot}:
L\times L' \rightarrow \Space{C}{}\) is \emph{left invariant}%
\index{invariant!pairing}%
\index{pairing!invariant} if
\begin{equation}
  \scalar{\Lambda(g)f}{\Lambda(g) f'}= \scalar{f}{f'}, \quad \textrm{ for all }
  \quad f\in L,\  f'\in L'.
\end{equation}
\begin{rem}
  \begin{enumerate}
  \item We do not require the pairing to be linear in general.
  \item If the pairing is invariant on space \(L\times L'\) it is not
    necessarily invariant (or even defined) on the whole
    \(\FSpace{C}{}(G)\times \FSpace{C}{}(G)\).
  \item In a more general setting we shall study an invariant pairing
    on a homogeneous spaces instead of the group. However due to length
    constraints we cannot consider it here beyond the Example~\ref{ex:hs-pairing}.
  \item An invariant pairing on \(G\) can be obtained from an \emph{invariant
    functional}%
  \index{invariant!functional}%
  \index{functional!invariant} \(l\) by the formula \(\scalar{f_1}{f_2}=l(f_1\bar{f}_2)\).
  \end{enumerate}
\end{rem}

For a representation \(\uir{}{}\) of \(G\) in \(V\) and \(v_0\in V\)
we fix a function \(w(g)=\uir{}{}(g)v_0\). We assume that the pairing
can be extended in its second component to this \(V\)-valued
functions, say, in the weak sense.
\begin{defn}
  \label{de:admissible}
  Let \(\scalar{\cdot}{\cdot}\) be a left invariant pairing on
  \(L\times L'\) as above, let \(\uir{}{}\) be a representation of
  \(G\) in a space \(V\), we define the function
  \(w(g)=\uir{}{}(g)v_0\) for \(v_0\in V\). The \emph{inverse
    covariant transform}%
  \index{inverse!covariant transform}%
  \index{covariant transform!inverse}%
  \index{transform!covariant!inverse} \(\oper{M}\) is a map \(L \rightarrow V\)
  defined by the pairing:
  \begin{equation}
    \label{eq:inv-cov-trans}
    \oper{M}: f \mapsto \scalar{f}{w}, \qquad \text{
      where } f\in L. 
  \end{equation}
\end{defn}

 \begin{example}
   Let \(G\) be a group with a unitary square integrable representation \(\rho\).%
  \index{square integrable!representation}%
  \index{representation!square integrable}
   An invariant pairing of two square integrable functions is obviously done by the
   integration over the Haar measure:%
   \index{Haar measure|see{invariant measure}}%
   \index{measure!Haar|see{invariant measure}}%
   \index{invariant!measure}%
   \index{measure!invariant}
   \begin{displaymath}
     \scalar{f_1}{f_2}=\int_G f_1(g)\bar{f}_2(g)\,dg.
   \end{displaymath}
   
   For an admissible vector%
   \index{admissible wavelet}%
   \index{wavelet!admissible} \(v_0\)~\citelist{\cite{DufloMoore}
   \amscite{AliAntGaz00}*{Chap.~8}} the inverse covariant transform is
   known in this setup as a \emph{reconstruction formula}%
   \index{reconstruction formula}%
   \index{formula!reconstruction}.
 \end{example}
 \begin{example}
   \label{ex:hs-pairing}
   Let \(\rho\) be a square integrable representation%
  \index{square integrable!representation}%
  \index{representation!square integrable} of \(G\) modulo a subgroup
   \(H\subset G\) and let \(X=G/H\) be the corresponding homogeneous
   space with a quasi-invariant measure \(dx\).  Then integration over
   \(dx\) with an appropriate weight produces an invariant pairing.
   The inverse covariant transform is a more general
   version~\amscite{AliAntGaz00}*{(7.52)} of the \emph{reconstruction
     formula} mentioned in the previous example.
 \end{example}
 

 Let \(\rho\) be not a square integrable representation (even modulo a subgroup) or
 let \(v_0\) be inadmissible vector of a square integrable  representation
 \(\rho\). An invariant pairing in this case is not associated with an
 integration over any non singular invariant measure on \(G\). In this
 case we have a \emph{Hardy pairing}%
\index{Hardy!pairing}%
\index{pairing!Hardy}. The following example explains
 the name.
\begin{example}
  Let \(G\) be the ``\(ax+b\)'' group%
  \index{$ax+b$ group}%
  \index{group!$ax+b$} and its representation
  \(\uir{}{}\)~\eqref{eq:ax+b-repr-1} from Ex.~\ref{ex:ax+b}. An
  invariant pairing on \(G\), which is not generated by the Haar
  measure%
   \index{invariant!measure}%
   \index{measure!invariant} \(a^{-2}da\,db\), is:
  \begin{equation}
    \label{eq:hardy-pairing}
    \scalar{f_1}{f_2}=
    \lim_{a\rightarrow 0}\int\limits_{-\infty}^{\infty}
    f_1(a,b)\,\bar{f}_2(a,b)\,db.
  \end{equation}
  For this pairing we can consider functions \(\frac{1}{2\pi i
    (x+i)}\) or \(e^{-x^2}\), which are not admissible vectors in the
  sense of square integrable representations. Then the inverse covariant transform
  provides an \emph{integral resolutions} of the identity.
\end{example}
Similar pairings can be defined for other semi-direct products of two
groups. We can also extend a Hardy pairing%
\index{Hardy!pairing}%
\index{pairing!Hardy} to a group, which has a
subgroup with such a pairing.
\begin{example}
  Let \(G\) be the group \(\SL\) from the Ex.~\ref{ex:sl2}. Then
  the ``\(ax+b\)'' group%
  \index{$ax+b$ group}%
  \index{group!$ax+b$} is a subgroup of \(\SL\), moreover we can
  parametrise \(\SL\) by triples \((a,b,\theta)\),
  \(\theta\in(-\pi,\pi]\) with the respective Haar
  measure~\amscite{Lang85}*{III.1(3)}. Then the Hardy
  pairing
  \begin{equation}
    \label{eq:hardy-pairing1}
    \scalar{f_1}{f_2}= \lim_{a\rightarrow 0}\int\limits_{-\infty}^{\infty}
    f_1(a,b,\theta)\,\bar{f}_2(a,b,\theta)\,db\,d\theta.
  \end{equation}
  is invariant on \(\SL\) as well.  The corresponding inverse
  covariant transform provides even a finer resolution of the identity
  which is invariant under conformal mappings of the Lobachevsky
  half-plane.%
  \index{Lobachevsky!geometry}%
  \index{geometry!Lobachevsky}
\end{example}

\section{Analytic Functions}
\label{sec:analytic-functions}

We saw in the first section that an inspiring geometry of cycles
can be recovered from the properties of \(\SL\).
In this section we consider a realisation of the function theory within
Erlangen approach~\cites{Kisil97c,Kisil97a,Kisil01a,Kisil02c}. The
covariant transform will be our principal tool in this construction.

\subsection{Induced Covariant Transform}
\label{sec:cauchy-transform}
The choice of a mother wavelet\index{wavelet} or fiducial operator \(F\) from
Section~\ref{sec:wavelet-transform}  can significantly influence the
behaviour of the covariant transform.  Let \(G\) be a group and
\({H}\) be its closed subgroup with the corresponding homogeneous
space \(X=G/{H}\). Let \(\uir{}{}\) be a representation of \(G\) by
operators on a space \(V\), we denote by \(\uir{}{H}\) the restriction
of \(\uir{}{}\) to the subgroup \(H\). 
\begin{defn}
  Let \(\chi\) be a representation of the subgroup \({H}\) in a space \(U\) and
  \(F: V\rightarrow U\) be an intertwining operator%
  \index{intertwining operator}%
  \index{operator!intertwining} between \(\chi\)
  and the representation \(\uir{}{H}\):
  \begin{equation}
    \label{eq:induced-mother-wavelet}
    F(\uir{}{}(h) v)=F(v)\chi(h), \qquad \text{ for all }h\in {H},\
    v\in V.
  \end{equation}
  Then the covariant transform~\eqref{eq:coheret-transf-gen} generated
  by \(F\) is called the \emph{induced covariant transform}%
  \index{induced!covariant transform}%
  \index{covariant!transform!induced}%
  \index{transform!covariant!induced}%
  \index{wavelet!transform!induced|see{induced covariant transform}}
\end{defn}
The following is the main motivating example.
\begin{example}
  Consider the traditional wavelet transform as outlined in
  Ex.~\ref{ex:wavelet}. Chose a vacuum vector \(v_0\) to be a joint
  eigenvector for all operators \(\uir{}{}(h)\), \(h\in H\), that is
  \(\uir{}{}(h) v_0=\chi(h) v_0\), where \(\chi(h)\) is a complex number
  depending of \(h\). Then \(\chi\) is obviously a character of \(H\).

  The image of wavelet transform~\eqref{eq:wavelet-transf} with such a
  mother wavelet will have a property:
  \begin{displaymath}
    \hat{v}(gh) = \scalar{ v}{\uir{}{}(gh)v_0} 
    = \scalar{v}{\uir{}{}(g)\chi(h)v_0}
    =\chi(h)\hat{v}(g).
  \end{displaymath}
  Thus the wavelet transform is uniquely defined by cosets on the
  homogeneous space \(G/H\).  In this case we previously spoke about the
  \emph{reduced wavelet transform}%
  \index{reduced wavelet transform}%
  \index{wavelet!transform!reduced}~\cite{Kisil97a}. 
  A representation \(\uir{}{0}\) is called \emph{square integrable}%
  \index{square integrable!representation}%
  \index{representation!square integrable}
  \(\mod H\) if the induced wavelet transform \([\oper{W}f_0](w)\) of
  the vacuum vector \(f_0(x)\) is square integrable on \(X\).
\end{example}
The image of induced covariant transform have the similar property:
\begin{equation}
  \label{eq:induced-covariant}
  \hat{v}(gh)=F(\uir{}{}((gh)^{-1})
  v)=F(\uir{}{}(h^{-1})\uir{}{}(g^{-1}) v)
  =F(\uir{}{}(g^{-1}) v)\chi{}{}(h^{-1}).
\end{equation}
Thus it is enough to know the value of the covariant transform only at a
single element in every coset \(G/H\) in order to reconstruct it for
the entire group \(G\) by the representation \(\chi\). Since coherent
states (wavelets) are now parametrised by points
homogeneous space \(G/H\) they are referred sometimes as coherent
states which are not connected to a group~\cite{Klauder96a}, however
this is true only in a very narrow sense as explained above.
\begin{example}
  To make it more specific we can consider the representation of \(\SL\)
  defined on \(\FSpace{L}{2}(\Space{R}{})\) by the formula, cf.~\eqref{eq:discrete}:
  \begin{displaymath}
    \uir{}{}(g): f(z) \mapsto \frac{1}{(c x + d)}\,
    f\left(\frac{a x+ b }{c x +d}\right), \qquad g^{-1}=\
    \begin{pmatrix}
      a & b \\ c & d 
    \end{pmatrix}.
  \end{displaymath}
  Let \(K\subset\SL\) be the compact subgroup 
  of matrices \(
  h_t=
  \begin{pmatrix}
    \cos t&\sin t\\-\sin t&\cos t
  \end{pmatrix}\). Then for the fiducial operator
  \(F_{\pm}\)~\eqref{eq:cauchy-pm} we have
  \(F_{\pm}\circ\uir{}{}(h_t)=e^{\mp\rmi t}F_{\pm}\). Thus we can
  consider the covariant transform only for points in the homogeneous
  space \(\SL/K\), moreover this set can be naturally identified with
  the \(ax+b\) group.%
  \index{$ax+b$ group}%
  \index{group!$ax+b$} Thus we do not obtain any advantage of
  extending the group in Ex.~\ref{ex:ax+b} from \(ax+b\) to \(\SL\) if
  we will be still using the fiducial operator
  \(F_\pm\)~\eqref{eq:cauchy-pm}.
\end{example}
Functions on the group \(G\), which have the property
\(\hat{v}(gh)=\hat{v}(g)\chi(h)\)~\eqref{eq:induced-covariant},
provide a space for the representation of \(G\) induced by the
representation \(\chi\) of the subgroup \(H\). This explains the
choice of the name for induced covariant transform.

\begin{rem}
  Induced covariant transform uses the fiducial operator \(F\) which
  passes through the action of the subgroup \({H}\). This reduces
  information which we obtained from this transform in some cases.
\end{rem}

There is also a simple connection between a covariant transform and
right shifts: 
\begin{prop}
  \label{pr:inducer-wave-intertw}
  Let \(G\) be a Lie group and \(\uir{}{}\) be a representation of
  \(G\) in a space \(V\). Let \([\oper{W}f](g)=F(\uir{}{}(g^{-1})f)\) be a
  covariant transform defined by the fiducial operator \(F: V \rightarrow U\).
  Then the right shift \([\oper{W}f](gg')\) by \(g'\) is the covariant transform
  \([\oper{W'}f](g)=F'(\uir{}{}(g^{-1})f)]\) defined by the fiducial operator
  \(F'=F\circ\uir{}{}(g^{-1})\). 

  In other words the covariant transform intertwines%
  \index{intertwining operator}%
  \index{operator!intertwining} right shifts on
  the group \(G\) with the associated action
  \(\uir{}{B}\)~\eqref{eq:oper-repres} on fiducial operators.
\end{prop}
Although the above result is obvious, its infinitesimal version has
interesting consequences.
\begin{cor}[\cite{Kisil10c}]
  \label{co:cauchy-riemann}
  Let \(G\) be a Lie group with a Lie algebra \(\algebra{g}\) and
  \(\uir{}{}\) be a smooth representation of \(G\). We denote by
  \(d\uir{}{B}\) the derived representation of the associated
  representation \(\uir{}{B}\)~\eqref{eq:oper-repres} on fiducial
  operators.

  Let a fiducial operator \(F\) be a null-solution, i.e. \(A F=0\),
  for the operator \(A=\sum_J a_j d\uir{X_j}{B}\), where
  \(X_j\in\algebra{g}\) and \(a_j\) are constants.  Then the covariant
  transform \([\oper{W} f](g)=F(\uir{}{}(g^{-1})f)\) for any \(f\)
  satisfies:
  \begin{displaymath}
    D F(g)= 0, \qquad \text{where} \quad
    D=\sum_j \bar{a}_j\linv{X_j}.
  \end{displaymath}
  Here \(\linv{X_j}\) are the left invariant fields (Lie derivatives) on
  \(G\) corresponding to \(X_j\).
\end{cor}
\begin{example}
  Consider the representation \(\uir{}{}\)~\eqref{eq:ax+b-repr-1} of the \(ax+b\)
  group with the \(p=1\). Let \(A\) and \(N\) be the basis of the
  corresponding Lie algebra generating one-parameter subgroups
  \((e^t,0)\) and \((0,t)\). Then the derived representations are:
  \begin{displaymath}
    [d\uir{A}{} f](x)= f(x)+xf'(x), \qquad [d\uir{N}{}f](x)=f'(x).
  \end{displaymath}
  The corresponding left invariant vector fields on \(ax+b\) group%
  \index{$ax+b$ group}%
  \index{group!$ax+b$} are:
  \begin{displaymath}
   \linv{A} =a \partial_a,\qquad \linv{N}=a\partial_b.
  \end{displaymath}
  The mother wavelet \(\frac{1}{x+\rmi}\) is a null solution of the
  operator \(d\uir{A}{} +\rmi d\uir{N}{}=I+(x+\rmi)\frac{d}{dx}\).
  Therefore the covariant transform with the fiducial operator
  \(F_+\)~\eqref{eq:cauchy-pm} will consist with the null solutions to
  the operator \(\linv{A}-\rmi\linv{N}=-\rmi a(\partial_b+\rmi\partial_a)\),
  that is in the essence the Cauchy-Riemann operator%
  \index{Cauchy-Riemann operator}%
  \index{operator!Cauchy-Riemann} in the upper
  half-plane. 
\end{example}
There is a statement which extends the previous Corollary from
differential operators to integro-differential ones. We will formulate
it for the wavelets setting.
\begin{cor}
  \label{co:cauchy-riemann-integ}
  Let \(G\) be a group and \(\uir{}{}\) be a unitary representation
  of \(G\), which can be extended to a vector space \(V\) of functions
  or distributions on \(G\).
   Let a mother wavelet \(w\in V'\)  satisfy the equation 
  \begin{displaymath}
    \int_{G} a(g)\, \uir{}{}(g) w\,dg=0,
  \end{displaymath}
  for a fixed distribution \(a(g) \in V\) and a (not necessarily
  invariant) measure \(dg\). Then  any wavelet transform
  \(F(g)= \oper{W} f(g)=\scalar{f}{\uir{}{}(g)w_0}\) obeys the condition:
  \begin{displaymath}
   DF=0,\qquad \text{where} \quad D=\int_{G} \bar{a}(g)\, R(g) \,dg,
  \end{displaymath}
  with \(R\) being the right regular representation of \(G\).
\end{cor}
Clearly, the Corollary~\ref{co:cauchy-riemann} is a particular case of
the Corollary~\ref{co:cauchy-riemann-integ} with a distribution \(a\),
which is a combination of derivatives of Dirac's delta functions. The
last Corollary will be illustrated at the end of
Section~\ref{sec:anoth-appr-analyt}. 
\begin{rem}
  \label{re:analyt-cor-expan}
  We note that Corollaries~\ref{co:cauchy-riemann}
  and~\ref{co:cauchy-riemann-integ} are true whenever we have an
  intertwining property between  \(\uir{}{}\) with the right regular
  representation of \(G\). 
\end{rem}

\subsection{Induced Wavelet Transform and Cauchy Integral}
\label{sec:wavel-transf-cauchy}

We again use the general scheme from
Subsection~\ref{sec:concl-induc-repr}.  The \(ax+b\) group%
\index{$ax+b$ group}%
\index{group!$ax+b$} is
isomorphic to a subgroups of \(\SL\) consisting of the
lower-triangular matrices:
\begin{displaymath}
  F=\left\{
      \frac{1}{\sqrt{a}}\begin{pmatrix}
        a &0\\b&1
      \end{pmatrix},\ a>0\right\}.
\end{displaymath}
The corresponding homogeneous space \(X=\SL/F\) is one-dimensional and
can be parametrised by a real number. The natural projection
\(p:\SL\rightarrow \Space{R}{}\) and its
left inverse \(s: \Space{R}{}\rightarrow \SL\) can
be defined as follows: 
\begin{equation}
  \label{eq:ps-maps-rl}
  p:
  \begin{pmatrix}
    a & b \\ c & d
  \end{pmatrix} \mapsto \frac{b}{d},\qquad  s: u \mapsto
  \begin{pmatrix}
    1 & u \\ 0 & 1
  \end{pmatrix}.
\end{equation}
Thus we  calculate the corresponding map \(r: \SL\rightarrow F\), see
Subsection~\ref{sec:hypercomplex-numbers}:
\begin{equation}
  \label{eq:r-map-rl}
  r:
  \begin{pmatrix}
    a & b \\ c & d
  \end{pmatrix} \mapsto \begin{pmatrix}
    d^{-1} & 0 \\ c & d
  \end{pmatrix}.
\end{equation}
Therefore the action of \(\SL\) on the real line is exactly the M\"obius
map%
\index{M\"obius map}%
\index{map!M\"obius}~\eqref{eq:moebius}:
\begin{displaymath}
  g:u\mapsto p(g^{-1}*s(u)) =\frac{au+b}{cu+d}, \quad \text{ where }
  g^{-1}=  \begin{pmatrix}
    a&b\\c&d
  \end{pmatrix}.
\end{displaymath}
We also calculate that
\begin{displaymath}
  r(g^{-1}*s(u)) =
  \begin{pmatrix}
    (cu+d)^{-1}&0\\
    c&cu+d
  \end{pmatrix}.
\end{displaymath}
 
To build an induced representation we need a character of the affine group. 
A generic character of \(F\) is a power of its diagonal element:
\begin{displaymath}
\uir{}{\kappa}
\begin{pmatrix}
  a&0\\c&a^{-1}
\end{pmatrix}=a^\kappa.  
\end{displaymath}

Thus the corresponding realisation of induced
representation~\eqref{eq:def-ind} is:
\begin{equation}
  \label{eq:induced-affine}
  \uir{}{\kappa }(g): f(u) \mapsto \frac{1}{(cu+d)^\kappa} \,
  f\left(\frac{au+b}{cu+d}\right)\qquad
  \text{ where }g^{-1}=  \begin{pmatrix}
    a&b\\c&d
  \end{pmatrix}.
\end{equation}
The only freedom remaining by the scheme is in a choice of a
value of \emph{number} \(\kappa\) and the corresponding functional
space where our representation acts. At this point we have a wider
choice of \(\kappa\) than it is usually assumed: it can belong to
different hypercomplex systems.

One of the important properties which would be nice to have is the
unitarity of the representation~\eqref{eq:induced-affine} with respect
to the standard inner product:
\begin{displaymath}
  \scalar{f_1}{f_2}=\int_{\Space{R}{2}} f_1(u)\bar{f}_2(u)\,du.
\end{displaymath}
A change of variables \(x=\frac{au+b}{cu+d}\) in the integral suggests
the following property is necessary and sufficient for that: 
\begin{equation}
  \label{eq:kappa-unitarity}
  \kappa+\bar{\kappa}=2.
\end{equation}

A mother wavelet\index{wavelet} for an induced wavelet transform shall be an
eigenvector for the action of a subgroup \(\tilde{H}\) of \(\SL\),
see~\eqref{eq:induced-mother-wavelet}. Let us consider the most common
case of \(\tilde{H}=K\) and take the infinitesimal condition with the
derived representation: \(d\uir{Z}{n}w_0 =\lambda w_0\), since
\(Z\)~\eqref{eq:sl2-basis} is the generator of the subgroup \(K\). In
other word the restriction of \(w_0\) to a \(K\)-orbit should be given
by \(e^{\lambda t}\) in the exponential coordinate \(t\) along the
\(K\)-orbit. However we usually need its expression in other ``more
natural'' coordinates.  For example~\cite{Kisil11b}, an eigenvector of
the derived representation of \(d\uir{Z}{n}\) should satisfy the
differential equation in the ordinary parameter \(x\in\Space{R}{}\):
\begin{equation}
  \label{eq:rl-k-deriv}
  -\kappa xf(x)-f'(x)(1+x^2)=\lambda f(x).
\end{equation}
The equation does not have singular points, the general solution
is globally defined (up to a constant factor) by:
\begin{equation}
  \label{eq:rl-k-eigenvector}
  w_{\lambda, \kappa }(x)= \frac{1}{(1+x^2)^{\kappa /2}}\left(\frac{x-\rmi }{x
      +\rmi}\right)^{\rmi\lambda/2}
  =\frac{(x-\rmi)^{(\rmi\lambda-\kappa )/2}
  }{(x+\rmi)^{(\rmi\lambda+\kappa )/2}}.
\end{equation}
To avoid multivalent functions we need \(2\pi\)-periodicity along the
exponential coordinate on \(K\).  This implies that the parameter
\(m=-\rmi\lambda\) is an integer.
Therefore the solution becomes:
\begin{equation}
  \label{eq:rl-k-eigenvector-m}
   w_{m,\kappa }(x) =\frac{(x+\rmi)^{(m-\kappa )/2}}{(x-\rmi)^{(m+\kappa )/2} }.
\end{equation}
The corresponding wavelets resemble the Cauchy kernel normalised to the
invariant metric in the Lobachevsky half-plane:%
\index{Lobachevsky!geometry}%
\index{geometry!Lobachevsky}
\begin{eqnarray*}
  w_{m,\kappa} (u,v;x) &=&  \uir{F}{\kappa }(s(u,v)) w_{m,\kappa }(x)
  =
  v^{\kappa/ 2}   \frac{\left(x -u +\rmi v\right)^{(m-\kappa)/2}}
  {\left(x- u -\rmi v\right)^{(m+\kappa)/2 }}
\end{eqnarray*}
Therefore the wavelet transform~\eqref{eq:wavelet-transf} from
function on the real line to functions on the upper half-plane is:
\begin{eqnarray*}
  \hat{f}(u,v)&=&\scalar{f}{\uir{F}{\kappa }(u,v)w_{m,\kappa }}
    =v^{\bar{\kappa} /2}\int_{\Space{R}{}} f(x)\, \frac{(x-(u+\rmi v))^{(m-\kappa )/2}}
  {(x-(u-\rmi  v))^{(m+\kappa)/2}}\, dx.
\end{eqnarray*}
Introduction of a complex variable \(z=u+\rmi v\) allows to write it
as:
\begin{equation}
  \label{eq:wavelet-trans-ell-hp}
  \hat{f}(z)=(\Im z)^{\bar{\kappa}/2} \int_{\Space{R}{}}f(x)  
  \frac{(x-{z})^{(m-\kappa)/2}}{(x-\bar{z})^{(m+\kappa)/2}}\, dx.
\end{equation}
According to the general theory this wavelet transform intertwines
representations \(\uir{F}{\kappa}\)~\eqref{eq:induced-affine} on the
real line (induced by the character \(a^{\kappa}\) of the subgroup
\(F\)) and \(\uir{K}{m}\)~\eqref{eq:discrete} on the upper half-plane
(induced by the character \(e^{\rmi m t}\) of the subgroup \(K\)).

\subsection{The Cauchy-Riemann (Dirac) and Laplace Operators}
\label{sec:dirac-cauchy-riemann}

Ladder operators \(\ladder{\pm}=\pm\rmi A +B\)%
\index{ladder operator}%
\index{operator!ladder} act by raising/lowering indexes of the
\(K\)-eigenfunctions \(w_{m,\kappa}\)~\eqref{eq:rl-k-eigenvector}, see
Subsection~\ref{sec:correspondence}. More explicitly~\cite{Kisil11b}:
\begin{equation}
  \label{eq:ladder-action-ell}
  d\uir{\ladder{\pm}}{\kappa }: w_{m,\kappa } 
  \mapsto -\frac{\rmi}{2}( m \pm \kappa) w_{m\pm 2,\kappa}.
\end{equation}
There are two possibilities here: \(m\pm\kappa\) is zero for some
\(m\) or not. In the first case the chain~\eqref{eq:ladder-action-ell}
of eigenfunction \(w_{m,\kappa}\) terminates on one side under the
transitive action~\eqref{eq:ladder-chain-1D} of the ladder operators;
otherwise the chain is infinite in both directions.  That is, the values
\(m=\mp\kappa \) and only those correspond to the maximal (minimal)
weight function \(w_{\mp\kappa ,\kappa }(x)=\frac{1}{(x\pm\rmi)^\kappa
} \in \FSpace{L}{2}(\Space{R}{})\), which are annihilated by
\(\ladder{\pm}\):
\begin{eqnarray}
  \label{eq:annihilated}
  d\uir{\ladder{\pm}}{\kappa } w_{\mp\kappa ,\kappa }= (\pm\rmi d\uir{A}{\kappa}
  +d\uir{B}{\kappa})\,w_{\mp\kappa ,\kappa } =0.
\end{eqnarray}

 By the Cor.~\ref{co:cauchy-riemann} for the mother
wavelets \(w_{\mp\kappa ,\kappa }\), which are annihilated
by~\eqref{eq:annihilated}, the images of the respective wavelet
transforms are null solutions to the left-invariant differential
operator \(D_{\pm}=\overline{\linv{{\ladder{\pm}}}}\):%
\index{ladder operator}%
\index{operator!ladder}
\begin{equation}
  \label{eq:cauchy-riemann-ell-conf}
  D_{\pm}=\mp\rmi\linv{A}+\linv{B}=\textstyle -\frac{\rmi \kappa 
  }{2}+v(\partial_u\pm\rmi\partial_v). 
\end{equation}

This is a conformal version of the Cauchy--Riemann equation%
\index{Cauchy-Riemann operator}%
\index{operator!Cauchy-Riemann}. The
second order conformal Laplace-type operators%
\index{Laplacian}
\(\Delta_+=\overline{\linv{\ladder{-}}\linv{\ladder{+}}}\) and
\(\Delta_-=\overline{\linv{\ladder{+}}\linv{\ladder{-}}}\) are:
\begin{equation}
  \label{eq:laplace-ell-conf}
  \Delta_\pm =   \textstyle(v\partial_u-\frac{\rmi \kappa  }{2})^2+v^{2}\partial_v^2
  \pm\frac{\kappa  }{2}. 
\end{equation}

For the mother wavelets \(w_{m,\kappa}\) in~\eqref{eq:annihilated}
such that \(m=\mp\kappa \) the unitarity condition
\(\kappa+\bar{\kappa}=2\), see~\eqref{eq:kappa-unitarity}, together
with \(m\in\Space{Z}{}\) implies \(\kappa=\mp m=1\). In such a case
the wavelet transforms~\eqref{eq:wavelet-trans-ell-hp} are:
\begin{equation}
  \label{eq:cauchy-int-ell}
  \hat{f}^+(z)=(\Im z)^{\frac{1}{2}} \int_{\Space{R}{}}  
  \frac{f(x)\, dx}{x-z}\quad\text{and}
  \quad   \hat{f}^-(z)=(\Im z)^{\frac{1}{2}} \int_{\Space{R}{}}  
  \frac{f(x)\, dx}{x-\bar{z}},
\end{equation}
for \(w_{-1,1}\) and \(w_{1,1}\) respectively. The first one is the Cauchy integral
formula%
  \index{integral!Cauchy}%
  \index{Cauchy!integral} up to the factor \(2\pi\rmi \sqrt{\Im z}\). Clearly, one integral is
the complex conjugation of another. Moreover, the minimal/maximal
weight cases can be intertwined by the following automorphism of the
Lie algebra \(\algebra{sl}_2\):
\begin{displaymath}
  A\rightarrow B, \quad B\rightarrow A,\quad Z\rightarrow -Z.
\end{displaymath}

As explained before
\(\hat{f}^\pm(w)\) are null solutions to the operators
\(D_\pm\)~\eqref{eq:cauchy-riemann-ell-conf}  and
\(\Delta_{\pm}\)~\eqref{eq:laplace-ell-conf}.  These transformations
intertwine unitary equivalent representations on the real line and
on the upper half-plane, thus they can be made unitary for proper
spaces. This is the source of two faces of the Hardy spaces:%
\index{space!Hardy}%
\index{Hardy!space} they can
be defined either as square-integrable on the real line with an analytic extension to
the half-plane, or analytic on the half-plane with
square-integrability on an infinitesimal displacement of the real line.

For the third possibility, \(m\pm\kappa\neq 0\), there is no an
operator spanned by the derived representation of the Lie algebra
\(\algebra{sl}_2\) which kills the mother wavelet \(w_{m,\kappa}\).
However the remarkable \emph{Casimir operator}%
\index{Casimir operator}%
\index{operator!Casimir}
\(C=Z^2-2(\ladder{-}\ladder{+}+\ladder{+}\ladder{-})\), which spans the
centre of the universal enveloping algebra of \(\algebra{sl}_2\)
\citelist{\amscite{MTaylor86}*{\S~8.1} \amscite{Lang85}*{\S~X.1}},
produces a second order operator which does the job. Indeed from the
identities~\eqref{eq:ladder-action-ell} we get:
\begin{equation}
  \label{eq:casimir-action-scalar-ell}
  d\uir{C}{\kappa} w_{m,\kappa } = ( 2\kappa  - \kappa^2) w_{m,\kappa}.
\end{equation}
Thus we get \(d\uir{C}{\kappa} w_{m,2}=0\) for \(\kappa=2\)  or \(0\).
The mother wavelet \(w_{0,2}\) turns to be the \emph{Poisson
  kernel}%
  \index{Poisson!kernel}%
  \index{kernel!Poisson}~\amscite{Grafakos08}*{Ex.~1.2.17}. 
The associated wavelet transform
\begin{equation}
  \label{eq:ell-poisson-int}
    \hat{f}(w)=\Im z \int_{\Space{R}{}}  
      \frac{f(x)\, dz}{\modulus{x-z}^2}
\end{equation}
consists of null solutions of the left-invariant second-order
Laplacian%
\index{Laplacian}%
\index{Laplace operator|see{Laplacian}}%
\index{operator!Laplace|see{Laplacian}}, image of the Casimir operator,
cf.~\eqref{eq:laplace-ell-conf}: 
\begin{displaymath}
    \Delta(:=\linv{C})    = v^2\partial^2_u+v^{2}\partial_v^2.
\end{displaymath}
Another integral formula producing solutions to this equation
delivered by the mother wavelet \(w_{m,0}\) with the value \(\kappa=0\)
in~\eqref{eq:casimir-action-scalar-ell}:
\begin{equation}
  \label{eq:unknown-harmonic-ell}
  \hat{f}(z)= \int_{\Space{R}{}}f(x)  
  \left(\frac{x-{z}}{x-\bar{z}}\right)^{m/2}\, dx.
\end{equation}

Furthermore, we can introduce higher order differential operators.
The functions \(w_{\mp 2m+1,1}\) are annihilated by \(n\)-th power of
operator \(d\uir{\ladder{\pm}}{\kappa}\) with \(1\leq m\leq n\). By
the Cor.~\ref{co:cauchy-riemann} the the image of wavelet
transform~\eqref{eq:wavelet-trans-ell-hp} from a mother wavelet
\(\sum_1^n a_m w_{\mp 2m,1}\) will consist of null-solutions of the
\(n\)-th power \(D^n_\pm\) of the conformal Cauchy--Riemann operator%
\index{Cauchy-Riemann operator}%
\index{operator!Cauchy-Riemann}~\eqref{eq:cauchy-riemann-ell-conf}.
They are a conformal flavour of \emph{polyanalytic}%
\index{polyanalytic function}%
\index{function!polyanalytic}
functions~\cite{Balk97a}.

We can similarly look for mother wavelets which are eigenvectors for other
types of one dimensional subgroups. Our consideration of subgroup
\(K\) is simplified by several facts: 
\begin{itemize}
\item The parameter \(\kappa\) takes only complex values.
\item The derived representation does not have singular points on the real line.
\end{itemize}
For both subgroups \(\Aprime\) and \(N'\) this will not be true. The
further consideration will be given in~\cite{Kisil11b}.

\subsection{The Taylor Expansion}
\label{sec:taylor-expansion}

Consider an induced wavelet transform generated by a Lie group \(G\),
its representation \(\uir{}{}\) and a mother wavelet \(w\) which is an
eigenvector of a one-dimensional subgroup \(\tilde{H}\subset G\). Then
by Prop.~\ref{pr:inducer-wave-intertw} the wavelet transform
intertwines \(\uir{}{}\) with a representation \(\uir{\tilde{H}}{}\)
induced by a character of \(\tilde{H}\).

If the mother wavelet is itself in the domain of the induced wavelet
transform then the chain~\eqref{eq:ladder-chain-1D} of
\(\tilde{H}\)-eigenvectors \(w_m\) will be mapped to the similar chain
of their images \(\hat{w}_m\). The corresponding derived induced
representation \(d\uir{\tilde{H}}{}\) produces ladder operators%
\index{ladder operator}%
\index{operator!ladder} with
the transitive action of the ladder operators on the chain of
\(\hat{w}_m\). Then the vector space of ``formal power series'':
\begin{equation}
  \label{eq:taylor-prototype}
  \hat{f}(z)=\sum_{m\in \Space{Z}{}} a_m \hat{w}_{m}(z)
\end{equation}
is a module for the Lie algebra of the group \(G\).

Coming back to the case of the group \(G=\SL\) and subgroup
\(\tilde{H}=K\).  Images \(\hat{w}_{m,1}\) of the
eigenfunctions~\eqref{eq:rl-k-eigenvector-m} under the Cauchy integral
transform%
  \index{integral!Cauchy}%
  \index{Cauchy!integral}~\eqref{eq:cauchy-int-ell} are:
\begin{displaymath}
  \hat{w}_{m,1}(z)=(\Im z)^{1/2} \frac{(z+\rmi)^{(m-1)/2}}{(z-\rmi)^{(m+1)/2}}.
\end{displaymath}
They are eigenfunctions of the derived representation on the upper
half-plane and the action of ladder operators is
given by the same expressions~\eqref{eq:ladder-action-ell}.
In particular, the \(\algebra{sl}_2\)-module generated by
\(\hat{w}_{1,1}\) will be one-sided since this vector is annihilated by
the lowering operator. 
Since the Cauchy integral%
  \index{integral!Cauchy}%
  \index{Cauchy!integral} produces an unitary intertwining
operator between two representations we get the following variant of
Taylor series:
\begin{displaymath}
  \hat{f}(z)=\sum_{m=0}^\infty c_m \hat{w}_{m,1}(z), \qquad
  \text{ where } \quad c_m=\scalar{f}{w_{m,1}}.
\end{displaymath}
For two other types of subgroups, representations and mother wavelets this
scheme shall be suitably adapted and detailed study will be presented
elsewhere~\cite{Kisil11b}.

\subsection{Wavelet Transform in the Unit Disk and Other Domains}
\label{sec:backgr-compl-analys}
We can similarly construct an analytic function theories in unit
disks, including parabolic and hyperbolic ones~\cite{Kisil05a}. This
can be done simply by an application of the \emph{Cayley transform}%
\index{Cayley transform}%
\index{transform!Cayley} to the
function theories in the upper half-plane. Alternatively we can apply
the full procedure for properly chosen groups and subgroups. We will
briefly outline such a possibility here, see also~\cite{Kisil97c}.

Elements of \(\SL\) could be also represented by \(2\times 2\)-matrices
with complex entries such that, cf.~Example~\ref{ex:su-group}:
\begin{displaymath}
  g= \matr{\alpha}{\bar{\beta}}{\beta}{\bar{\alpha}},
  \qquad 
  g^{-1}= \matr{\bar{\alpha}}{-\bar{\beta}}{-\beta}{\alpha}, 
  \qquad
  \modulus{\alpha}^2-\modulus{\beta}^2=1.
\end{displaymath}
This realisations of \(\SL\) (or rather \(SU(2,\Space{C}{})\)) is more
suitable for function theory in the unit disk. It is obtained from the
form, which we used before for the upper half-plane, by means of the
Cayley transform~\amscite{Kisil05a}*{\S~8.1}.

We may identify the \emph{unit disk}%
\index{unit!disk} \(\Space{D}{}\) with the homogeneous space
\(\SL/\Space{T}{}\) for the \emph{unit circle}%
\index{unit!circle} \(\Space{T}{}\) through
the important decomposition \(\SL\sim \Space{D}{}\times\Space{T}{}\)\
with \(K=\Space{T}{}\)---the compact subgroup of \(\SL\):
\begin{eqnarray}
\label{eq:sl2-u-psi-coord}
  \matr{\alpha}{\bar{\beta}}{\beta}{\bar{\alpha}} 
  & =& \modulus{\alpha} 
  \matr{1}{\bar{\beta}\bar{\alpha}^{-1}}{{\beta}{\alpha}^{-1}}{1}
  \matr{ 
    \frac{{\alpha}}{ \modulus{\alpha} } }{0}{0}{\frac{\bar{\alpha}}{ 
      \modulus{\alpha} } }\\
  &=& \frac{1}{\sqrt{1- \modulus{u}^2 }}
  \matr{1}{u}{\bar{u}}{1}, \nonumber 
  \matr{e^{ix}}{0}{0}{e^{-ix}}
  \end{eqnarray}
where
\begin{displaymath}
  x=\arg \alpha,\qquad 
  u=\bar{\beta}\bar{\alpha}^{-1},\qquad \modulus{u}<1.
\end{displaymath}
Each element \(g\in\SL\) acts by the linear-fractional transformation
(the M\"obius map) on \(\Space{D}{}\)\ 
and \(\Space{T}{}\)
\(\FSpace{H}{2}(\Space{T}{})\) as follows: 
\begin{equation}
  \label{eq:moebius-su}
  g: z \mapsto \frac{\alpha  z +\beta }{{\bar{\beta}} z+\bar{\alpha}},
\qquad \textrm{ where } \quad
g=\matr{\alpha}{\beta}{\bar{\beta}}{\bar{\alpha}}.
\end{equation}
In the decomposition~\eqref{eq:sl2-u-psi-coord} the first matrix on
the right hand side acts by transformation~\eqref{eq:moebius-su} as an
orthogonal rotation of \(\Space{T}{}\) or \(\Space{D}{}\); and the
second one---by transitive family of maps of the unit disk onto
itself.

The representation  induced by a complex-valued character 
\(\chi_k(z)=z^{-k}\) of \(\Space{T}{}\) according to the
Section~\ref{sec:concl-induc-repr} is:
\begin{equation}
  \label{eq:rho-1-1}
  \rho_k(g): f(z) \mapsto
  \frac{1}{(\alpha-{\beta}{z})^k} \,
  f\left(
    \frac{\bar{\alpha} z - \bar{\beta}}{\alpha-{\beta} z} 
  \right)
  \qquad  \textrm{ where } \quad 
  g=\matr{\alpha}{\beta}{\bar{\beta}}{\bar{\alpha}}.
\end{equation}
The representation \(\uir{}{1}\) is unitary on square-integrable
functions and irreducible on the \emph{Hardy space}%
\index{space!Hardy}%
\index{Hardy!space} on the unit circle.
 
We choose~\cites{Kisil98a,Kisil01a} \(K\)-invariant function \(v_0(z)\equiv 1\) 
to be  a \emph{vacuum vector}.
Thus the associated \emph{coherent states}
\begin{displaymath}
  v(g,z)=\rho_1(g)v_0(z)= (u-z)^{-1}
\end{displaymath} 
are completely determined by the point on the unit disk \(
u=\bar{\beta}\bar{\alpha}^{-1}\). The family of coherent states considered
as a function of both \(u\) and \(z\) is obviously the \emph{Cauchy
  kernel}~\cite{Kisil97c}. The \emph{wavelet transform}~\cite{Kisil97c,Kisil98a}
\(\oper{W}:\FSpace{L}{2}(\Space{T}{})\rightarrow
\FSpace{H}{2}(\Space{D}{}): f(z)\mapsto
\oper{W}f(g)=\scalar{f}{v_g}\)\ is the \emph{Cauchy integral}%
  \index{integral!Cauchy}%
  \index{Cauchy!integral}:
\begin{equation}
  \label{eq:cauchy}
  \oper{W} f(u)=\frac{1}{2\pi i}\int_{\Space{T}{}}f(z)\frac{1}{u-z}\,dz.
\end{equation}

This approach can be extended to arbitrary connected simply-connected
domain. Indeed, it is known that M\"obius maps%
\index{M\"obius map}%
\index{map!M\"obius} is the whole group of
biholomorphic automorphisms of the unit disk or upper half-plane. Thus
we can state the following corollary from the \emph{Riemann mapping
  theorem}%
\index{Riemann!mapping theorem}%
\index{theorem!Riemann mapping}:
\begin{cor}
  \label{co:riemann-mapping}
  The group of biholomorphic automorphisms of a connected
  simply connected domain%
  \index{simply connected domain}%
  \index{domain!simply connected} with at least two points on its boundary is
  isomorphic to \(\SL\).
\end{cor}

If a domain is non-simply connected%
  \index{non-simply connected domain}%
  \index{domain!non-simply connected}, then the group of its
biholomorphic mapping can be
trivial~\cites{MityushevRogosin00a,Beardon07a}.  However we may look
for a rich group acting on function spaces rather than on geometric sets.
Let a connected non-simply connected domain \(D\) be
bounded by a finite collection of non-intersecting contours
\(\mathsf{a}mma_i\), \(i=1,\ldots,n\). For each \(\mathsf{a}mma_i\) consider the
isomorphic image \(G_i\) of the \(\SL\) group which is defined by the
Corollary~\ref{co:riemann-mapping}.  Then define the group
\(G=G_1\times G_2\times \ldots \times G_n\) and its action on
\(\FSpace{L}{2}(\partial D)= \FSpace{L}{2}(\mathsf{a}mma_1)\oplus
\FSpace{L}{2}(\mathsf{a}mma_2)\oplus \ldots \oplus \FSpace{L}{2}(\mathsf{a}mma_n)\)
through the Moebius action of \(G_i\) on \(\FSpace{L}{2}(\mathsf{a}mma_i)\).
\begin{example}
  Consider an \emph{annulus}\index{annulus} defined by
  \(r<\modulus{z}<R\). It is bounded by two circles: \(\mathsf{a}mma_1=\{z:
  \modulus{z}=r\}\) and \(\mathsf{a}mma_2=\{z:
  \modulus{z}=R\}\). For \(\mathsf{a}mma_1\) the M\"obius action of \(\SL\)
  is
  \begin{displaymath}
    \begin{pmatrix}
      \alpha&\bar\beta\\
      \beta &\bar\alpha
    \end{pmatrix}: z\mapsto
    \frac{\alpha z +\bar\beta/r}{\beta z/r + \bar\alpha},\qquad
    \text{where} \quad \modulus{\alpha}^2-\modulus{\beta}^2=1,
  \end{displaymath}
  with the respective action on \(\mathsf{a}mma_2\). Those action can be
  linearised in the spaces \(\FSpace{L}{2}(\mathsf{a}mma_{1})\) and
  \(\FSpace{L}{2}(\mathsf{a}mma_{2})\). If we consider a subrepresentation
  reduced to analytic function on the annulus, then one copy of
  \(\SL\) will act on the part of functions analytic outside of
  \(\mathsf{a}mma_1\) and another copy---on the part of functions analytic
  inside of \(\mathsf{a}mma_2\).
\end{example}

Thus all classical objects of complex analysis (the Cauchy-Riemann
equation, the Taylor series, the Bergman space%
\index{space!Bergman}%
\index{Bergman!space}, etc.) for a rather generic domain \(D\) can
be also obtained from suitable representations similarly
to the case of the upper half-plane~\cite{Kisil97c,Kisil01a}.

\section{Covariant and Contravariant Calculi}
\label{sec:functional-calculus}

United in the trinity functional calculus%
\index{functional!calculus}%
\index{calculus!functional}, spectrum, and spectral
mapping theorem play the exceptional r\^ole in functional
analysis and could not be substituted by anything else. 
Many traditional  definitions of functional calculus  are covered by the
following rigid template based on the \emph{algebra homomorphism}%
\index{algebra!homomorphism}%
\index{homomorphism!algebraic} property:
\begin{defn}
\label{de:calculus-old}
 An \emph{functional calculus} for an element
 \(a\in\algebra{A}\)\ is a continuous 
linear mapping
\(\Phi: \mathcal{ A}\rightarrow \algebra{A}\)\ such that
\begin{enumerate} 
\item 
 \(\Phi\)\ is a unital \emph{algebra homomorphism}
 \begin{displaymath}
   \Phi(f \cdot g)=\Phi(f) \cdot \Phi (g).
\end{displaymath}
\item 
 There is an initialisation condition: \(\Phi[v_0]=a\)\ for
 for a fixed function \(v_0\), e.g. \(v_0(z)=z\). 
\end{enumerate}
\end{defn}

The most typical definition of the spectrum is seemingly independent
and uses the important notion of resolvent:
\begin{defn}
  \label{de:spectrum}
  A \emph{resolvent}\index{resolvent} of element \(a\in\algebra{A}\)\
  is the function \(R(\lambda)=(a-\lambda e)^{-1}\), which is the
  image under \(\Phi\)\ of the Cauchy kernel \((z-\lambda)^{-1}\).

  A \emph{spectrum}\index{spectrum} of \(a\in\algebra{A}\)\ is the set
  \(\spec a\)\ of singular points of its resolvent \(R(\lambda)\).
\end{defn}
Then the following important theorem links spectrum and functional calculus
together. 
\begin{thm}[Spectral Mapping]
  \index{theorem!spectral mapping}
  \index{spectrum!mapping}
  \label{th:spectral-mapping}
  For a function \(f\) suitable for the   functional calculus:
   \begin{equation}
     \label{eq:spectral-mapping}
     f(\spec a)=\spec  f(a).
   \end{equation}
\end{thm}

However the power of the classic spectral theory rapidly decreases if
we move beyond the study of one normal operator (e.g. for
quasinilpotent ones) and is virtually nil if we consider several
non-commuting ones.  Sometimes these severe limitations are seen to be
irresistible and alternative constructions, i.e. model theory%
\index{functional!model}%
\index{model!functional} cf.
Example~\ref{ex:functional-model} and~\cite{Nikolskii86}, were
developed.

Yet the spectral theory can be revived from a fresh start. While three
compon\-ents---functional calculus, spectrum, and spectral mapping
theorem---are highly interdependent in various ways 
we will nevertheless arrange them as follows: 

\begin{enumerate}
\item Functional  calculus is an \emph{original} notion defined in
  some independent terms;
\item Spectrum (or more specifically \emph{contravariant spectrum})%
  \index{spectrum!contravariant}%
  \index{contravariant!spectrum} (or spectral decomposition) is derived
  from previously defined functional calculus as its \emph{support}
  (in some appropriate sense);
\item Spectral mapping theorem then should drop out naturally in the
  form~\eqref{eq:spectral-mapping} or some its variation.
\end{enumerate}

Thus the entire scheme depends from the notion of the functional
calculus and our ability to escape limitations of
Definition~\ref{de:calculus-old}.  The first known to the present
author definition of functional calculus not linked to algebra
homomorphism property was the Weyl functional calculus defined by an
integral formula~\cite{Anderson69}. Then its intertwining property
with affine transformations of Euclidean space was proved as a
theorem. However it seems to be the only ``non-homomorphism'' calculus
for decades.

The different approach to whole range of calculi was given
in~\cite{Kisil95i} and developed in~\cites{Kisil98a,Kisil02a,Kisil04d,Kisil10c}
in terms of \emph{intertwining operators}%
\index{intertwining operator}%
\index{operator!intertwining} for group representations.
It was initially targeted for several non-commuting operators because
no non-trivial algebra homomorphism is possible with a commutative
algebra of function in this case.  However it emerged later that the
new definition is a useful replacement for classical one across all
range of problems.

In the following Subsections we will support the last claim by
consideration of the simple known problem: characterisation a \(n
\times n\)\ matrix up to similarity. Even that ``freshman'' question
could be only sorted out by the classical spectral theory for a small
set of diagonalisable matrices. Our solution in terms of new spectrum
will be full and thus unavoidably coincides with one given by the
Jordan normal form of matrices. Other more difficult questions are the
subject of ongoing research.

\subsection{Intertwining Group Actions on Functions and Operators}
\label{sec:anoth-appr-analyt}

Any functional calculus uses properties of \emph{functions} to model
properties of \emph{operators}. Thus changing our viewpoint on
functions, as was done in Section~\ref{sec:analytic-functions}, we
could get another approach to operators. The two main possibilities
are encoded in Definitions~\ref{de:covariant-calculus}
and~\ref{de:conravariant-calculus}: we can assign a certain function
to the given operator or wise verse. Here we consider the second
possibility and treat the first in the
Subsection~\ref{sec:funct-model-spectr}.

The representation \(\uir{}{1}\)~\eqref{eq:rho-1-1} is unitary
irreducible when acts on the Hardy space \(\FSpace{H}{2}\).%
\index{space!Hardy}%
\index{Hardy!space}
Consequently we have one more reason to abolish the template
definition~\ref{de:calculus-old}: \(\FSpace{H}{2}\) is \emph{not} an
algebra. Instead we replace the \emph{homomorphism property} by a
\emph{symmetric covariance}:
\begin{defn}[\cite{Kisil95i}]
  \label{de:functional-calculus-new}
  An \emph{contravariant analytic calculus}%
  \index{contravariant!calculus!analytic}%
  \index{calculus!contravariant!analytic}%
  \index{analytic!contravariant calculus} for an element
  \(a\in\algebra{A}\)\ and  an
  \(\algebra{A}\)-module \(M\)\ is a \emph{continuous 
  linear} mapping
  \(\Phi:\FSpace{A}{}(\Space{D}{})\rightarrow \FSpace{A}{}(\Space{D}{},M)\)\ such that 
  \begin{enumerate} 
  \item \(\Phi\)\ is an \emph{intertwining operator}%
    \index{intertwining operator}%
    \index{operator!intertwining}
    \begin{displaymath}
      \Phi\rho_1=\rho_a\Phi
    \end{displaymath}
    between two representations of the
    \(\SL\)\ group \(\rho_1\)~\eqref{eq:rho-1-1} and \(\rho_a\)\
    defined below in~\eqref{eq:rho-a}.
  \item There is an initialisation condition: \(\Phi[v_0]=m\)\ for
    \(v_0(z)\equiv 1\) and \(m\in M\), where \(M\) is a left
    \(\algebra{A}\)-module.  
  \end{enumerate}
\end{defn} 
Note that our functional calculus released from the homomorphism
condition can take value in any left \(\algebra{A}\)-module \(M\),
which however could be \(\algebra{A}\) itself if suitable. This add
much flexibility to our construction.

The earliest functional calculus, which is \emph{not} an algebraic
homomorphism, was the Weyl functional calculus and
was defined just by an integral formula as an operator valued
distribution~\cite{Anderson69}. In that paper
(joint) spectrum%
\index{joint!spectrum}%
\index{spectrum!joint} was defined as support of the Weyl calculus, i.e. as
the set of point where this operator valued distribution does not
vanish. We also define
the spectrum as a support of functional calculus, but due to our
Definition~\ref{de:functional-calculus-new} it will means the set of
non-vanishing intertwining operators with primary subrepresentations.
\begin{defn}
  \label{de:spectrum-new}
  A corresponding \emph{spectrum}%
  \index{spectrum!contravariant}%
  \index{contravariant!spectrum} of \(a\in\algebra{A}\) is the
  \emph{support} of the functional
  calculus \(\Phi\), i.e. the collection of intertwining operators of
  \(\rho_a\) with \emph{primary representations}%
  \index{primary!representation}%
  \index{representation!primary}~\amscite{Kirillov76}*{\S~8.3}.
\end{defn}
More variations of contravariant functional calculi are obtained from other groups and their
representations~\cites{Kisil95i,Kisil98a,Kisil02a,Kisil04d,Kisil10c}. 

A simple but important observation is that the M\"obius
transformations~\eqref{eq:moebius} can be easily extended to any 
Banach algebra.
Let \(\algebra{A}\) be a Banach algebra with the unit \(e\), 
an element \(a\in\algebra{A}\) with \(\norm{a}<1\) be fixed, then%
\index{representation!$\SL$ group!in Banach space}%
\index{$\SL$ group!representation!in Banach space}
\begin{equation}
  \label{eq:sl2-on-A}
  g: a \mapsto g\cdot a=(\bar{\alpha} a -\bar{\beta} e)(\alpha e-\beta a)^{-1}, \qquad
  g\in\SL
\end{equation}
is a well defined \(\SL\) action on a subset \(\Space{A}{}=\{g\cdot
a \such g\in 
\SL\}\subset\algebra{A}\), i.e. \(\Space{A}{}\) is a \(\SL\)-homogeneous
space. Let us define the \emph{resolvent}\index{resolvent} function
\(R(g,a):\Space{A}{}\rightarrow \algebra{A}\):
\begin{displaymath}
  R(g, a)=(\alpha e-\beta a)^{-1} \quad 
\end{displaymath}
then 
\begin{equation}
  \label{eq:ind-rep-multipl}
  R(g_1,\mathsf{a})R(g_2,g_1^{-1}\mathsf{a})=R(g_1g_2,\mathsf{a}).
\end{equation}
The last identity is well known in representation
theory~\amscite{Kirillov76}*{\S~13.2(10)} and is a key ingredient of
\emph{induced representations}%
\index{induced!representation}%
\index{representation!induced}. Thus we can again
linearise~\eqref{eq:sl2-on-A}, cf.~\eqref{eq:rho-1-1}, in
the space of continuous functions \(\FSpace{C}{}(\Space{A}{},M)\)
with values in  a left
\(\algebra{A}\)-module \(M\), e.g. \(M=\algebra{A}\):
\begin{eqnarray}
  \rho_a(g_1): f(g^{-1}\cdot a ) &\mapsto&
  R(g_1^{-1}g^{-1}, a)f(g_1^{-1}g^{-1}\cdot a) \label{eq:rho-a}\\
  &&\quad =
  (\alpha' e-\beta'a)^{-1} \,
  f\left(
    \frac{\bar{\alpha}' \cdot a - \bar{\beta}' e}{\alpha'  e -\beta' a} 
  \right).  \nonumber
\end{eqnarray}  
For any \(m\in M\) we can define a \(K\)-invariant
\emph{vacuum vector}%
\index{mother wavelet}%
\index{wavelet!mother} as \(v_m(g^{-1}\cdot 
a)=m\otimes v_0(g^{-1}\cdot a) \in \FSpace{C}{}(\Space{A}{},M)\). 
It generates the associated with \(v_m\) family of \emph{coherent
  states}\index{wavelet} \(v_m(u,a)=(ue-a)^{-1}m\), where \(u\in\Space{D}{}\).

The \emph{wavelet transform}%
\index{wavelet!transform}  defined by
the same common formula based on coherent states
(cf.~\eqref{eq:cauchy}):
\begin{equation}
  \label{eq:cauchy-int-oper}
  \oper{W}_m f(g)= \scalar{f}{\rho_a(g) v_m},
\end{equation}
is a version of Cauchy integral%
\index{integral!Cauchy}%
\index{Cauchy!integral}, which maps
\(\FSpace{L}{2}(\Space{A}{})\) to \(\FSpace{C}{}(\SL,M)\). It is
 closely related (but not identical!) to the
Riesz-Dunford functional calculus:  the traditional functional
calculus is given by the case:
\begin{displaymath}
  \Phi: f \mapsto \oper{W}_m f(0) \qquad\textrm{ for } M=\algebra{A}
  \textrm{ and } m=e.
\end{displaymath}

The both conditions---the intertwining property and initial
value---required by Definition~\ref{de:functional-calculus-new} easily
follows from our construction.  Finally, we wish to provide an example
of application of the Corollary~\ref{co:cauchy-riemann-integ}.
\begin{example}
  \label{ex:min-poly-integral}
  Let \(a\) be an operator and \(\phi\) be a function which annihilates
  it, i.e. \(\phi(a)=0\). For example, if \(a\) is a matrix \(\phi\) can be its minimal
  polynomial. From the integral representation of the contravariant 
  calculus on \(G=\SL\) we can rewrite the annihilation property like this:
  \begin{displaymath}
    \int_G \phi(g) R(g,a)\,dg=0.
  \end{displaymath}
  Then the vector-valued function \([\oper{W}_m f](g)\) defined
  by~\eqref{eq:cauchy-int-oper} shall satisfy to the following condition:
  \begin{displaymath}
        \int_G \phi(g')\, [\oper{W}_m f] (gg')\,dg'=0
  \end{displaymath}
  due to the Corollary~\ref{co:cauchy-riemann-integ}.
\end{example}

\subsection[Jet Bundles and Prolongations]{Jet Bundles and Prolongations of $\rho_1$}
\label{sec:jet-bundl-prol-1}
Spectrum was defined in~\ref{de:spectrum-new} as
the \emph{support}%
\index{functional!calculus!support|see{spectrum}}%
\index{calculus!functional!support|see{spectrum}}%
\index{support!functional calculus|see{spectrum}}%
\index{spectrum}
 of our functional calculus. To elaborate its meaning we
need the notion of a \emph{prolongation}\index{prolongation} of representations introduced by
S.~Lie, see  \cite{Olver93,Olver95} for a detailed exposition.

\begin{defn} \textup{\amscite{Olver95}*{Chap.~4}}
  Two holomorphic functions have \(n\)th \emph{order contact} in a point
  if their value and their first \(n\) derivatives agree at that point,
  in other words their Taylor expansions are the same in first \(n+1\)
  terms. 

  A point \((z,u^{(n)})=(z,u,u_1,\ldots,u_n)\) of the \emph{jet space}\index{jet}
  \(\Space{J}{n}\sim\Space{D}{}\times\Space{C}{n}\)
  \index{$\Space{J}{n}$ (jet space)} is the equivalence
  class of holomorphic functions having \(n\)th contact at the point \(z\)
  with the polynomial:
  \begin{equation}\label{eq:Taylor-polynom}
    p_n(w)=u_n\frac{(w-z)^n}{n!}+\cdots+u_1\frac{(w-z)}{1!}+u.
  \end{equation}
\end{defn}

For a fixed \(n\) each holomorphic function
\(f:\Space{D}{}\rightarrow\Space{C}{}\) has \(n\)th \emph{prolongation}
(or \emph{\(n\)-jet}) \(\obj[_n]{j}f: \Space{D}{} \rightarrow
\Space{C}{n+1}\): 
\begin{equation}\label{eq:n-jet}
  \obj[_n]{j}f(z)=(f(z),f'(z),\ldots,f^{(n)}(z)).
\end{equation}The graph \(\mathsf{a}mma^{(n)}_f\) of \(\obj[_n]{j}f\) is a
submanifold of \(\Space{J}{n}\) which is section of the \emph{jet
bundle}%
\index{jet!bundle} over \(\Space{D}{}\) with a fibre \(\Space{C}{n+1}\). We also
introduce a notation \(J_n\) for the map \(
  J_n:f\mapsto\mathsf{a}mma^{(n)}_f
\) of a holomorphic \(f\) to the graph \(\mathsf{a}mma^{(n)}_f\) of its \(n\)-jet
\(\obj[_n]{j}f(z)\)~\eqref{eq:n-jet}.

One can prolong any map of functions \(\psi: f(z)\mapsto [\psi f](z)\) to
a map \(\psi^{(n)}\) of \(n\)-jets by the formula
\begin{equation}\label{eq:prolong-def}
  \psi^{(n)} (J_n f) = J_n(\psi f).
\end{equation} For example such a prolongation \(\rho_1^{(n)}\) of the
representation \(\rho_1\) of the group \(\SL\) in
\(\FSpace{H}{2}(\Space{D}{})\) (as any other representation of a Lie
group~\cite{Olver95}) will be again a representation of
\(\SL\). Equivalently we can say that \(J_n\) \emph{intertwines} \(\rho_1\) and
\(\rho^{(n)}_1\):
\begin{displaymath}
   J_n \rho_1(g)= \rho_1^{(n)}(g) J_n \quad
  \textrm{ for all } g\in\SL.
\end{displaymath}
Of course, the representation \(\rho^{(n)}_1\) is not irreducible: any
jet subspace \(\Space{J}{k}\), \(0\leq k \leq n\) is
\(\rho^{(n)}_1\)-invariant subspace of \(\Space{J}{n}\).  However the
representations \(\rho^{(n)}_1\) are \emph{primary}%
\index{primary!representation}%
\index{representation!primary}~\amscite{Kirillov76}*{\S~8.3} in the
sense that they are not sums of two subrepresentations.

The following statement explains why jet spaces appeared in our study
of functional calculus.
\begin{prop}
  \label{pr:Jordan-zero}
  Let matrix \(a\) be a Jordan block of a length \(k\) with the
  eigenvalue\index{eigenvalue} \(\lambda=0\),
  and \(m\) be its root vector of order \(k\), i.e. \(a^{k-1}m\neq
  a^k m =0\). Then the restriction of \(\rho_a\) on the subspace
  generated by \(v_m\) is equivalent to the representation
  \(\rho_1^{k}\).
\end{prop}

\subsection{Spectrum and Spectral Mapping Theorem}
\label{sec:spectr-jord-norm}

Now we are prepared to describe a spectrum of a matrix. Since the
functional calculus is an intertwining operator its support is
a decomposition into intertwining operators with primary
representations%
\index{primary!representation}%
\index{representation!primary} (we could not expect generally 
that these primary subrepresentations are irreducible).

Recall the transitive on \(\Space{D}{}\) group of inner
automorphisms of \(\SL\), which can send any
\(\lambda\in\Space{D}{}\) to \(0\) and are actually parametrised by
such a \(\lambda\).
This group extends Proposition~\ref{pr:Jordan-zero} to the complete
characterisation of \(\rho_a\) for matrices.
\begin{prop}  
  Representation \(\rho_a\) is equivalent to a direct sum of the
  prolongations \(\rho_1^{(k)}\) of \(\rho_1\) in the \(k\)th jet space
  \(\Space{J}{k}\) intertwined with inner automorphisms. Consequently
  the \emph{spectrum}%
  \index{spectrum!contravariant}%
  \index{contravariant!spectrum} of \(a\) (defined via the functional calculus
  \(\Phi=\oper{W}_m\)) labelled exactly by \(n\) pairs of numbers
  \((\lambda_i,k_i)\), \(\lambda_i\in\Space{D}{}\),
  \(k_i\in\Space[+]{Z}{}\), \(1\leq i \leq n\) some of whom could
  coincide.
\end{prop}
Obviously this spectral theory is a fancy restatement of the \emph{Jordan
  normal form}%
\index{Jordan!normal form}%
\index{matrix!Jordan normal form} of matrices.

\begin{figure}[tb]
  \begin{center}
    (a) \includegraphics[scale=.85]{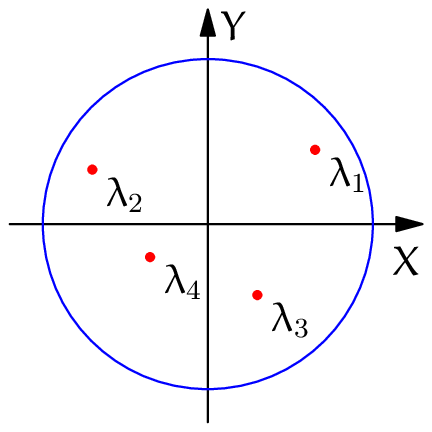}\hfill
    (b)\includegraphics[scale=.85]{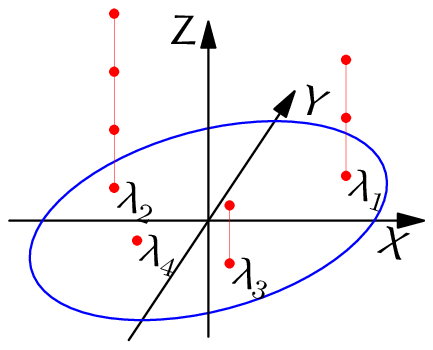}\hfill
    (c)\includegraphics[scale=.85]{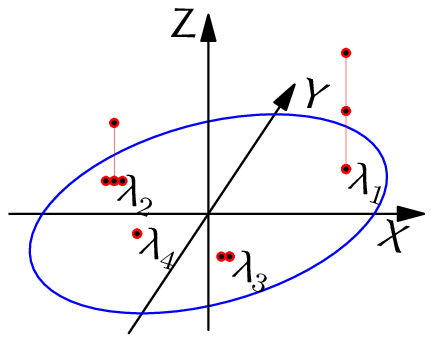}
    \caption[Three dimensional spectrum]{Classical spectrum of the matrix
      from the Ex.~\ref{ex:3dspectrum} is shown at (a). Contravariant
      spectrum of the same matrix in the jet space is drawn at (b).  The
      image of the contravariant spectrum under the map from
      Ex.~\ref{ex:spectral-mapping} is presented at (c).} 
    \label{fig:3dspectrum}
  \end{center}
\end{figure}

\begin{example}
  \label{ex:3dspectrum}
  Let \(J_k(\lambda)\) denote the Jordan block of the length \(k\) for the
  eigenvalue \(\lambda\). In Fig.~\ref{fig:3dspectrum} there are two
  pictures of the spectrum for the matrix
  \begin{displaymath}
    a=J_3\left(\lambda_1\right)\oplus     J_4\left(\lambda_2\right) 
    \oplus J_1\left(\lambda_3\right) \oplus      J_2\left(\lambda_4\right),
  \end{displaymath} 
  where
  \begin{displaymath}
    \lambda_1=\frac{3}{4}e^{i\pi/4}, \quad
    \lambda_2=\frac{2}{3}e^{i5\pi/6}, \quad
    \lambda_3=\frac{2}{5}e^{-i3\pi/4}, \quad
    \lambda_4=\frac{3}{5}e^{-i\pi/3}.
  \end{displaymath} Part (a) represents the conventional two-dimensional
  image of the spectrum, i.e. eigenvalues of \(a\), and
  \href{http://www.maths.leeds.ac.uk/~kisilv/calc1vr.gif}{(b) describes
  spectrum \(\spec{} a\) arising from the wavelet construction}. The
  first image did not allow to distinguish \(a\) from many other
  essentially different matrices, e.g. the diagonal matrix
  \begin{displaymath}
    \diag\left(\lambda_1,\lambda_2,\lambda_3,\lambda_4\right),
  \end{displaymath}
  which even have a different dimensionality.
  At the same time Fig.~\ref{fig:3dspectrum}(b)
  completely characterise \(a\) up to a similarity. Note that each point of
  \(\spec a\) in Fig.~\ref{fig:3dspectrum}(b) corresponds to a particular
  root vector, which spans a primary subrepresentation.
\end{example}


As was mentioned in the beginning of this section a resonable spectrum
should be linked to the corresponding functional calculus by an
appropriate spectral mapping theorem. The new version of spectrum is
based on prolongation of \(\rho_1\) into jet\index{jet} spaces (see
Section~\ref{sec:jet-bundl-prol-1}). Naturally a correct version of
spectral mapping theorem should also operate in jet spaces.

Let \(\phi: \Space{D}{} \rightarrow \Space{D}{}\) be a holomorphic
map, let us define its action on functions \([\phi_*
f](z)=f(\phi(z))\). According to the general formula~\eqref{eq:prolong-def}
we can define the prolongation
\(\phi_*^{(n)}\) onto the jet space \(\Space{J}{n}\). Its associated
action \(\rho_1^k \phi_*^{(n)}=\phi_*^{(n)}\rho_1^n\) on the pairs
\((\lambda,k)\) is given by the formula:
\begin{equation}
  \label{eq:phi-star-action}
  \phi_*^{(n)}(\lambda,k)=\left(\phi(\lambda),
    \left[\frac{k}{\deg_\lambda \phi}\right]\right),
\end{equation}
where \(\deg_\lambda \phi\) denotes the degree of zero of the function
\(\phi(z)-\phi(\lambda)\) at the point \(z=\lambda\) and \([x]\) denotes
the integer part of \(x\). 

\begin{thm}[Spectral mapping] 
  Let \(\phi\) be a holomorphic mapping  \(\phi: \Space{D}{}
  \rightarrow \Space{D}{}\) and its prolonged action \(\phi_*^{(n)}\) defined
  by~\eqref{eq:phi-star-action}, then
  \begin{displaymath}
    \spec \phi(a) = \phi_*^{(n)} \spec a. 
  \end{displaymath}
\end{thm} 

The explicit expression of~\eqref{eq:phi-star-action} for
\(\phi_*^{(n)}\), which involves derivatives of \(\phi\) upto \(n\)th order,
is known, see for example~\amscite{HornJohnson94}*{Thm.~6.2.25}, but was
not recognised before as form of spectral mapping.

\begin{example}
  \label{ex:spectral-mapping}
  Let us continue with Example~\ref{ex:3dspectrum}. Let \(\phi\) map
  all four eigenvalues \(\lambda_1\), \ldots, \(\lambda_4\) of the
  matrix \(a\) into themselves. Then Fig.~\ref{fig:3dspectrum}(a) will
  represent the classical spectrum of \(\phi(a)\) as well as \(a\).

  However Fig.~\ref{fig:3dspectrum}(c) shows mapping of the new
  spectrum for the case
  \(\phi\)  has
  \emph{orders of zeros}%
  \index{orders!of zero}%
  \index{zero!order, of} at these points as follows: the order \(1\)
  at \(\lambda_1\), exactly the order \(3\) at \(\lambda_2\), an order
  at least \(2\) at \(\lambda_3\), and finally any order at
  \(\lambda_4\).
\end{example}

\subsection{Functional Model and Spectral Distance}
\label{sec:funct-model-spectr}

Let \(a\) be a matrix and \(\mu(z)\) be its \emph{minimal
  polynomial}%
\index{minimal!polynomial}%
\index{polynomial!minimal}:
\begin{displaymath}
  \mu_a(z)=(z-\lambda_1)^{m_1}\cdot \ldots\cdot (z-\lambda_n)^{m_n}.
\end{displaymath}
If all eigenvalues\index{eigenvalue} \(\lambda_i\) of \(a\) (i.e. all roots of
\(\mu(z)\) belong to the unit disk we can consider the respective
\emph{Blaschke product}%
\index{Blaschke!product}%
\index{product!Blaschke}
\begin{displaymath}
  B_a(z)=\prod_{i=1}^n\left(\frac{z-\lambda _i}{1-\overline{\lambda_i}z}\right)^{m_i},
\end{displaymath}
such that its numerator coincides with the minimal polynomial
\(\mu(z)\).  Moreover, for an unimodular \(z\) we have
\(B_a(z)=\mu_a(z)\overline{\mu^{-1}_a(z)}z^{-m}\), where
\(m=m_1+\ldots +m_n\). We also have the following covariance property:
\begin{prop}
  \label{pr:covariance-funct-model}
  The above correspondence \(a\mapsto B_a\) intertwines the \(\SL\)
  action~\eqref{eq:sl2-on-A} on the matrices with the
  action~\eqref{eq:rho-1-1} with \(k=0\) on functions. 
\end{prop}
The result follows from the observation that every elementary product
\(\frac{z-\lambda _i}{1-\overline{\lambda_i}z}\) is the Moebius
transformation%
\index{M\"obius map}%
\index{map!M\"obius} of \(z\) with the matrix \(\begin{pmatrix}
  1&-\lambda _i\\-\overline{\lambda_i}&1
\end{pmatrix}\). Thus the correspondence \(a\mapsto  B_a(z)\) is a
covariant (symbolic) calculus%
\index{covariant!calculus}%
\index{calculus!covariant}
in the sense of the Defn.~\ref{de:covariant-calculus}.  See also the
Example~\ref{ex:functional-model}. 

The Jordan normal form of a matrix provides a description, which is
equivalent to its contravariant spectrum%
\index{contravariant!spectrum}%
\index{spectrum!contravariant}.  From various viewpoints, e.g.
numerical approximations, it is worth to consider its stability under
a perturbation. It is easy to see, that an arbitrarily small
disturbance breaks the Jordan structure of a matrix. However, the
result of random small perturbation will not be random, its nature is
described by the following remarkable theorem:
\begin{thm}[Lidskii~\cite{Lidskii66a}, see
  also~\cite{MoroBurkeOverton97a}]
  \index{Lidskii theorem}
  \index{theorem!Lidskii}
  \label{th:lidskii}
  Let \(J_n\) be a Jordan block of a length \(n>1\) with zero
  eigenvalues and \(K\) be an arbitrary matrix. 
  Then eigenvalues of the perturbed matrix \(J_n+\varepsilon^n  K\) admit
  the expansion
  \begin{displaymath}
    \lambda_j=  \varepsilon \xi^{1/n}  +o(\varepsilon), \qquad j=1,\ldots,n,  
  \end{displaymath}
  where \(\xi^{1/n}\)  represents all \(n\)-th complex roots of certain \(\xi\in\Space{C}{}\).
\end{thm}
\begin{figure}
  \centering
  (a)\includegraphics[scale=.9]{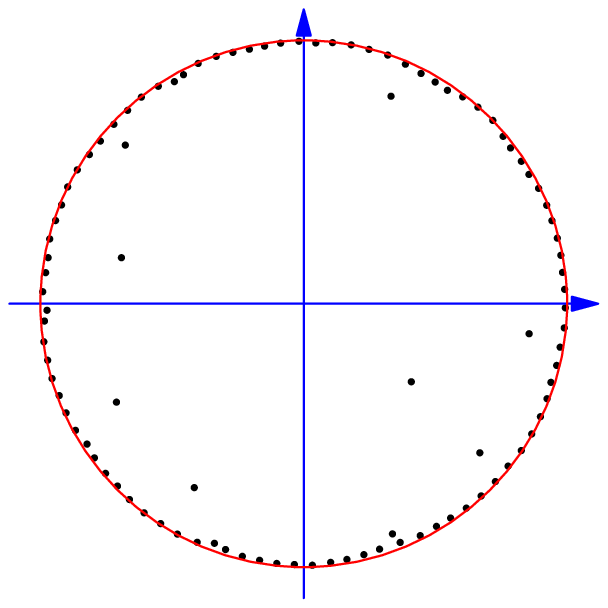}\hfill
  (b)\includegraphics[scale=.9]{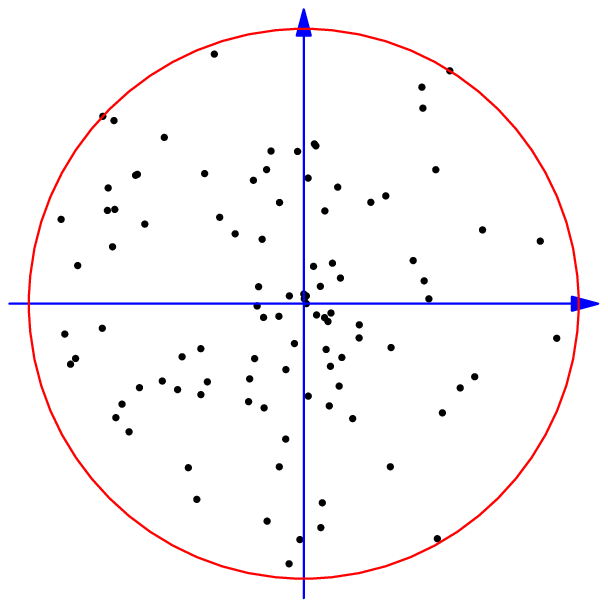}\\
  \caption[Spectral stability]{Perturbation of the Jordan block's
    spectrum: 
    (a) The spectrum of the perturbation \(J_{100}+\varepsilon^100 K\) of the Jordan
    block \(J_{100}\) by a random matrix \(K\). 
    (b) The spectrum of the random matrix \(K\).}
  \label{fig:jord-perturb}
\end{figure}
The left picture in Fig.~\ref{fig:jord-perturb} presents a
perturbation of a Jordan block \(J_{100}\) by a random matrix.
Perturbed eigenvalues are close to vertices of a right polygon with
\(100\) vertices. Those regular arrangements occur despite of the fact
that eigenvalues of the matrix \(K\) are dispersed through the unit
disk (the right picture in Fig.~\ref{fig:jord-perturb}). In a sense it
is rather the Jordan block regularises eigenvalues of \(K\) than \(K\)
perturbs the eigenvalue of the Jordan block.

Although the Jordan structure itself is extremely fragile, it still
can be easily guessed from a perturbed eigenvalues. Thus there exists
a certain characterisation of matrices which is stable under small
perturbations. We will describe a sense, in which the covariant
spectrum%
\index{spectrum!contravariant!stability}%
\index{contravariant!spectrum!stability}%
\index{stability!contravariant spectrum} of the matrix
\(J_n+\varepsilon^n K\) is stable for small \(\varepsilon\).  For this
we introduce the covariant version of spectral distances motivated by
the functional model%
\index{functional!model}%
\index{model!functional}. Our definition is different from other types
known in the literature~\amscite{Tyrtyshnikov97a}*{Ch.~5}.
\begin{defn}
  Let \(a\) and \(b\) be two matrices with all their eigenvalues
  sitting inside of the unit disk and \(B_a(z)\) and \(B_b(z)\) be
  respective Blaschke products as defined above. The \emph{(covariant)
    spectral distance}%
  \index{spectral!covariant distance}%
  \index{distance!covariant spectral}%
  \index{covariant!spectral distance} \(d(a,b)\) between \(a\) and
  \(b\) is equal to the distance \(\norm[2]{B_a-B_b}\) between
  \(B_a(z)\) and \(B_b(z)\) in the Hardy space%
  \index{Hardy!space}%
  \index{space!Hardy} on the unit circle.
\end{defn}
Since the spectral distance is defined through the distance in 
\(\FSpace{H}{2}\) all standard axioms of a distance are automatically
satisfied.  For a Blaschke products we have
\(\modulus{B_a(z)}=1\) if \(\modulus{z}=1\), thus \(\norm[p]{B_a}=1\)
in any \(\FSpace{L}{p}\) on the unit circle. Therefore an alternative
expression for the spectral distance is:
\begin{displaymath}
  d(a,b)=2(1-\scalar{B_a}{B_b}).
\end{displaymath}
In particular, we always have \(0\leq d(a,b) \leq 2\). 
We get an obvious consequence of Prop.~\ref{pr:covariance-funct-model},
which justifies the name of the covariant spectral distance:
\begin{cor}
  For any \(g\in\SL\) we have  \(d(a,b)=d(g\cdot a, g\cdot a)\), where
  \(\cdot\) denotes the M\"obius action~\eqref{eq:sl2-on-A}.
\end{cor}
An important property of the covariant spectral distance is its
stability under small perturbations.
\begin{thm}
  For \(n=2\) let \(\lambda_1(\varepsilon )\) and
  \(\lambda_2(\varepsilon )\) be eigenvalues of the matrix
  \(J_2+\varepsilon ^2\cdot K\) for some matrix \(K\). Then
  \begin{equation}
    \label{eq:eigenvalue-est}
    \modulus{\lambda_1(\varepsilon )}+\modulus{\lambda_2(\varepsilon
      )}=O(\varepsilon ),\quad
  \text{ however }\quad
    \modulus{\lambda_1(\varepsilon )+\lambda_2(\varepsilon
      )}=O(\varepsilon ^2).
  \end{equation}
  The spectral distance from the \(1\)-jet\index{jet} at \(0\) to two
  \(0\)-jets at points \(\lambda_1\) and \(\lambda_2\) bounded only by
  the first condition in~\eqref{eq:eigenvalue-est} is \(O(\varepsilon
  ^2)\).  However the spectral distance between \(J_2\) and
  \(J_2+\varepsilon ^2\cdot K\) is \(O(\varepsilon ^4)\).
\end{thm}
In other words, a matrix with eigenvalues satisfying to the Lisdkii
condition from the Thm.~\ref{th:lidskii} is much closer to the Jordan
block \(J_2\) than a generic one with eigenvalues of the same order.
Thus the covariant spectral distance is more stable under perturbation
that magnitude of eigenvalues. For \(n=2\) a proof can be
forced through a direct calculation. We also conjecture that the
similar statement is true for any \(n\geq 2\).

\subsection{Covariant Pencils of Operators}
\label{sec:quad-penc-oper}

Let \(H\) be a real Hilbert space, possibly of finite dimensionality.
For bounded linear operators \(A\) and \(B\) consider the
\emph{generalised eigenvalue problem}%
\index{generalised!eigenvalue}%
\index{eigenvalue!generalised}, that is finding a scalar \(\lambda \)
and a vector \(x\in H\) such that:
\begin{equation}
  \label{eq:eigenvalue-generalised}
  Ax=\lambda Bx \qquad \text{or equivalently} \qquad
  (A-\lambda B)x=0.
\end{equation}
The standard eigenvalue problem corresponds to the case \(B=I\),
moreover for an invertible \(B\) the generalised problem can be
reduced to the standard one for the operator \(B^{-1}A\). Thus it is
sensible to introduce the equivalence relation on the pairs of
operators:
\begin{equation}
  \label{eq:pairs-op-equivalence}
  (A,B)\sim(DA,DB) \quad \text{for any invertible operator } D. 
\end{equation}

We may treat the pair \((A,B)\) as a column vector \(
\begin{pmatrix}
  A\\B
\end{pmatrix}\).  Then there is an action of the \(\SL\) group on the
pairs:
\begin{equation}
  \label{eq:SL2-act-pairs-op}
  g\cdot
  \begin{pmatrix}
    A\\B
  \end{pmatrix}=
  \begin{pmatrix}
    aA+bB\\cA+dB
  \end{pmatrix},\qquad \text{where }
  g=
  \begin{pmatrix}
    a&b\\c&d
  \end{pmatrix}
  \in\SL.
\end{equation}
If we consider this \(\SL\)-action subject to the equivalence
relation~\eqref{eq:pairs-op-equivalence} then we will arrive to a
version of the linear-fractional transformation of the operator
defined in~\eqref{eq:sl2-on-A}. There is a connection of the
\(\SL\)-action~\eqref{eq:SL2-act-pairs-op} to the
problem~\eqref{eq:eigenvalue-generalised} through the following
intertwining relation%
\index{intertwining operator}%
\index{operator!intertwining}:
\begin{prop}
  Let \(\lambda\) and \(x\in H\) solve the generalised eigenvalue
  problem~\eqref{eq:eigenvalue-generalised} for the pair
  \((A,B)\). Then the pair \((C,D)=g\cdot (A,B)\), \(g\in\SL\) has a
  solution \(\mu\) and \(x\), where
  \begin{displaymath}
    \mu=g\cdot \lambda=\frac{a\lambda +b}{c\lambda +d},\qquad\text{for
    } g=
    \begin{pmatrix}
      a&b\\c&d
    \end{pmatrix}\in\SL,
  \end{displaymath}
  is defined by the M\"obius transformation~\eqref{eq:moebius}%
  \index{M\"obius map}%
  \index{map!M\"obius}. 
\end{prop}

In other words the correspondence 
\begin{displaymath}
  (A,B)\mapsto \text{all generalised eigenvalues}
\end{displaymath}
is another realisation of a covariant calculus%
\index{covariant!calculus}%
\index{calculus!covariant} in the sense of
Defn.~\ref{de:covariant-calculus}.  The collection of all pairs
\(g\cdot (A,B)\), \(g\in\SL\) is an example of \emph{covariant
  pencil}%
\index{covariant!pencil}%
\index{pencil!covariant} of operators. This set is a
\(\SL\)-homogeneous spaces, thus it shall be within the classification
of such homogeneous spaces provided in the
Subsection~\ref{sec:hypercomplex-numbers}.
\begin{example} 
  \label{ex:eph-op-pairs}
  It is easy to demonstrate that all existing homogeneous spaces can
  be realised by matrix pairs.
  \begin{enumerate}
  \item Take the pair \((O, I)\) where \(O\) and \(I\) are the zero
    and identity \(n\times n\) matrices respectively. Then any
    transformation of this pair by a lower-triangular matrix from
    \(\SL\) is equivalent to  \((O, I)\). The respective homogeneous
    space is isomorphic to the real line with the M\"obius
    transformations~\eqref{eq:moebius}. 
  \item Consider \(H=\Space{R}{2}\). Using the notations \(\alli\)
    from Subsection~\ref{sec:make-guess-three} we define three
    realisations (elliptic, parabolic and hyperbolic) of an operator
    \(A_{\alli}\):
    \begin{equation}
      \label{eq:A-alli-matr}
      A_{\rmi}=\begin{pmatrix}
        0&1\\-1&0
      \end{pmatrix}, \qquad
      A_{\rmp}=\begin{pmatrix}
        0&1\\0&0
      \end{pmatrix},\qquad
      A_{\rmh}=\begin{pmatrix}
        0&1\\1&0
      \end{pmatrix}.
    \end{equation}
    Then for an arbitrary element \(h\) of the subgroup \(K\), \(N\)
    or \(A\) the respective (in the sense of the
    Principle~\ref{pr:simil-corr-principle}) pair \(h\cdot
    (A_{\alli},I)\) is equivalent to \((A_{\alli},I)\) itself.  Thus
    those three homogeneous spaces are isomorphic to the elliptic,
    parabolic and hyperbolic half-planes%
    \index{half-plane} under respective actions of
    \(\SL\).  Note, that \(A_{\alli}^2 =\alli^2 I\), that is \(A_\alli\) is
    a model for hypercomplex units.
    \item \label{it:oper-SL2-homog}
      Let \(A\) be a direct sum of any two different matrices out of the
      three \(A_{\alli}\) from~\eqref{eq:A-alli-matr}, then the fix
      group of the equivalence class of the pair \((A,I)\) is the
      identity of \(\SL\). Thus the corresponding homogeneous space
      coincides with the group itself.
  \end{enumerate}
\end{example}
Hawing homogeneous spaces generated by pairs of operators we can define
respective functions on those spaces. The special attention is due the
following paraphrase of the resolvent\index{resolvent}:
\begin{displaymath}
  R_{(A,B)}(g)=(cA+d B)^{-1}\qquad\text{where} \quad 
  g^{-1}=\begin{pmatrix}a&b\\c&d
  \end{pmatrix}\in SL_2(\mathbb{R}).
\end{displaymath}
Obviously \(R_{(A,B)}(g)\) contains the essential information about
the pair \((A,B)\). Probably, the function \(R_{(A,B)}(g)\) contains
too much simultaneous information, we may restrict it to get a more
detailed view.  For vectors \(u\), \(v\in H\) we also consider vector
and scalar-valued functions related to the generalised resolvent:
\begin{displaymath}
  R^u_{(A,B)}(g)=(cA+dB)^{-1}u, \qquad\text{ and }\qquad
  R^{(u,v)}_{(A,B)}(g)=\scalar{(cA+dB)^{-1}u}{v}, 
\end{displaymath}
where \((cA+dB)^{-1}u\) is understood as a solution \(w\) of the
equation \(u=(cA+dB)w\) if it exists and is unique, this does not
require the full invertibility of \(cA+dB\).

It is easy to see that the map \((A,B)\mapsto R^{(u,v)}_{(A,B)}(g)\)
is a covariant calculus%
\index{covariant!calculus}%
\index{calculus!covariant} as well. It worth to notice that function
\(R_{(A,B)}\) can again fall into three EPH cases.
\begin{example}
  For the three matrices \(A_{\alli}\) considered in the previous
  Example we denote by \(R_{\alli}(g)\) the resolvetn-type function of
  the pair \((A_{\alli}, I)\). Then:
  \begin{displaymath}
    R_{\rmi}(g)=\frac{1}{c^2+d^2}\begin{pmatrix}
      d&-c\\c&d
    \end{pmatrix},\quad
    R_{\rmp}(g)=\frac{1}{d^2}\begin{pmatrix}
      d&-c\\0&d
    \end{pmatrix},\quad
    R_{\rmh}(g)=\frac{1}{d^2-c^2}\begin{pmatrix}
      d&-c\\-c&d
    \end{pmatrix}.
  \end{displaymath}
  Put \(u=(1,0)\in H\), then \(R_{\alli}(g) u\) is a two-dimensional
  real vector valued functions with components equal to real and
  imaginary part of hypercomplex Cauchy kernel considered
  in~\cite{Kisil11b}.
\end{example}

Consider the space \(L(G)\) of functions spanned by all left
translations of \(R_{(A,B)}(g)\). As usual, a closure in a suitable
metric, say \(\FSpace{L}{p}\), can be taken. The left action \(g:
f(h)\mapsto f(g^{-1}h)\) of \(SL_2(\mathbb{R})\) on this space is a
linear representation of this group. Afterwards the representation can
be decomposed into a sum of primary subrepresentations%
\index{primary!representation}%
\index{representation!primary}.
\begin{example}
  For the matrices \(A_\alli\) the irreducible components are
  isomorphic to analytic spaces of hypercomplex functions under the
  fraction-linear transformations build in
  Subsection~\ref{sec:concl-induc-repr}.
\end{example}
An important observation is that a decomposition into irreducible or
primary components can reveal an EPH structure even in the cases
hiding it on the homogeneous space level.
\begin{example}
  Take the operator \(A=A_{\rmi}\oplus A_{\rmh}\) from the
  Example~\ref{ex:eph-op-pairs}(\ref{it:oper-SL2-homog}). The
  corresponding homogeneous space coincides with the entire
  \(\SL\). However if we take two vectors
  \(u_{\rmi}=(1,0)\oplus(0,0)\) and \(u_{\rmh}=(0,0)\oplus(1,0)\) then
  the respective linear spaces generated by functions \(R_A(g)u_{\rmi}\) and
  \(R_A(g)u_{\rmh}\) will be of elliptic and hyperbolic types
  respectively. 
\end{example}

Let us briefly consider a \emph{quadratic eigenvalue}%
\index{eigenvalue!quadratic}%
\index{quadratic!eigenvalue} problem: for given operators (matrices)
\(A_0\), \(A_1\) and \(A_2\) from \(B(H)\) find a scalar \(\lambda\)
and a vector \(x\in H\) such that
\begin{equation}
  \label{eq:quadratic-ev}
  Q(\lambda)x=0, \qquad \text{where}\quad
  Q(\lambda)=\lambda^2 A_2 + \lambda A_1 + A_0.
\end{equation}
There is a connection with our study of conic sections from
Subsection~\ref{sec:cycles-as-invariant} which we will only hint for
now.  Comparing \eqref{eq:quadratic-ev} with the equation of the
cycle~\eqref{eq:cycle-eq} we can associate the respective
Fillmore--Springer--Cnops--type matrix%
\index{operator!Fillmore--Springer--Cnops construction}%
\index{Fillmore--Springer--Cnops construction!operator}%
\index{construction!Fillmore--Springer--Cnops!operator} to
\(Q(\lambda)\), cf.~\eqref{eq:FSCc-matrix}:
\begin{equation}
  \label{eq:FSC-operator}
  Q(\lambda)=\lambda^2 A_2 + \lambda A_1 + A_0
  \quad
  \longleftrightarrow
  \quad C_Q=
  \begin{pmatrix}
    A_1 & A_0\\
    A_2 & -A_1
  \end{pmatrix}.
\end{equation}
Then we can state the following analogue of
Thm.~\ref{th:FSCc-intertwine} for the quadratic eigenvalues:
\begin{prop}
  Let two quadratic matrix polynomials \(Q\) and \(\tilde{Q}\) are
  such that their FSC matrices~\eqref{eq:FSC-operator} are conjugated
  \(C_{\tilde{Q}}=gC_{{Q}} g^{-1}\) by an element \(g\in\SL\). Then
  \(\lambda\) is a solution of the quadratic eigenvalue problem for
  \(Q\) and \(x\in H\) if and only if \(\mu=g\cdot \lambda\) is a
  solution of the quadratic eigenvalue problem for \(\tilde{Q}\) and \(x\).
  Here \(\mu=g\cdot \lambda\) is the M\"obius
  transformation~\eqref{eq:moebius}%
  \index{M\"obius map}%
  \index{map!M\"obius} associated to \(g\in\SL\).
\end{prop}
So quadratic matrix polynomials are non-commuting%
\index{non-commutative geometry}%
\index{geometry!non-commutative} 
analogues of the cycles\index{cycle} and it would be exciting to
extend the geometry from Section~\ref{sec:geometry} to this
non-com\-mu\-ta\-tive setting as much as possible.
\begin{rem}
  It is beneficial to extend a notion of a scalar in an (generalised)
  eigenvalue problem to an abstract field or ring. For example, we can
  consider pencils of operators/matrices with polynomial coefficients.
  In many circumstances we may factorise the polynomial ring by an
  ideal generated by a collection of algebraic equations. Our work
  with hypercomplex units is the most elementary realisation of this
  setup. Indeed, the algebra of hypercomplex numbers with the
  hypercomplex unit \(\alli\) is a realisation of the polynomial ring
  in a variable \(t\) factored by the single quadratic relation
  \(t^2+\sigma=0\), where \(\sigma=\alli^2\).
\end{rem}

\section{Quantum Mechanics}
\label{sec:quantum-mechanics-1}

Complex valued representations of the Heisenberg group%
\index{Heisenberg!group}%
\index{group!Heisenberg} (also known as Weyl or Heisenberg-Weyl group)
provide a natural framework for quantum
mechanics~\cites{Howe80b,Folland89}. This is the most fundamental
example of the Kirillov orbit method%
\index{orbit!method}%
\index{method!orbits, of}, induced representations%
\index{representation!induced}%
\index{induced!representation} and geometrical quantisation%
\index{geometrical quantisation}%
\index{quantisation!geometrical}
technique~\cites{Kirillov99,Kirillov94a}.  Following the presentation
in Section~\ref{sec:induc-repr} we will consider representations of
the Heisenberg group which are induced by hypercomplex characters of
its centre: complex (which correspond to the elliptic case), dual
(parabolic) and double (hyperbolic).

To describe dynamics of a physical system we use a universal equation
based on inner derivations%
\index{inner!derivation}%
\index{derivation!inner} (commutator%
\index{commutator}) of the convolution
algebra~\citelist{\cite{Kisil00a} \cite{Kisil02e}}.  The complex
valued representations produce the standard framework for quantum
mechanics with the Heisenberg dynamical equation%
\index{Heisenberg!equation}%
\index{equation!Heisenberg}~\cite{Vourdas06a}.

The double number%
\index{number!double}%
\index{double!number} valued representations, with
the hyperbolic unit \(\rmh^2=1\), is a natural source of hyperbolic
quantum mechanics developed for a
while~\cites{Hudson04a,Hudson66a,Khrennikov03a,Khrennikov05a,Khrennikov08a}.
The universal dynamical equation employs hyperbolic commutator in this
case. This can be seen as a \emph{Moyal bracket}%
\index{Moyal bracket!hyperbolic}%
\index{bracket!Moyal!hyperbolic}%
\index{hyperbolic!Moyal bracket} based on the hyperbolic sine
function. The hyperbolic observables act as operators on a Krein
space%
\index{Krein!space}%
\index{space!Krein} with an indefinite inner product. Such spaces are
employed in study of \(\mathcal{PT}\)-symmetric%
\index{PT-symmetry@\(\mathcal{PT}\)-symmetry} Hamiltonians and
hyperbolic unit \(\rmh^2=1\) naturally appear in this
setup~\cite{GuentherKuzhel10a}.

The representations with values in dual numbers%
\index{number!dual}%
\index{dual!number}
provide a convenient description of the classical mechanics. For this
we do not take any sort of semiclassical limit%
  \index{semiclassical!limit}%
  \index{limit!semiclassical}, rather the nilpotency
of the parabolic unit (\(\rmp^2=0\)) do the task. This removes the
vicious necessity to consider the Planck \emph{constant}%
\index{Planck!constant}%
\index{constant!Planck} tending to
zero. 
The dynamical equation takes the Hamiltonian%
\index{Hamilton!equation}%
\index{equation!Hamilton} form. We also describe
classical non-commutative representations of the Heisenberg group
which acts in the first jet space\index{jet}.

\begin{rem}
  It is worth to note that our technique is different from contraction
  technique in the theory of Lie
  groups~\cites{LevyLeblond65a,GromovKuratov05b}. Indeed a contraction
  of the Heisenberg group \(\Space{H}{n}\) is the commutative
  Euclidean group \(\Space{R}{2n}\) which does not recreate neither
  quantum nor classical mechanics.
\end{rem}

The approach provides not only three different types of dynamics, it
also generates the respective rules for addition of probabilities%
\index{probability!quantum}%
\index{quantum!probability} as well.
For example, the quantum interference is the consequence of the same
complex-valued structure, which directs the Heisenberg equation. The
absence of an interference (a particle behaviour) in the classical
mechanics is again the consequence the nilpotency of the parabolic
unit. Double numbers creates the hyperbolic law of additions of
probabilities, which was extensively
investigates~\cites{Khrennikov03a,Khrennikov05a}. There are still
unresolved issues with positivity of the probabilistic interpretation
in the hyperbolic case~\cites{Hudson04a,Hudson66a}. 

\begin{rem}
  \label{re:quantum-complex}
  It is commonly accepted since the Dirac's paper~\cite{Dirac26a} that
  the striking (or even \emph{the only}) difference between quantum
  and classical mechanics is non-commutativity of observables in the
  first case. In particular the Heisenberg commutation
  relations~\eqref{eq:heisenberg-comm} imply the uncertainty
  principle, the Heisenberg equation of motion and other quantum
  features. However, the entire book of Feynman on
  QED~\cite{Feynman1990qed} does not contains any reference to
  non-commutativity. Moreover, our work shows that there is a
  non-commutative formulation of classical mechanics. Non-commutative
  representations of the Heisenberg group in dual numbers implies the
  Poisson dynamical equation and local addition of probabilities in
  Section\ref{sec:class-repr-phase}, which are completely classical.

  This entirely dispels any illusive correlation between
  classical/quantum and commutative/non-commutative. Instead we show
  that quantum mechanics is fully determined by the properties of
  complex numbers. In Feynman's exposition~\cite{Feynman1990qed}
  complex numbers are presented by a clock, rotations of its arm
  encode multiplications by unimodular complex numbers. Moreover,
  there is no a presentation of quantum mechanics, which does not
  employ complex phases (numbers) in one or another form. Analogous
  parabolic and hyperbolic phases (or characters produced by
  associated hypercomplex numbers, see Section~\ref{sec:hyperc-char})
  lead to classical and hypercomplex mechanics respectively.
\end{rem}

This section clarifies foundations of quantum and classical mechanics.
We recovered the existence of three non-isomorphic models of mechanics
from the representation theory. They were already derived
in~\cites{Hudson04a,Hudson66a} from translation invariant formulation,
that is from the group theory as well.
It also hinted that hyperbolic counterpart is (at least theoretically)
as natural as classical and quantum mechanics are. The approach
provides a framework for a description of aggregate system which have
say both quantum and classical components. This can be used to model
quantum computers with classical terminals~\cite{Kisil09b}.

Remarkably, simultaneously with the work \cite{Hudson66a}
group-invariant axiomatics of geometry leaded
R.I.~Pimenov~\cite{Pimenov65a} to description of \(3^n\) Cayley--Klein
constructions. The connection between group-invariant geometry and respective
mechanics were explored in many works of N.A.~Gromov, see for
example~\cites{Gromov90a,Gromov90b,GromovKuratov05b}. They already
highlighted the r\^ole of three types of hypercomplex units for the
realisation of elliptic, parabolic and hyperbolic geometry and
kinematic.

There is a further connection between representations of the
Heisenberg group and hypercomplex numbers. The symplectomorphism of
phase space%
\index{phase!space}%
\index{space!phase} are also automorphism of the Heisenberg
group~\cite{Folland89}*{\S~1.2}. We recall that the symplectic group
\(\Sp[2]\)%
\index{group!$\Sp[2]$}%
\index{group!$\Sp[2]$|see{also $\SL$}}~\amscite{Folland89}*{\S~1.2} is
isomorphic to the group \(\SL\)~\citelist{ \cite{Lang85}
  \cite{HoweTan92} \cite{Mazorchuk09a}} and provides linear
symplectomorphisms%
\index{symplectic!transformation}%
\index{transformation!symplectic} of the two-dimensional phase space.
It has three types of non-iso\-mor\-phic one-dimensional continuous
subgroups~(\ref{eq:k-subgroup}-\ref{eq:ap-subgroup}) with symplectic
action on the phase space%
\index{phase!space}%
\index{space!phase} illustrated by Fig.~\ref{fig:rotations}.
Hamiltonians, which produce those symplectomorphism, are of
interest~\citelist{\amscite{Wulfman10a}*{\S~3.8} \cite{ATorre08a}
  \cite{ATorre10a}}. An analysis of those Hamiltonians from
Subsection~\ref{sec:correspondence} by means of ladder operators
recreates hypercomplex coefficients as well~\cites{Kisil11a}.

Harmonic oscillators%
\index{harmonic!oscillator}%
\index{oscillator!harmonic}, which we shall use as the main
illustration here, are treated in most textbooks on quantum mechanics%
\index{quantum mechanics}%
\index{quantum mechanics}. This is efficiently done through
creation/annihilation (ladder) operators,
cf.~\S~\ref{sec:correspondence} and~\citelist{\cite{Gazeau09a}
  \cite{BoyerMiller74a}}. The underlying structure is the
representation theory of the Heisenberg%
\index{Heisenberg!group|(}%
\index{group!Heisenberg|(} and symplectic groups%
\index{symplectic!group}%
\index{group!symplectic}~\citelist{\amscite{Lang85}*{\S~VI.2}
  \amscite{MTaylor86}*{\S~8.2} \cite{Howe80b} \cite{Folland89}}.  As
we will see, they are naturally connected with respective hypercomplex
numbers.  As a result we obtain further illustrations to the
Similarity and Correspondence Principle~\ref{pr:simil-corr-principle}.

We work with the simplest case of a particle with only one
degree of freedom. Higher dimensions and the respective group of
symplectomorphisms \(\Sp[2n]\) may require consideration of Clifford
algebras%
\index{Clifford!algebra}%
\index{algebra!Clifford}~\citelist{\cite{Kisil93c}
  \cite{ConstalesFaustinoKrausshar11a} \cite{CnopsKisil97a}
  \cite{GuentherKuzhel10a} \cite{Porteous95}}.

\subsection{The Heisenberg Group and Its Automorphisms}
\label{sec:heis-group-auto}

\subsubsection{The Heisenberg group and induced representations}
\label{sec:heis-group-induc}

Let \((s,x,y)\), where \(s\), \(x\), \(y\in \Space{R}{}\), be an
element of the one-dimensional \emph{Heisenberg group}
\(\Space{H}{1}\)%
\index{Heisenberg!group}%
\index{group!Heisenberg}%
\index{$\Space{H}{1}$, Heisenberg group}~\cites{Folland89,Howe80b}.
Consideration of the general case of \(\Space{H}{n}\) will be similar,
but is beyond the scope of present paper. The group law on
\(\Space{H}{1}\) is given as follows:
\begin{equation}
  \label{eq:H-n-group-law}
  \textstyle
  (s,x,y)\cdot(s',x',y')=(s+s'+\frac{1}{2}\omega(x,y;x',y'),x+x',y+y'), 
\end{equation} 
where the non-commutativity is due to \(\omega\)---the
\emph{symplectic form}%
\index{symplectic!form}%
\index{form!symplectic} on \(\Space{R}{2n}\), which is the central
object of the classical mechanics%
\index{classical mechanics}%
\index{mechanics!classical}~\amscite{Arnold91}*{\S~37}:
\begin{equation}
  \label{eq:symplectic-form}
  \omega(x,y;x',y')=xy'-x'y.
\end{equation}
The Heisenberg group is a non-commutative Lie
group with the centre
\begin{displaymath}
  Z=\{(s,0,0)\in \Space{H}{1}, \ s \in \Space{R}{}\}.
\end{displaymath}
The left shifts
\begin{equation}
  \label{eq:left-right-regular}
  \Lambda(g): f(g') \mapsto f(g^{-1}g')  
\end{equation}
act as a representation of \(\Space{H}{1}\) on a certain linear space
of functions. For example, an action on \(\FSpace{L}{2}(\Space{H}{},dg)\) with
respect to the Haar measure%
\index{Heisenberg!group!invariant measure}%
\index{group!Heisenberg!invariant measure}%
\index{invariant!measure}%
\index{measure!invariant} \(dg=ds\,dx\,dy\) is the \emph{left regular}
representation%
\index{left regular representation}%
\index{representation!left regular}, which is unitary.

The Lie algebra \(\algebra{h}^n\) of \(\Space{H}{1}\) is spanned by
left-(right-)invariant vector fields
\begin{equation}
\textstyle  S^{l(r)}=\pm{\partial_s}, \quad
  X^{l(r)}=\pm\partial_{ x}-\frac{1}{2}y{\partial_s},  \quad
 Y^{l(r)}=\pm\partial_{y}+\frac{1}{2}x{\partial_s}
  \label{eq:h-lie-algebra}
\end{equation}
on \(\Space{H}{1}\) with the Heisenberg \emph{commutator relation}%
\index{Heisenberg!commutator relation}%
\index{commutator relation}
\begin{equation}
  \label{eq:heisenberg-comm}
  [X^{l(r)},Y^{l(r)}]=S^{l(r)} 
\end{equation}
and all other commutators vanishing. We will sometimes omit the
superscript \(l\) for left-invariant field.

We can construct linear representations of \(\Space{H}{1}\) by
induction%
\index{representation!induced!Heisenberg group,of}%
\index{induced!representation!Heisenberg group,of}%
\index{Heisenberg!group!induced representation}%
\index{group!Heisenberg!induced
  representation}~\cite{Kirillov76}*{\S~13} from a character \(\chi\)
of the centre \(Z\). Here we prefer the following one,
cf.~\S~\ref{sec:concl-induc-repr}
and~\citelist{\cite{Kirillov76}*{\S~13} \cite{MTaylor86}*{Ch.~5}}. Let
\(\FSpace[\chi]{F}{2}(\Space{H}{n})\) be the space of functions on
\(\Space{H}{n}\) having the properties:
\begin{equation}
  \label{eq:induced-prop}
  f(gh)=\chi(h)f(g), \qquad \text{ for all } g\in \Space{H}{n},\ h\in Z
\end{equation}
and
\begin{equation}
  \label{eq:L2-condition}
  \int_{\Space{R}{2n}} \modulus{f(0,x,y)}^2dx\,dy<\infty.
\end{equation}
Then \(\FSpace[\chi]{F}{2}(\Space{H}{n})\) is invariant under the left
shifts and those shifts restricted to
\(\FSpace[\chi]{F}{2}(\Space{H}{n})\) make a representation
\(\uir{}{\chi}\) of \(\Space{H}{n}\) induced by \(\chi\).

If the character \(\chi\) is unitary, then the induced representation
is unitary as well. However the representation \(\uir{}{\chi}\) is not
necessarily irreducible. Indeed, left shifts are commuting with the
right action of the group. Thus any subspace of null-solutions of a
linear combination \(aS+\sum_{j=1}^n (b_jX_j+c_jY_j)\) of
left-invariant vector fields is left-invariant and we can restrict
\(\uir{}{\chi}\) to this subspace. The left-invariant differential
operators define analytic condition for functions, cf. Cor.~\ref{co:cauchy-riemann}.
\begin{example}
  The function \(f_0(s,x,y)=e^{\rmi \myh s -\myh(x^2 +y^2)/4}\), where
  \(\myh=2\pi\myhbar\), belongs to
  \(\FSpace[\chi]{F}{2}(\Space{H}{n})\) for the character
  \(\chi(s)=e^{\rmi \myh s}\). It is also a null solution for all the
  operators \(X_j-\rmi Y_j\). The closed linear span of functions
  \(f_g=\Lambda(g) f_0\) is invariant under left shifts and provide a
  model for Fock--Segal--Bargmann (FSB) type representation of the Heisenberg
  group%
  \index{Fock--Segal--Bargmann!representation}%
  \index{representation!Fock--Segal--Bargmann}%
  \index{Heisenberg!group!Fock--Segal--Bargmann representation}%
  \index{group!Heisenberg!Fock--Segal--Bargmann representation}, which
  will be considered below.
\end{example}

\subsubsection{Symplectic Automorphisms of the Heisenberg Group}
\label{sec:sympl-autom-heis}

The group of outer automorphisms of \(\Space{H}{1}\), which trivially
acts on the centre of \(\Space{H}{1}\), is the symplectic group
\(\Sp[2]\)\index{$\Sp[2]$}. It is the group of symmetries of the
symplectic form \(\omega\)
in~\eqref{eq:H-n-group-law}~\citelist{\amscite{Folland89}*{Thm.~1.22}
  \amscite{Howe80a}*{p.~830}}. The symplectic group is isomorphic to
\(\SL\) considered in the first half of this work. The
explicit action of \(\Sp[2]\) on the Heisenberg group is:
\begin{equation}
  \label{eq:sympl-auto}
  g: h=(s,x,y)\mapsto g(h)=(s,x',y'), 
\end{equation}
where 
\begin{displaymath}
  g=\begin{pmatrix}
    a&b\\
    c&d
  \end{pmatrix}\in\Sp[2], \quad\text{ and }\quad
  \begin{pmatrix}
    x'\\y'
  \end{pmatrix}
  =\begin{pmatrix}
    a&b\\
    c&d
  \end{pmatrix}
  \begin{pmatrix}
    x\\y
  \end{pmatrix}.
\end{displaymath}
The Shale--Weil theorem~\citelist{\amscite{Folland89}*{\S~4.2}
  \amscite{Howe80a}*{p.~830}} states that any representation
\(\uir{}{\myhbar}\) of the Heisenberg groups generates a unitary
\emph{oscillator}%
\index{oscillator!representation}%
\index{representation!oscillator!} (or \emph{metaplectic}) representation%
\index{representation!metaplectic|see{oscillator representation}}%
\index{metaplectic representation|see{oscillator representation}}%
\index{representation!Shale--Weil|see{oscillator representation}}%
\index{Shale--Weil representation|see{oscillator representation}}
\(\uir{\text{SW}}{\myhbar}\) of the \(\widetilde{\mathrm{Sp}}(2)\),
the two-fold cover of the symplectic group~\amscite{Folland89}*{Thm.~4.58}.

We can consider the semidirect product
\(G=\Space{H}{1}\rtimes\widetilde{\mathrm{Sp}}(2)\) with the standard group law:
\begin{equation}
  \label{eq:schrodinger-group}
  (h,g)*(h',g')=(h*g(h'),g*g'), \qquad \text{where } 
  h,h'\in\Space{H}{1}, \quad g,g'\in\widetilde{\mathrm{Sp}}(2),
\end{equation}
and the stars denote the respective group operations while the action
\(g(h')\) is defined as the composition of the projection map
\(\widetilde{\mathrm{Sp}}(2)\rightarrow {\mathrm{Sp}}(2)\) and the
action~\eqref{eq:sympl-auto}. This group is sometimes called the
\emph{Schr\"odinger group}%
\index{Schr\"odinger!group}%
\index{group!Schr\"odinger} and it is known as the maximal kinematical
invariance group of both the free Schr\"odinger equation and the
quantum harmonic oscillator~\cite{Niederer73a}. This group is of
interest not only in quantum mechanics but also in
optics\index{optics}~\cites{ATorre10a,ATorre08a}. 
The Shale--Weil theorem allows us to expand any representation
\(\uir{}{\myhbar}\) of the Heisenberg group to the representation
\(\uir{2}{\myhbar}=\uir{}{\myhbar}\oplus\uir{\text{SW}}{\myhbar}\) of the
group \(G\).

Consider the Lie algebra \(\algebra{sp}_2\) of the group \(\Sp[2]\).
We again use the basis \(A\), \(B\),
\(Z\)~\eqref{eq:sl2-basis} with
commutators~\eqref{eq:sl2-commutator}.  Vectors \(Z\), \(B-Z/2\) and
\(B\) are generators of the one-parameter subgroups \(K\), \(N'\) and
\(\Aprime\) (\ref{eq:k-subgroup}--\ref{eq:ap-subgroup}) respectively.
Furthermore we can consider the basis \(\{S, X, Y, A, B, Z\}\) of the
Lie algebra \(\algebra{g}\) of the Lie group
\(G=\Space{H}{1}\rtimes\widetilde{\mathrm{Sp}}(2)\). All non-zero
commutators besides those already listed in~\eqref{eq:heisenberg-comm}
and~\eqref{eq:sl2-commutator} are: 
\begin{align}
  \label{eq:cross-comm}
  [A,X]&=\textstyle\frac{1}{2}X,&
  [B,X]&=\textstyle-\frac{1}{2}Y,&
  [Z,X]&=Y;\\
  \label{eq:cross-comm1}
  [A,Y]&=\textstyle-\frac{1}{2}Y,&
  [B,Y]&=\textstyle-\frac{1}{2}X,&
  [Z,Y]&=-X.
\end{align}
Of course, there is the derived form of the Shale--Weil representation
for \(\algebra{g}\). It can often be explicitly written in contrast
to the Shale--Weil representation.
\begin{example}
  Let \(\uir{}{\myhbar}\) be the Schr\"odinger
  representation%
  \index{representation!Heisenberg group!Schr\"odinger}%
  \index{Schr\"odinger!representation}%
  \index{representation!Schr\"odinger}~\amscite{Folland89}*{\S~1.3} of \(\Space{H}{1}\) in
  \(\FSpace{L}{2}(\Space{R}{})\), that is \amscite{Kisil10a}*{(3.5)}:
  \begin{displaymath}
    [\uir{}{\chi}(s,x,y) f\,](q)=e^{2\pi\rmi\myhbar (s-xy/2)
      +2\pi\rmi x q}\,f(q-\myhbar y).  
  \end{displaymath}
  Thus the action of the derived representation on the Lie algebra
  \(\algebra{h}_1\) is:
  \begin{equation}
    \label{eq:schroedinger-rep-conf-der}
    \uir{}{\myhbar}(X)=2\pi\rmi q,\qquad \uir{}{\myhbar}(Y)=-\myhbar \frac{d}{dq},
    \qquad
    \uir{}{\myhbar}(S)=2\pi\rmi\myhbar I.
  \end{equation}
  Then the associated Shale--Weil representation of \(\Sp[2]\) in
  \(\FSpace{L}{2}(\Space{R}{})\) has the
  derived action,  cf.~\citelist{\amscite{ATorre08a}*{(2.2)} \amscite{Folland89}*{\S~4.3}}:
  \begin{equation}
    \label{eq:shale-weil-der}
    \uir{\text{SW}}{\myhbar}(A) =-\frac{q}{2}\frac{d}{dq}-\frac{1}{4},\quad
    \uir{\text{SW}}{\myhbar}(B)=-\frac{\myhbar\rmi}{8\pi}\frac{d^2}{dq^2}-\frac{\pi\rmi q^2}{2\myhbar},\quad
    \uir{\text{SW}}{\myhbar}(Z)=\frac{\myhbar\rmi}{4\pi}\frac{d^2}{dq^2}-\frac{\pi\rmi q^2}{\myhbar}.
  \end{equation}
  We can verify commutators~\eqref{eq:heisenberg-comm} and
  \eqref{eq:sl2-commutator}, \eqref{eq:cross-comm1} for
  operators~(\ref{eq:schroedinger-rep-conf-der}--\ref{eq:shale-weil-der}).
  It is also obvious that in this representation the following
  algebraic relations hold:
  \begin{eqnarray}
    \label{eq:quadratic-A}
    \qquad\uir{\text{SW}}{\myhbar}(A) &=&
    \frac{\rmi}{4\pi\myhbar}(\uir{}{\myhbar}(X)\uir{}{\myhbar}(Y)-{\textstyle\frac{1}{2}}\uir{}{\myhbar}(S))\\
    &=&\frac{\rmi}{8\pi\myhbar}(\uir{}{\myhbar}(X)\uir{}{\myhbar}(Y)+\uir{}{\myhbar}(Y)\uir{}{\myhbar}(X) ), \nonumber\\ 
    \label{eq:quadratic-B}
    \uir{\text{SW}}{\myhbar}(B) &=&
    \frac{\rmi}{8\pi\myhbar}(\uir{}{\myhbar}(X)^2-\uir{}{\myhbar}(Y)^2), \\
    \label{eq:quadratic-Z}
    \uir{\text{SW}}{\myhbar}(Z)
    &=&\frac{\rmi}{4\pi\myhbar}(\uir{}{\myhbar}(X)^2+\uir{}{\myhbar}(Y)^2). 
  \end{eqnarray}
  Thus it is common in quantum optics to name \(\algebra{g}\) as a Lie
  algebra with  quadratic generators%
  \index{generator!quadratic}%
  \index{quadratic!generator}, see~\amscite{Gazeau09a}*{\S~2.2.4}.
\end{example}
Note that \(\uir{\text{SW}}{\myhbar}(Z)\) is the Hamiltonian of the
harmonic oscillator%
\index{harmonic!oscillator}%
\index{oscillator!harmonic} (up to a factor). Then we can consider
\(\uir{\text{SW}}{\myhbar}(B)\) as the Hamiltonian of a repulsive
(hyperbolic) oscillator. The operator
\(\uir{\text{SW}}{\myhbar}(B-Z/2)=\frac{\myhbar\rmi}{4\pi}\frac{d^2}{dq^2}\)
is the parabolic analog. A graphical representation of all three
transformations defined by those Hamiltonian is given in
Fig.~\ref{fig:rotations} and a further discussion of these
Hamiltonians can be found in~\amscite{Wulfman10a}*{\S~3.8}.

An important observation, which is often missed, is that the
three linear symplectic transformations are unitary rotations in the
corresponding hypercomplex algebra, cf.~\amscite{Kisil09c}*{\S~3}. This
means, that the symplectomorphisms generated by operators \(Z\),
\(B-Z/2\), \(B\) within time \(t\) coincide with the
multiplication of hypercomplex number \(q+\alli p\) by \(e^{\alli
  t}\), see Subsection~\ref{sec:hyperc-char} and
Fig.~\ref{fig:rotations}, which is just another illustration of the
Similarity and Correspondence Principle~\ref{pr:simil-corr-principle}.
\begin{example}
  There are many advantages of considering representations of the
  Heisenberg group on the phase
  space~\citelist{\amscite{Howe80b}*{\S~1.7}
    \amscite{Folland89}*{\S~1.6} \cite{deGosson08a}}. A convenient
  expression for Fock--Segal--Bargmann%
  \index{Fock--Segal--Bargmann!representation}%
  \index{representation!Fock--Segal--Bargmann}%
  \index{FSB!representation|see{Fock--Segal--Bargmann representation}}%
  \index{representation!FSB|see{Fock--Segal--Bargmann representation}}
  (FSB) representation on the phase space%
  \index{phase!space}%
  \index{space!phase} is,
  cf.~\S~\ref{sec:schr-segal-bargm} and~\citelist{\cite{Kisil02e}*{(2.9)}
  \cite{deGosson08a}*{(1)}}:
  \begin{equation}
    \label{eq:stone-inf}
    \textstyle
    [\uir{}{F}(s,x,y) f] (q,p)=
    e^{-2\pi\rmi(\myhbar s+qx+py)}
    f \left(q-\frac{\myhbar}{2} y, p+\frac{\myhbar}{2} x\right).
  \end{equation}
  Then the derived representation of \(\algebra{h}_1\) is:
  \begin{equation}
    \label{eq:fock-rep-conf-der-par1}
    \textstyle
    \uir{}{F}(X)=-2\pi\rmi q+\frac{\myhbar}{2}\partial_{p},\qquad
    \uir{}{F}(Y)=-2\pi\rmi p-\frac{\myhbar}{2}\partial_{q},
    \qquad
    \uir{}{F}(S)=-2\pi\rmi\myhbar I.
  \end{equation}
  This produces the derived form of the Shale--Weil representation:
  \begin{equation}
    \label{eq:shale-weil-der-ell}
    \textstyle
    \uir{\text{SW}}{F}(A) =\frac{1}{2}\left(q\partial_{q}-p\partial_{p}\right),\quad
    \uir{\text{SW}}{F}(B)=-\frac{1}{2}\left(p\partial_{q}+q\partial_{p}\right),\quad
    \uir{\text{SW}}{F}(Z)=p\partial_{q}-q\partial_{p}.
  \end{equation}
  Note that this representation does not contain the parameter
  \(\myhbar\) unlike the equivalent
  representation~\eqref{eq:shale-weil-der}. Thus the FSB model%
  \index{Fock--Segal--Bargmann!space}%
  \index{space!Fock--Segal--Bargmann} explicitly shows the equivalence
  of \(\uir{\text{SW}}{\myhbar_1}\) and \(\uir{\text{SW}}{\myhbar_2}\)
  if \(\myhbar_1 \myhbar_2>0\)~\amscite{Folland89}*{Thm.~4.57}.

  As we will also see below the FSB-type representations%
  \index{Fock--Segal--Bargmann!representation}%
  \index{representation!Fock--Segal--Bargmann} in
  hypercomplex numbers produce almost the same Shale--Weil
  representations.
\end{example}

\subsection{p-Mechanic Formalism}
Here we briefly outline a formalism~\cites{Kisil96a,Prezhdo-Kisil97,Kisil00a,BrodlieKisil03a,Kisil02e}, which allows to unify quantum and
classical mechanics.\index{p-mechanics}

\subsubsection{Convolutions (Observables) on $\Space{H}{n}$ and Commutator}
\label{sec:conv-algebra-hg}

Using a invariant measure \(dg=ds\,dx\,dy\) on \(\Space{H}{n}\) we can
define the convolution\index{convolution} of two functions: 
\begin{eqnarray}
  (k_1 * k_2) (g) &=& \int_{\Space{H}{n}} k_1(g_1)\,
  k_2(g_1^{-1}g)\,dg_1  .
  \label{eq:de-convolution}
\end{eqnarray}
This is a non-commutative operation, which is meaningful for functions
from various spaces including \(\FSpace{L}{1}(\Space{H}{n},dg)\), the
Schwartz space \(\FSpace{S}{}\) and many classes of distributions,
which form algebras under convolutions.  Convolutions on
\(\Space{H}{n}\) are used as \emph{observables}%
\index{observable}%
\index{p-mechanics!observable} in
\(p\)-mechanic~\cites{Kisil96a,Kisil02e}.

A unitary representation \(\uir{}{} \) of \(\Space{H}{n}\) extends
 to \(\FSpace{L}{1}(\Space{H}{n} ,dg)\) by the formula:
\begin{equation}
  \label{eq:rho-extended-to-L1}
  \uir{}{} (k) = \int_{\Space{H}{n}} k(g)\uir{}{}  (g)\,dg .
\end{equation}
This is also an algebra homomorphism of convolutions to linear
operators. 

For a dynamics of observables we need inner \emph{derivations}%
\index{inner!derivation}%
\index{derivation!inner} \(D_k\) of
the convolution algebra \(\FSpace{L}{1}(\Space{H}{n})\), which are
given by the \emph{commutator}\index{commutator}:
\begin{eqnarray}
  \qquad D_k: f \mapsto [k,f]&=&k*f-f*k   \label{eq:commutator}
   \\  &=&
  \int_{\Space{H}{n}} k(g_1)\left( f(g_1^{-1}g)-f(gg_1^{-1})\right)\,dg_1
, \quad f,k\in\FSpace{L}{1}(\Space{H}{n}).
\nonumber 
\end{eqnarray}

To describe dynamics of a time-dependent observable \(f(t,g)\) we use
the universal equation%
\index{p-mechanics!dynamic equation}%
\index{equation!dynamics in p-mechanics}, cf.~\cites{Kisil94d,Kisil96a}:
\begin{equation}
  \label{eq:universal}
  S\dot{f}=[H,f],
\end{equation}
where \(S\) is the left-invariant vector
field~\eqref{eq:h-lie-algebra} generated by the centre of
\(\Space{H}{n}\). The presence of operator \(S\) fixes the
dimensionality of both sides of the equation~\eqref{eq:universal} if
the observable \(H\) (Hamiltonian) has the dimensionality of
energy~\cite{Kisil02e}*{Rem~4.1}. If we apply a right inverse
\(\anti\) of \(S\) to both sides of the equation~\eqref{eq:universal}
we obtain the equivalent equation
\begin{equation}
  \label{eq:universal-bracket}
  \dot{f}=\ub{H}{f},
\end{equation}
based on the universal bracket \(\ub{k_1}{k_2}=k_1*\anti k_2-k_2*\anti
k_1\)~\cite{Kisil02e}. 
\begin{example}[Harmonic oscillator]
  \label{ex:p-harmonic}
  Let \(H=\frac{1}{2} (mk^2 q^2 + \frac{1}{m}p^2)\) be the
  Hamiltonian of a one-dimensional harmonic oscillator%
  \index{harmonic!oscillator|indef}%
  \index{oscillator!harmonic|indef}, where \(k\) is a constant
  frequency and \(m\) is a constant mass.  Its \emph{p-mechanisation}%
  \index{p-mechanisation} will be the second order differential
  operator on \(\Space{H}{n}\)~\cite{BrodlieKisil03a}*{\S~5.1}:
  \begin{displaymath}
    \textstyle
    H=\frac{1}{2} (mk^2 X^2
    + \frac{1}{m}Y^2),
  \end{displaymath}
  where we dropped sub-indexes of vector
  fields~\eqref{eq:h-lie-algebra} in one dimensional setting. We can
  express the commutator as a difference between the left and the
  right action of the vector fields:
  \begin{displaymath}
    \textstyle
    [H,f]=\frac{1}{2} (mk^2 ((X^{r})^2-(X^{l})^2)
    + \frac{1}{m}((Y^{r})^2-(Y^{l})^2))f.
  \end{displaymath}
   Thus the equation~\eqref{eq:universal} becomes~\cite{BrodlieKisil03a}*{(5.2)}:
   \begin{equation}
     \label{eq:p-harm-osc-dyn}
     \frac{\partial }{\partial s}\dot{f}= \frac{\partial }{\partial s}
    \left(m k^2 y\frac{\partial}{\partial x}-\frac{1}{m} x
      \frac{\partial}{\partial y} \right) f. 
   \end{equation}
   Of course, the derivative \(\frac{\partial }{\partial s}\) can be
   dropped from both sides of the equation and the general solution
   is found to be:
   \begin{equation}
     \label{eq:p-harm-sol}
     \textstyle
     f(t;s,x,y)  =  f_0\left(s, x\cos(k
    t) +
      m k y\sin( k t), 
      -\frac{x}{mk} \sin(k t) + y\cos (k t)\right),
   \end{equation}
   where \(f_0(s,x,y)\) is the initial value of an observable on \(\Space{H}{n}\).
\end{example}
\begin{example}[Unharmonic oscillator] 
  \label{ex:p-unharmonic}
  We consider unharmonic
  oscillator%
  \index{unharmonic!oscillator}%
  \index{oscillator!unharmonic} with cubic potential, see~\cite{CalzettaVerdaguer06a} and
  references therein:
  \begin{equation}
    \label{eq:unharmonic-hamiltonian}
    H=\frac{mk^2}{2} q^2+\frac{\lambda}{6} q^3
  + \frac{1}{2m}p^2.
  \end{equation}
  Due to the absence of non-commutative products p-mechanisation is straightforward: 
  \begin{displaymath}
    H=\frac{mk^2}{2}  X^2+\frac{\lambda}{6} X^3
    + \frac{1}{m}Y^2.
  \end{displaymath}
  Similarly to the harmonic case the dynamic equation, after
  cancellation of \(\frac{\partial }{\partial s}\) on both sides,
  becomes:
   \begin{equation}
     \label{eq:p-unharm-osc-dyn}
     \dot{f}=     \left(m k^2 y\frac{\partial}{\partial x}
       +\frac{\lambda}{6}\left(3y\frac{\partial^2}{\partial x^2} 
         +\frac{1}{4}y^3\frac{\partial^2}{\partial s^2}\right)-\frac{1}{m} x
      \frac{\partial}{\partial y} \right) f. 
   \end{equation}
Unfortunately, it cannot be solved analytically as easy as in the
harmonic case.
\end{example}

\subsubsection{States and Probability}
\label{sec:states-probability}

Let an observable \(\uir{}{}(k)\)~\eqref{eq:rho-extended-to-L1} is
defined by a kernel \(k(g)\) on the Heisenberg group and its
representation \(\uir{}{}\) at a Hilbert space \(\mathcal{H}\). A
\emph{state}%
\index{state}%
\index{p-mechanics!state} on the convolution algebra is given by a vector
\(v\in\mathcal{H}\). A simple calculation:
\begin{eqnarray*}
  \scalar[\mathcal{H}]{\uir{}{}(k)v}{v}&=& \scalar[\mathcal{H}]{\int_{\Space{H}{n}} k(g)
    \uir{}{}(g)v\,dg}{v}\\
  &=& \int_{\Space{H}{n}} k(g) \scalar[\mathcal{H}]{\uir{}{}(g)v}{v}dg\\
  &=& \int_{\Space{H}{n}} k(g) \overline{\scalar[\mathcal{H}]{v}{\uir{}{}(g)v}}\,dg
\end{eqnarray*}
can be restated as:
\begin{displaymath}
  \scalar[\mathcal{H}]{\uir{}{}(k)v}{v}=\scalar[]{k}{l}, \qquad \text{where} \quad
  l(g)=\scalar[\mathcal{H}]{v}{\uir{}{}(g)v}.
\end{displaymath}
Here the left-hand side contains the inner product on
\(\mathcal{H}\), while the right-hand side uses a skew-linear pairing
between functions on \(\Space{H}{n}\) based on the Haar measure
integration. In other words we obtain,
cf.~\cite{BrodlieKisil03a}*{Thm.~3.11}:
\begin{prop}
  \label{pr:state-functional}
  A state defined by a vector \(v\in\mathcal{H}\) coincides with the
  linear functional given by the wavelet transform
  \begin{equation}
    \label{eq:kernel-state}
    l(g)=\scalar[\mathcal{H}]{v}{\uir{}{}(g)v}
  \end{equation}
  of \(v\) used as the mother wavelet as well.
\end{prop}
The addition of vectors in \(\mathcal{H}\) implies the following
operation on states:
\begin{eqnarray}
  \scalar[\mathcal{H}]{v_1+v_2}{\uir{}{}(g)(v_1+v_2)}&=&
    \scalar[\mathcal{H}]{v_1}{\uir{}{}(g)v_1}
    +\scalar[\mathcal{H}]{v_2}{\uir{}{}(g)v_2}\nonumber \\
    &&{}+
    \scalar[\mathcal{H}]{v_1}{\uir{}{}(g)v_2}
   + \overline{\scalar[\mathcal{H}]{v_1}{\uir{}{}(g^{-1})v_2}}
   \label{eq:kernel-add}
\end{eqnarray}
The last expression can be conveniently rewritten for kernels of the
functional as 
\begin{equation}
  \label{eq:addition-functional}
  l_{12}=l_1+l_2+2 A\sqrt{l_1l_2}
\end{equation}
for some real number \(A\). This formula is behind the contextual law
of addition of conditional probabilities%
\index{probability!quantum}%
\index{quantum!probability}~\cite{Khrennikov01a} and will
be illustrated below. Its physical interpretation is an interference,%
\index{interference} say, from two slits. Despite of a common belief,
the mechanism of such interference can be both causal\index{causal}
and local, see~\citelist{\cite{Kisil01c} \cite{KhrenVol01}}.

\subsection{Elliptic characters and Quantum Dynamics}
\label{sec:ellipt-char-moyal}
In this subsection we consider the representation \(\uir{}{\myh}\) of
\(\Space{H}{n}\) induced by the elliptic character
\(\chi_\myh(s)=e^{\rmi\myh s}\) in complex numbers parametrised by
\(\myh\in\Space{R}{}\). We also use the convenient agreement
\(\myh=2\pi\myhbar\) borrowed from physical literature.

\subsubsection{Fock--Segal--Bargmann and Schr\"odinger Representations}
\label{sec:schr-segal-bargm}
The realisation of \(\uir{}{\myh}\) by the left
shifts~\eqref{eq:left-right-regular} on
\(\FSpace[\myh]{L}{2}(\Space{H}{n})\) is rarely used in quantum
mechanics. Instead two unitary equivalent forms are more common: the
Schr\"odinger%
\index{Schr\"odinger!representation}%
\index{representation!Schr\"odinger} and Fock--Segal--Bargmann (FSB) representations.%
\index{Fock--Segal--Bargmann!representations|indef}%
\index{representations!Fock--Segal--Bargmann|indef}

The FSB representation can be obtained from the orbit
method of Kirillov~\cite{Kirillov94a}. It allows spatially separate
irreducible components of the left regular representation, each of
them become located on the orbit of the co-adjoint representation,
see~\citelist{\cite{Kisil02e}*{\S~2.1} \cite{Kirillov94a}} for
details, we only present a brief summary here.

We identify \(\Space{H}{n}\) and its Lie algebra \(\algebra{h}_n\)
through the exponential map~\cite{Kirillov76}*{\S~6.4}. The dual
\(\algebra{h}_n^*\) of \(\algebra{h}_n\) is presented by the Euclidean
space \(\Space{R}{2n+1}\) with coordinates \((\myhbar,q,p)\).  The
pairing \(\algebra{h}_n^*\) and \(\algebra{h}_n\) given by
\begin{displaymath}
  \scalar{(s,x,y)}{(\myhbar,q,p)}=\myhbar s + q \cdot x+p\cdot y.
\end{displaymath}
This pairing defines the Fourier
transform \(\hat{\ }: \FSpace{L}{2}(\Space{H}{n})\rightarrow
\FSpace{L}{2}(\algebra{h}_n^*)\) given by~\cite{Kirillov99}*{\S~2.3}:
\begin{equation}
  \label{eq:fourier-transform}
  \hat{\phi}(F)=\int_{\algebra{h}^n} \phi(\exp X) 
  e^{-2\pi\rmi  \scalar{X}{F}}\,dX \qquad \textrm{ where }
  X\in\algebra{h}^n,\ F\in\algebra{h}_n^*. 
\end{equation}
For a fixed \(\myhbar\) the left regular
representation~\eqref{eq:left-right-regular} is mapped by the Fourier
transform to the FSB type representation~\eqref{eq:stone-inf}.  The
collection of points \((\myhbar,q,p)\in\algebra{h}_n^*\) for a fixed
\(\myhbar\) is naturally identified with the \emph{phase space}%
\index{phase!space}%
\index{space!phase} of the system.
\begin{rem}
  It is possible to identify the case of \(\myhbar=0\) with classical
  mechanics~\cite{Kisil02e}. Indeed, a substitution of the zero value of \(\myhbar\)
  into~\eqref{eq:stone-inf} produces the commutative representation:
  \begin{equation}
    \label{eq:commut-repres}
      \uir{}{0}(s,x,y): f (q,p) \mapsto 
      e^{-2\pi\rmi(qx+py)}
      f \left(q, p\right).
  \end{equation}
  It can be decomposed into the direct integral of one-dimensional
  representations parametrised by the points \((q,p)\) of the phase
  space. The classical mechanics%
  \index{classical mechanics}%
  \index{mechanics!classical}, including the Hamilton equation, can
  be recovered from those representations~\cite{Kisil02e}. However 
  the condition \(\myhbar=0\) (as well as the \emph{semiclassical limit}%
  \index{semiclassical!limit}%
  \index{limit!semiclassical} \(\myhbar\rightarrow 0\)) is
  not completely physical. Commutativity (and subsequent relative
  triviality) of those representation is the main reason why they are
  oftenly neglected. The commutativity can be outweighed by special
  arrangements, e.g. an antiderivative~\cite{Kisil02e}*{(4.1)}, but the
  procedure is not straightforward, see discussion in~\citelist{\cite{Kisil05c}
    \cite{AgostiniCapraraCiccotti07a}
    \cite{Kisil09a}}. A direct approach using dual
  numbers will  be shown below, cf. Rem.~\ref{re:hamilton-from-nc}.
\end{rem}

To recover the Schr\"odinger representation we use notations and
technique of induced representations from
\S~\ref{sec:concl-induc-repr}, see also~\cite{Kisil98a}*{Ex.~4.1}.
The subgroup \(H=\{(s,0,y)\such s\in\Space{R}{},
y\in\Space{R}{n}\}\subset\Space{H}{n}\) defines the homogeneous space
\(X=G/H\), which coincides with \(\Space{R}{n}\) as a manifold. The
natural projection \(\mathbf{p}:G\rightarrow X\) is \(\mathbf{p}(s,x,y)=x\) and its left
inverse \(\mathbf{s}:X\rightarrow G\) can be as simple as \(\mathbf{s}(x)=(0,x,0)\).
For the map \(\mathbf{r}:G\rightarrow H\), \(\mathbf{r}(s,x,y)=(s-xy/2,0,y)\) we have
the decomposition
\begin{displaymath}
  (s,x,y)=\mathbf{s}(p(s,x,y))*\mathbf{r}(s,x,y)=(0,x,0)*(s-\textstyle\frac{1}{2}xy,0,y).
\end{displaymath}
For a character
\(\chi_{\myh}(s,0,y)=e^{\rmi\myh s}\) of \(H\) the lifting
\(\oper{L}_\chi: \FSpace{L}{2}(G/H) \rightarrow
\FSpace[\chi]{L}{2}(G)\) is as follows:
\begin{displaymath}
  [\oper{L}_\chi f](s,x,y)=\chi_{\myh}(\mathbf{r}(s,x,y))\, 
  f(\mathbf{p}(s,x,y))=e^{\rmi\myh (s-xy/2)}f(x).  
\end{displaymath}
Thus the representation \(\uir{}{\chi}(g)=\oper{P}\circ\Lambda
(g)\circ\oper{L}\) becomes:
\begin{equation}
  \label{eq:schroedinger-rep}
  [\uir{}{\chi}(s',x',y') f](x)=e^{-2\pi\rmi\myhbar (s'+xy'-x'y'/2)}\,f(x-x').  
\end{equation}
After the Fourier transform \(x\mapsto q\) we get the Schr\"odinger
representation%
\index{Schr\"odinger!representation|indef}%
\index{representation!Schr\"odinger|indef} on the \emph{configuration
  space}%
\index{configuration!space}%
\index{space!configuration}:
\begin{equation}
  \label{eq:schroedinger-rep-conf}
  [\uir{}{\chi}(s',x',y') \hat{f}\,](q)=e^{-2\pi\rmi\myhbar (s'+x'y'/2)
    -2\pi\rmi x' q}\,\hat{f}(q+\myhbar y').  
\end{equation}
Note that this again turns into a commutative representation
(multiplication by an unimodular function) if \(\myhbar=0\). To get
the full set of commutative representations in this way we need to use the
character \(\chi_{(\myh,p)}(s,0,y)=e^{2\pi\rmi(\myhbar+ py)}\) in the
above consideration. 

\subsubsection{Commutator and the Heisenberg Equation}
\label{sec:comm-heis-equat}

The property~\eqref{eq:induced-prop} of
\(\FSpace[\chi]{F}{2}(\Space{H}{n})\) implies that the restrictions of
two operators \(\uir{}{\chi} (k_1)\) and \(\uir{}{\chi} (k_2)\) to
this space are equal if
\begin{displaymath}
  \int_{\Space{R}{}} k_1(s,x,y)\,\chi(s)\, ds = \int_{\Space{R}{}} k_2(s,x,y)\,\chi(s)\,ds.
\end{displaymath}
In other words, for a character \(\chi(s)=e^{2\pi\rmi \myhbar s}\) the
operator \(\uir{}{\chi} (k)\) depends only on
\begin{displaymath}
  \hat{k}_s(\myhbar,x,y)=\int_{\Space{R}{}} k(s,x,y)\,e^{-2\pi\rmi \myhbar s}\,ds,
\end{displaymath}
which is the partial Fourier transform \(s\mapsto \myhbar\) of
\(k(s,x,y)\). The restriction to \(\FSpace[\chi]{F}{2}(\Space{H}{n})\)
of the composition formula for convolutions is~\cite{Kisil02e}*{(3.5)}:
\begin{equation}
  \label{eq:composition-ell}
  (k'*k)\hat{_s}
  =
  \int_{\Space{R}{2n}} e^{ {\rmi \myh}{}(xy'-yx')/2}\, \hat{k}'_s(\myhbar ,x',y')\,
 \hat{k}_s(\myhbar ,x-x',y-y')\,dx'dy'. 
\end{equation}
Under the Schr\"odinger representation~\eqref{eq:schroedinger-rep-conf} the
convolution~\eqref{eq:composition-ell} defines a rule for composition
of two pseudo-differential operators (PDO) in the Weyl
calculus~\citelist{\cite{Howe80b} \cite{Folland89}*{\S~2.3}}.

Consequently the representation~\eqref{eq:rho-extended-to-L1} of
commutator~\eqref{eq:commutator} depends only on its partial Fourier
transform~\cite{Kisil02e}*{(3.6)}: 
\begin{eqnarray}
  [k',k]\hat{_s}
  &=&   2 \rmi  \int_{\Space{R}{2n}}\!\! \sin(\textstyle\frac{\myh}{2}
  (xy'-yx'))\,\label{eq:repres-commutator}\\
   && \qquad\times 
  \hat{k}'_s(\myhbar, x', y')\,
  \hat{k}_s(\myhbar, x-x', y-y')\,dx'dy'. \nonumber 
\end{eqnarray}
Under the Fourier transform~\eqref{eq:fourier-transform} this
commutator is exactly the \emph{Moyal bracket}%
\index{Moyal bracket}%
\index{bracket!Moyal}~\cite{Zachos02a} for of
\(\hat{k}'\) and \(\hat{k}\) on the phase space.

For observables in the space
\(\FSpace[\chi]{F}{2}(\Space{H}{n})\) the action of \(S\) is reduced to
multiplication
, e.g. for \(\chi(s)=e^{\rmi \myh s}\) the action of \(S\) is
multiplication by \(\rmi \myh\). Thus the
equation~\eqref{eq:universal} reduced to the space
\(\FSpace[\chi]{F}{2}(\Space{H}{n})\) becomes the Heisenberg type
equation~\cite{Kisil02e}*{(4.4)}:
\begin{equation}
  \label{eq:heisenberg-eq}
  \dot{f}=\frac{1}{\rmi\myh}  [H,f]\hat{_s},
\end{equation}
based on the above bracket~\eqref{eq:repres-commutator}.  The
Schr\"odinger representation~\eqref{eq:schroedinger-rep-conf} transforms
this equation to the original Heisenberg equation%
\index{Heisenberg!equation}%
\index{equation!Heisenberg}.

\begin{example}
  \label{ex:quntum-oscillators}
  \begin{enumerate}
  \item 
    \label{it:q-harmonic}
    Under the Fourier transform \((x,y)\mapsto(q,p)\) the
    p-dynamic equation~\eqref{eq:p-harm-osc-dyn} of the harmonic
    oscillator%
    \index{harmonic!oscillator}%
    \index{oscillator!harmonic} becomes: 
    \begin{equation}
      \label{eq:harmic-osc-dyn}
      \dot{f}=     \left(m k^2 q\frac{\partial}{\partial p}-\frac{1}{m} p
      \frac{\partial}{\partial q} \right) f. 
    \end{equation}
    The same transform creates its solution out
    of~\eqref{eq:p-harm-sol}.
  \item 
    \label{it:q-unharmonic}
    Since \(\frac{\partial}{\partial s}\) acts on
    \(\FSpace[\chi]{F}{2}(\Space{H}{n})\) as multiplication by \(\rmi
    \myhbar\), the quantum representation of unharmonic dynamics%
    \index{unharmonic!oscillator}%
    \index{oscillator!unharmonic}
    equation~\eqref{eq:p-unharm-osc-dyn} is:
    \begin{equation}
      \label{eq:q-unhar-dyn}
      \dot{f}=     \left(m k^2 q\frac{\partial}{\partial
          p}+\frac{\lambda}{6}\left(3q^2\frac{\partial}{\partial p}
          -\frac{\myhbar^2}{4}\frac{\partial^3}{\partial p^3}\right)-\frac{1}{m}
        p \frac{\partial}{\partial q} \right) f. 
    \end{equation}
    This is exactly the equation for the Wigner function obtained
    in~\cite{CalzettaVerdaguer06a}*{(30)}. 
  \end{enumerate}
\end{example}

\subsubsection{Quantum Probabilities}
\label{sec:quantum-probabilities}
For the elliptic character \(\chi_\myh(s)=e^{\rmi\myh s }\) we can use
the Cauchy--Schwartz inequality to demonstrate that the real number
\(A\) in the identity~\eqref{eq:addition-functional} is between \(-1\)
and \(1\). Thus we can put \(A=\cos \alpha\) for some angle (phase)
\(\alpha\) to get the formula for counting quantum
probabilities%
\index{probability!quantum}%
\index{quantum!probability}, cf.~\cite{Khrennikov03a}*{(2)}:
\begin{equation}
  \label{eq:addition-functional-ell}
  l_{12}=l_1+l_2+2 \cos\alpha \,\sqrt{l_1l_2}
\end{equation}

\begin{rem}
  \label{re:sine-cosine}
  It is interesting to note that the both trigonometric functions are
  employed in quantum mechanics: sine is in the heart of the Moyal
  bracket%
  \index{Moyal bracket}%
  \index{bracket!Moyal}~\eqref{eq:repres-commutator} and cosine is responsible for
  the addition of probabilities~\eqref{eq:addition-functional-ell}. In
  the essence the commutator and probabilities took respectively the
  odd and even parts of the elliptic character \(e^{\rmi\myh s}\).
\end{rem}

\begin{example}
Take a  vector \(v_{(a,b)}\in\FSpace[\myh]{L}{2}(\Space{H}{n})\) defined by a
Gaussian\index{Gaussian} with mean value \((a,b)\) in the phase space%
\index{phase!space}%
\index{space!phase} for a harmonic oscillator of the
mass \(m\) and the frequency \(k\):
\begin{equation}
  \label{eq:gauss-state}
  v_{(a,b)}(q,p)=\exp\left(-\frac{2\pi k m}{\myhbar}(q-a)^2-\frac{2\pi}{\myhbar
      k m}(p-b)^2\right).
\end{equation}
A direct calculation shows:
\begin{eqnarray*}
  \lefteqn{\scalar{v_{(a,b)}}{\uir{}{\myhbar}(s,x,y)v_{(a',b')}}=\frac{4}{\myhbar}
  \exp\left(
    \pi \rmi \left(2s\myhbar+x (a+a')+y (b+b')\right)\frac{}{}\right.}\\
  &&\left.{} -\frac{\pi}{2 \myhbar k m }((\myhbar x+b-b')^2
    +(b-b')^2)
-\frac{\pi k m}{2\myhbar} ((\myhbar y+a'-a)^2
  + (a'-a)^2)
  \right)\\
  &=&\frac{4}{\myhbar}
  \exp\left(
    \pi \rmi \left(2s\myhbar+x (a+a')+y (b+b')\right)\frac{}{}\right.\\
  &&\left.{}  -\frac{\pi}{\myhbar k m }((b-b'+{\textstyle\frac{\myhbar x}{2}})^2
    +({\textstyle\frac{\myhbar x}{2}})^2)
    -\frac{\pi k m}{\myhbar} ((a-a'-{\textstyle\frac{\myhbar y}{2}})^2
    + ({\textstyle\frac{\myhbar y}{2}})^2) 
  \right)
\end{eqnarray*}
Thus the kernel
\(l_{(a,b)}=\scalar{v_{(a,b)}}{\uir{}{\myhbar}(s,x,y)v_{(a,b)}}\)~\eqref{eq:kernel-state}
for a state \(v_{(a,b)}\) is:
\begin{eqnarray}
  l_{(a,b)}&=&\frac{4}{\myhbar}
  \exp\left(
    2\pi \rmi (s\myhbar+xa+yb)\frac{}{}
    -\frac{\pi\myhbar}{2 k m }x^2
    -\frac{\pi k m \myhbar}{2\myhbar} y^2
  \right)
  \label{eq:single-slit}
\end{eqnarray}
An observable registering a particle at a point \(q=c\) of the
configuration space%
\index{configuration!space}%
\index{space!configuration} is \(\delta(q-c)\). On the Heisenberg group this
observable is given by the kernel:
\begin{equation}
  \label{eq:coordinate}
  X_c(s,x,y)=e^{2\pi\rmi (s\myhbar+x  c)}\delta(y).
\end{equation}
The measurement of \(X_c\) on the
state~\eqref{eq:gauss-state} (through the
kernel~\eqref{eq:single-slit}) predictably is: 
\begin{displaymath}
  \scalar{X_c}{l_{(a,b)}}=\sqrt{\frac{2k
  m}{\myhbar}}\exp\left(-\frac{2\pi k m}{\myhbar}(c-a)^2\right).
\end{displaymath}
\end{example}
\begin{example}
  Now take two states \(v_{(0,b)}\) and \(v_{(0,-b)}\), where for the
  simplicity we assume the mean values of coordinates vanish in the
  both cases.  Then the corresponding kernel~\eqref{eq:kernel-add} has
  the interference\index{interference} terms:
\begin{eqnarray*}
  l_i&=&  \scalar{v_{(0,b)}}{\uir{}{\myhbar}(s,x,y)v_{(0,-b)}}\\
  &=&\frac{4}{\myhbar}
  \exp\left(2\pi \rmi s\myhbar
    -\frac{\pi}{2 \myhbar k m }((\myhbar x+2b)^2
    +4b^2)
    -\frac{\pi\myhbar k m}{2} y^2
  \right).
\end{eqnarray*}
The measurement of \(X_c\)~\eqref{eq:coordinate} on this term contains
the oscillating part:
\begin{displaymath}
  \scalar{X_c}{l_i}=\sqrt{\frac{2k m}{\myhbar}} \exp\left(-\frac{2\pi k
    m }{\myhbar} c^2
  -\frac{2\pi}{ k
    m \myhbar}b^2+\frac{4\pi\rmi}{\myhbar} cb\right)
\end{displaymath}
Therefore on the kernel \(l\) corresponding to
the state \(v_{(0,b)}+v_{(0,-b)}\) the measurement is 
\begin{eqnarray*}
  \scalar{X_c}{l}
  &=&2\sqrt{\frac{2 k
      m}{\myhbar}}\exp\left(-\frac{2\pi k m}{\myhbar}c^2\right)
  \left(1+\exp\left(    -\frac{2\pi}{ k
      m \myhbar}b^2\right)\cos\left(\frac{4\pi}{\myhbar} cb\right)\right).
\end{eqnarray*}
\begin{figure}[htbp]
  \centering
  (a)\includegraphics[scale=.75]{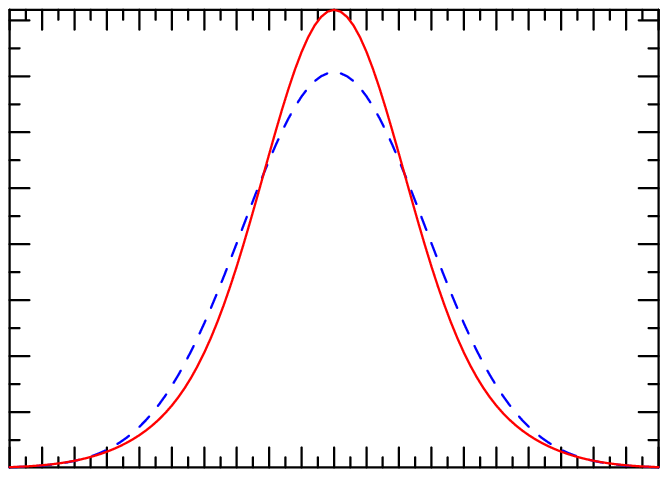}\hfill
  (b)\includegraphics[scale=.75]{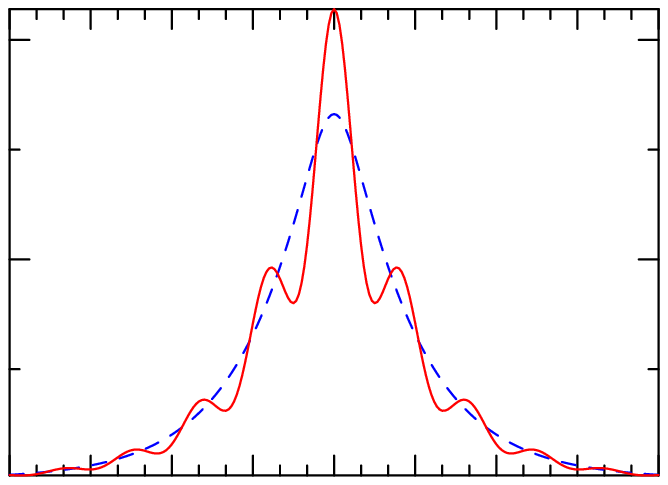}
  \caption{Quantum probabilities: the blue (dashed) graph shows the
    addition of probabilities without interaction, the red (solid)
    graph present the quantum interference. Left picture shows the
    Gaussian state~\eqref{eq:gauss-state}, the right---the rational
    state~\eqref{eq:poly-state}}
  \label{fig:quant-prob}
\end{figure}
The presence of the cosine term in the last expression can generate an
interference picture. In practise it does not happen for the minimal
uncertainty state~\eqref{eq:gauss-state} which we are using here: it
rapidly vanishes outside of the neighbourhood of zero, where
oscillations of the cosine occurs, see Fig.~\ref{fig:quant-prob}(a).
\end{example}

\begin{example}
To see a traditional interference pattern one can use a state which is far
from the minimal uncertainty. For example, we can consider the state:
\begin{equation}
  \label{eq:poly-state}
  u_{(a,b)}(q,p)=\frac{\myhbar^2}{((q-a)^2+\myhbar/ k
    m)((p-b)^2+\myhbar k m)}.
\end{equation}
To evaluate the observable \(X_c\)~\eqref{eq:coordinate} on the
state \(l(g)=\scalar{u_1}{\uir{}{h}(g)u_2}\)~\eqref{eq:kernel-state}
we use the following formula: 
\begin{displaymath}
  \scalar{X_c}{l}=\frac{2}{\myhbar}\int_{\Space{R}{n}} \hat{u}_1(q,
  2(q-c)/\myhbar)\, 
  \overline{\hat{u}_2(q, 2(q-c)/\myhbar)}\,dq,
\end{displaymath}
where \(\hat{u}_{i}(q,x)\) denotes  the partial Fourier transform
\(p\mapsto x\) of \(u_{i}(q,p)\). The formula is obtained by 
swapping order of integrations.  The numerical evaluation of the state
obtained by the addition \(u_{(0,b)}+u_{(0,-b)}\) is plotted on
Fig.~\ref{fig:quant-prob}(b), the red curve shows the canonical
interference pattern. 
\end{example}

\subsection{Ladder Operators and Harmonic Oscillator}
\label{sec:ladder-operators}

Let \(\uir{}{}\) be a representation of the Schr\"odinger group%
\index{Schr\"odinger!group}%
\index{group!Schr\"odinger}
\(G=\Space{H}{1}\rtimes\widetilde{\mathrm{Sp}}(2)\)~\eqref{eq:schrodinger-group}
in a space \(V\).  Consider the derived representation of the Lie
algebra \(\algebra{g}\)~\amscite{Lang85}*{\S~VI.1} and denote
\(\tilde{X}=\uir{}{}(X)\) for \(X\in\algebra{g}\). To see the
structure of the representation \(\uir{}{}\) we can decompose the
space \(V\) into eigenspaces of the operator \(\tilde{X}\) for some
\(X\in \algebra{g}\). The canonical example is the Taylor series in
complex analysis.

We are going to consider three cases corresponding to three
non-isomorphic subgroups~(\ref{eq:k-subgroup}--\ref{eq:ap-subgroup}) of
\(\Sp[2]\) starting from the compact case. Let \(H=Z\) be a
generator of the compact subgroup \(K\).  Corresponding
symplectomorphisms~\eqref{eq:sympl-auto} of the phase space%
\index{phase!space}%
\index{space!phase} are given
by orthogonal rotations with matrices \(\begin{pmatrix} \cos t & \sin
  t\\ -sin t& \cos t
\end{pmatrix}\). The Shale--Weil
representation~\eqref{eq:shale-weil-der} coincides with the
Hamiltonian of the harmonic oscillator%
\index{harmonic!oscillator}%
\index{oscillator!harmonic} in Schr\"odinger representation.

Since \(\widetilde{\mathrm{Sp}}(2)\) is a two-fold cover the
corresponding eigenspaces of a compact group \(\tilde{Z} v_k=\rmi k
v_k\) are parametrised by a half-integer \(k\in\Space{Z}{}/2\).
Explicitly for a half-integer \(k\) eigenvectors are:
\begin{equation}
  \label{eq:hermit-poly}
  v_k(q)=H_{k+\frac{1}{2}}\left(\sqrt{\frac{2\pi}{\myhbar}}q\right) e^{-\frac{\pi}{\myhbar}q^2},
\end{equation}
where \(H_k\) is the \emph{Hermite polynomial}%
\index{Hermite!polynomial}%
\index{polynomial!Hermite}%
\index{function!Hermite|see{Hermite polynomial}}%
~\citelist{\amscite{Folland89}*{\S~1.7}  \amscite{ErdelyiMagnusII}*{8.2(9)}}.  

From the point of view of quantum mechanics as well as the representation
theory it is beneficial to introduce the
ladder operators \(\ladder{\pm}\)~\eqref{eq:raising-lowering}, known
also as \emph{creation/annihilation} in quantum
mechanics~\citelist{\amscite{Folland89}*{p.~49} \cite{BoyerMiller74a}}.%
\index{ladder operator|(}%
\index{operator!ladder|(}
There are two ways to search for ladder operators: in
(complexified) Lie algebras \(\algebra{h}_1\) and \(\algebra{sp}_2\).
The later coincides with our consideration in
Section~\ref{sec:correspondence} in the essence.  

\subsubsection{Ladder Operators from the Heisenberg Group}
\label{sec:heis-group-oper}

Assuming \(\ladder{+}=a\tilde{X}+b\tilde{Y}\) we obtain from the
relations~(\ref{eq:cross-comm}--\ref{eq:cross-comm1})
and~\eqref{eq:raising-lowering} the linear equations with unknown
\(a\) and \(b\):
\begin{displaymath}
  a=\lambda_+ b, \qquad -b=\lambda_+ a.
\end{displaymath}
The equations have a solution if and only if \(\lambda_+^2+1=0\), and
the raising/lowering operators are \(\ladder{\pm}=
\tilde{X}\mp\rmi\tilde{Y}\). 
\begin{rem}
  Here we have an interesting asymmetric response: due to the
  structure of the semidirect product
  \(\Space{H}{1}\rtimes\widetilde{\mathrm{Sp}}(2)\) it is the
  symplectic group
  \index{symplectic!group}%
  \index{group!symplectic}%
  which acts on \(\Space{H}{1}\), not vise versa.
  However the Heisenberg group has a weak action in the opposite
  direction: it shifts eigenfunctions of \(\Sp[2]\).
\end{rem}

In the Schr\"odinger representation%
\index{Schr\"odinger!representation}%
\index{representation!Schr\"odinger}~\eqref{eq:schroedinger-rep-conf-der}
the ladder operators are
\begin{equation}
  \label{eq:ell-ladder-heisen-rep}
  \uir{}{\myhbar}(\ladder{\pm})= 2\pi\rmi q\pm\rmi\myhbar \frac{d}{dq}.
\end{equation}
The standard treatment of the harmonic oscillator in quantum mechanics,
which can be found in many textbooks, e.g.~\citelist{
  \amscite{Folland89}*{\S~1.7} \amscite{Gazeau09a}*{\S~2.2.3}}, 
is as follows. The vector  \(v_{-1/2}(q)=e^{-\pi q^2/\myhbar}\) is an
eigenvector of \(\tilde{Z}\) with the eigenvalue
\(-\frac{\rmi}{2}\). In addition \(v_{-1/2}\) is annihilated by
\(\ladder{+}\). Thus the chain~\eqref{eq:ladder-chain-1D} terminates to
the right and the complete set of eigenvectors of the harmonic
oscillator Hamiltonian is presented by \((\ladder{-})^k v_{-1/2}\)
with \(k=0, 1, 2, \ldots\).

We can make a wavelet transform%
\index{wavelet!transform} generated by the Heisenberg group with
the mother wavelet \(v_{-1/2}\), and the image will be the
Fock--Segal--Bargmann (FSB) space%
\index{Fock--Segal--Bargmann!space}%
\index{space!Fock--Segal--Bargmann}%
\index{FSB!space|see{Fock--Segal--Bargmann space}}%
\index{space!FSB|see{Fock--Segal--Bargmann space}} \citelist{\cite{Howe80b}
  \amscite{Folland89}*{\S~1.6}}. Since \(v_{-1/2}\) is the null solution
of \(\ladder{+}=\tilde{X}-\rmi \tilde{Y}\), then by
Cor.~\ref{co:cauchy-riemann} the image of
the wavelet transform will be null-solutions of the corresponding linear
combination of the Lie derivatives~\eqref{eq:h-lie-algebra}:
\begin{equation}
  \label{eq:CR-Bargmann}
  D=\overline{X^{r} -\rmi  Y^{r}}=(\partial_{ x} +\rmi\partial_{y})-\pi\myhbar(x-\rmi
y),
\end{equation}
which turns out to be the Cauchy--Riemann equation%
\index{Cauchy-Riemann operator}%
\index{operator!Cauchy-Riemann} on a weighted FSB-type space.

\subsubsection{Symplectic Ladder Operators}
\label{sec:sympl-ladd-oper}
We can also look for ladder operators within the Lie algebra
\(\algebra{sp}_2\), see~\S~\ref{sec:ellipt-ladd-oper} and~\amscite{Kisil09c}*{\S~8}.
Assuming \(\ladder[2]{+}=a\tilde{A}+b\tilde{B}+c\tilde{Z}\) from the
relations~\eqref{eq:sl2-commutator} and defining
condition~\eqref{eq:raising-lowering} we obtain the linear equations
with unknown \(a\), \(b\) and \(c\): 
\begin{displaymath}
  c=0, \qquad 2a=\lambda_+ b, \qquad -2b=\lambda_+ a.
\end{displaymath}
The equations have a solution if and only if \(\lambda_+^2+4=0\), and
the raising/lowering operators are \(\ladder[2]{\pm}=\pm\rmi
\tilde{A}+\tilde{B}\). In the Shale--Weil
representation~\eqref{eq:shale-weil-der} they turn out to be:
\begin{equation}
  \label{eq:ell-ladder-symplect}
  \ladder[2]{\pm}=\pm\rmi\left(\frac{q}{2}\frac{d}{dq}+\frac{1}{4}\right)-\frac{\myhbar\rmi}{8\pi}\frac{d^2}{dq^2}-\frac{\pi\rmi q^2}{2\myhbar}=-\frac{\rmi}{8\pi\myhbar}\left(\mp2\pi q+\myhbar\frac{d}{dq}\right)^2.
\end{equation}
Since this time \(\lambda_+=2\rmi\) the ladder operators
\(\ladder[2]{\pm}\) produce a shift on the
diagram~\eqref{eq:ladder-chain-1D} twice bigger than the operators
\(\ladder{\pm}\) from the Heisenberg group. After all, this is not
surprising since from the explicit
representations~\eqref{eq:ell-ladder-heisen-rep} and~\eqref{eq:ell-ladder-symplect} we get:
\begin{displaymath}
  \ladder[2]{\pm}=-\frac{\rmi}{8\pi\myhbar}(\ladder{\pm})^2.
\end{displaymath}
\index{ladder operator|)}%
\index{operator!ladder|)}%

\subsection{Hyperbolic Quantum Mechanics}
\label{sec:hyperb-repr-addt}

Now we turn to double numbers%
\index{number!double|indef}%
\index{double!number|indef} also known as hyperbolic, split-complex,
etc. numbers~\citelist{\cite{Yaglom79}*{App.~C} \cite{Ulrych05a}
  \cite{KhrennikovSegre07a}}. They form a two dimensional algebra \(\Space{O}{}\)
spanned by \(1\) and \(\rmh\) with the property \(\rmh^2=1\).  There
are zero divisors:
\begin{displaymath}
  \rmh_\pm=\textstyle\frac{1}{\sqrt{2}}(1\pm j), \qquad\text{ such that }\quad
  \rmh_+ \rmh_-=0 
  \quad
  \text{ and }
  \quad
  \rmh_\pm^2=\rmh_\pm.
\end{displaymath}
Thus double numbers algebraically isomorphic to two copies of
\(\Space{R}{}\) spanned by \(\rmh_\pm\). Being algebraically dull
double numbers are nevertheless interesting as a homogeneous
space~\cites{Kisil05a,Kisil09c} and they are relevant in
physics~\cites{Khrennikov05a,Ulrych05a,Ulrych08a}.  The combination of
p-mechanical approach with hyperbolic quantum mechanics was already
discussed in~\cite{BrodlieKisil03a}*{\S~6}.

For the hyperbolic character \(\chi_{\rmh \myh}(s)=e^{\rmh \myh
  s}=\cosh \myh s +\rmh\sinh \myh s\)
of \(\Space{R}{}\) one can define
the hyperbolic Fourier-type transform:
\begin{displaymath}
  \hat{k}(q)=\int_{\Space{R}{}} k(x)\,e^{-\rmh q x}dx.
\end{displaymath}
It can be understood in the sense of distributions on the space dual
to the set of analytic functions~\cite{Khrennikov08a}*{\S~3}. Hyperbolic
Fourier transform intertwines the derivative \(\frac{d}{dx}\) and
multiplication by \(\rmh q\)~\cite{Khrennikov08a}*{Prop.~1}.
\begin{example}
  For the Gaussian\index{Gaussian} the hyperbolic Fourier transform is the ordinary
  function (note the sign  difference!):
  \begin{displaymath}
    \int_{\Space{R}{}} e^{-x^2/2} e^{-\rmh q x}dx= \sqrt{2\pi}\, e^{q^2/2}.
  \end{displaymath}
  However the opposite identity:
  \begin{displaymath}
    \int_{\Space{R}{}} e^{x^2/2} e^{-\rmh q x}dx= \sqrt{2\pi}\, e^{-q^2/2}
  \end{displaymath}
  is true only in a suitable distributional sense. To this end we may
  note that \(e^{x^2/2}\) and \(e^{-q^2/2}\) are null solutions to the
  differential operators \(\frac{d}{dx}-x\) and \(\frac{d}{dq}+q\)
  respectively, which are intertwined (up to the factor \(\rmh\)) by
  the hyperbolic Fourier transform. The above differential operators
  \(\frac{d}{dx}-x\) and \(\frac{d}{dq}+q\) are images of the ladder
  operators%
  \index{ladder operator}%
  \index{operator!ladder}~\eqref{eq:ell-ladder-heisen-rep} in the Lie algebra of the Heisenberg group.
  They are intertwining by the Fourier transform, since this is an
  automorphism of the Heisenberg group~\cite{Howe80a}.
\end{example}
An elegant theory of hyperbolic Fourier transform may be achieved by a
suitable adaptation of~\cite{Howe80a}, which uses representation
theory of the Heisenberg group.

\subsubsection{Hyperbolic Representations of the Heisenberg Group} 
\label{sec:segre-quatern-hyperb}

Consider the space
\(\FSpace[\rmh]{F}{\myh}(\Space{H}{n})\) of \(\Space{O}{}\)-valued
functions on \(\Space{H}{n}\) with the property:
\begin{equation}
  \label{eq:induced-prop-h}
  f(s+s',h,y)=e^{\rmh \myh s'} f(s,x,y), \qquad \text{ for all }
  (s,x,y)\in \Space{H}{n},\ s'\in \Space{R}{} ,
\end{equation}
and the square integrability condition~\eqref{eq:L2-condition}. Then
the hyperbolic representation is obtained by the restriction of the
left shifts to \(\FSpace[\rmh]{F}{\myh}(\Space{H}{n})\).
To obtain an equivalent representation on the phase space%
\index{phase!space}%
\index{space!phase} we take
\(\Space{O}{}\)-valued functional of the Lie algebra \(\algebra{h}_n\): 
\begin{equation}
  \label{eq:hyp-character}
  \chi^j_{(\myh,q,p)}(s,x,y)=e^{\rmh(\myh s +qx+ py)}
  =\cosh (\myh s +qx+ py) + \rmh\sinh(\myh s +qx+ py).
\end{equation}
The hyperbolic Fock--Segal--Bargmann type representation%
  \index{Fock--Segal--Bargmann!representation!hyperbolic}%
  \index{representation!Fock--Segal--Bargmann!hyperbolic}%
  \index{hyperbolic!Fock--Segal--Bargmann representation}%
  \index{Heisenberg!group!Fock--Segal--Bargmann representation!hyperbolic}%
  \index{group!Heisenberg!Fock--Segal--Bargmann representation!hyperbolic} is intertwined
with the left group action by means of the Fourier
transform~\eqref{eq:fourier-transform} with the hyperbolic
functional~\eqref{eq:hyp-character}. Explicitly this representation is:
\begin{equation}
  \label{eq:segal-bargmann-hyp}
  \uir{}{\myhbar}(s,x,y): f (q,p) \mapsto 
  \textstyle e^{-\rmh(\myh s+qx+py)}
  f \left(q-\frac{\myh}{2} y, p+\frac{\myh}{2} x\right).
\end{equation}
For a hyperbolic Schr\"odinger type representation%
\index{representation!Heisenberg group!hyperbolic}%
\index{Schr\"odinger!representation!hyperbolic}%
\index{hyperbolic!Schr\"odinger representation!}%
\index{representation!Schr\"odinger!hyperbolic} we again use the
scheme described in \S~\ref{sec:concl-induc-repr}. Similarly to the
elliptic case one obtains the formula,
resembling~\eqref{eq:schroedinger-rep}:
\begin{equation}
    \label{eq:schroedinger-rep-hyp}
    [\uir{\rmh}{\chi}(s',x',y') f](x)=e^{-\rmh\myh (s'+xy'-x'y'/2)}f(x-x').
\end{equation}
Application of the hyperbolic Fourier transform produces a
Schr\"odinger type representation on the configuration space%
\index{configuration!space}%
\index{space!configuration},
cf.~\eqref{eq:schroedinger-rep-conf}: 
\begin{displaymath}
  [\uir{\rmh}{\chi}(s',x',y') \hat{f}\,](q)=e^{-\rmh\myh (s'+x'y'/2)
    -\rmh x' q}\,\hat{f}(q+\myh y').  
\end{displaymath}
The extension of this representation to kernels according
to~\eqref{eq:rho-extended-to-L1} generates hyperbolic
pseudodifferential operators introduced
in~\cite{Khrennikov08a}*{(3.4)}.

\subsubsection{Hyperbolic Dynamics}
\label{sec:hyperbolic-dynamics}

Similarly to the elliptic (quantum) case we consider a convolution
of two kernels on \(\Space{H}{n}\) restricted to
\(\FSpace[\rmh]{F}{\myh}(\Space{H}{n})\). The composition law becomes,
cf.~\eqref{eq:composition-ell}:
\begin{equation}
  \label{eq:composition-par}
  (k'*k)\hat{_s}
  =
  \int_{\Space{R}{2n}} e^{ {\rmh \myh}{}(xy'-yx')}\, \hat{k}'_s(\myh ,x',y')\,
 \hat{k}_s(\myh ,x-x',y-y')\,dx'dy'. 
\end{equation}
This is close to the calculus of hyperbolic PDO obtained
in~\cite{Khrennikov08a}*{Thm.~2}.
Respectively for the commutator of two convolutions we get,
cf.~\eqref{eq:repres-commutator}:
\begin{equation}
  \label{eq:commut-par}
  [k',k]\hat{_s}
  = 
  \int_{\Space{R}{2n}}\!\! \sinh(\myh
   (xy'-yx'))\, \hat{k}'_s(\myh ,x',y')\,
 \hat{k}_s(\myh ,x-x',y-y')\,dx'dy'. 
\end{equation}
This the hyperbolic version of the Moyal bracket%
\index{Moyal bracket!hyperbolic}%
\index{bracket!Moyal!hyperbolic}%
\index{hyperbolic!Moyal bracket},
cf.~\cite{Khrennikov08a}*{p.~849}, which generates the corresponding
image of the dynamic equation~\eqref{eq:universal}.
\begin{example}
  \begin{enumerate}
  \item 
    \label{it:hyp-harm-oscil} 
    For a quadratic Hamiltonian, e.g.  harmonic oscillator%
    \index{harmonic!oscillator}%
    \index{oscillator!harmonic} from
    Example~\ref{ex:p-harmonic}, the hyperbolic equation and
    respective dynamics is identical to quantum considered before.
  \item Since \(\frac{\partial}{\partial s}\) acts on
    \(\FSpace[\rmh]{F}{2}(\Space{H}{n})\) as multiplication by \(\rmh
    \myh\) and \(\rmh^2=1\), the hyperbolic image of the unharmonic
    equation%
    \index{unharmonic!oscillator!hyperbolic}%
    \index{hyperbolic!unharmonic oscillator}%
    \index{oscillator!unharmonic!hyperbolic}~\eqref{eq:p-unharm-osc-dyn} becomes:
    \begin{displaymath}
      \dot{f}=     \left(m  k^2 q\frac{\partial}{\partial
          p}+\frac{\lambda}{6}\left(3q^2\frac{\partial}{\partial p}
          +\frac{\myhbar^2}{4}\frac{\partial^3}{\partial p^3}\right)-\frac{1}{m}
        p \frac{\partial}{\partial q} \right) f. 
    \end{displaymath}
    The difference with quantum mechanical
    equation~\eqref{eq:q-unhar-dyn} is in the sign of the cubic
    derivative. 
  \end{enumerate}
\end{example}

\subsubsection{Hyperbolic Probabilities}
\label{sec:hyperb-prob}
\begin{figure}[htbp]
  \centering
  (a)\includegraphics[scale=.75]{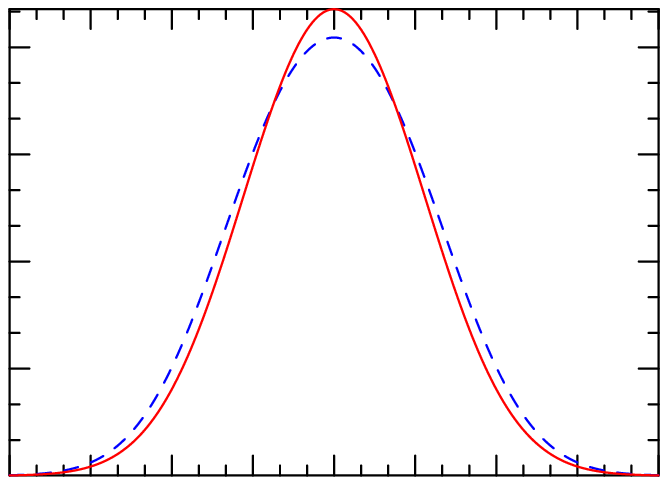}\hfill
  (b)\includegraphics[scale=.75]{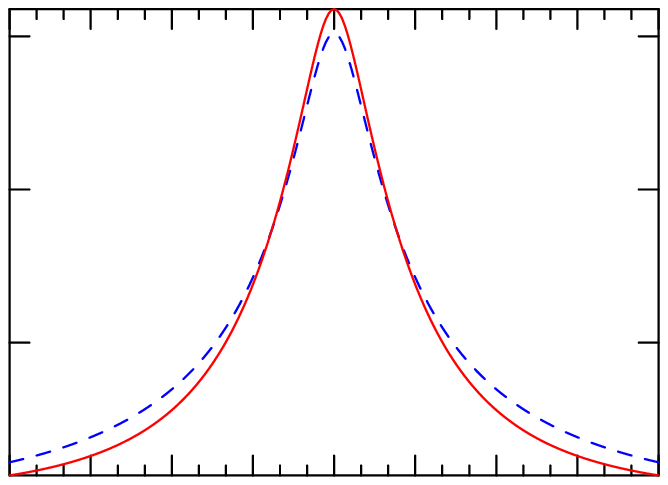}
  \caption{Hyperbolic probabilities: the blue (dashed) graph shows the
    addition of probabilities without interaction, the red (solid)
    graph present the quantum interference. Left picture shows the
    Gaussian state~\eqref{eq:gauss-state}, with the same distribution
    as in quantum mechanics, cf.~Fig.~\ref{fig:quant-prob}(a). The
    right picture shows the rational
    state~\eqref{eq:poly-state}, note the absence of interference
    oscillations in comparison with the quantum
    state on~Fig.~\ref{fig:quant-prob}(b).} 
  \label{fig:hyp-prob}
\end{figure}

To calculate probability%
\index{probability!hyperbolic}%
\index{hyperbolic!probability} distribution generated by a hyperbolic state we are
using the general procedure from
Section~\ref{sec:states-probability}. The main differences with the
quantum case are as follows:
\begin{enumerate}
\item The real number \(A\) in the
  expression~\eqref{eq:addition-functional} for the addition of
  probabilities is bigger than \(1\) in absolute value. Thus it can
  be associated with the hyperbolic cosine \(\cosh \alpha \),
  cf.~Rem.~\ref{re:sine-cosine}, for certain phase
  \(\alpha\in\Space{R}{}\)~\cite{Khrennikov08a}.
\item The nature of hyperbolic interference on two slits is affected
  by the fact that \(e^{\rmh \myh s}\) is not periodic and the
  hyperbolic exponent \(e^{\rmh t}\) and cosine \(\cosh t\) do not
  oscillate. It is worth to notice that for Gaussian\index{Gaussian} states the
  hyperbolic interference is exactly the same as quantum one,
  cf.~Figs.~\ref{fig:quant-prob}(a) and~\ref{fig:hyp-prob}(a). This is
  similar to coincidence of quantum and hyperbolic dynamics of
  harmonic oscillator.

  The contrast between two types of interference is prominent for
  the rational state~\eqref{eq:poly-state}, which is far from the
  minimal uncertainty, see the different patterns
  on Figs.~\ref{fig:quant-prob}(b) and~\ref{fig:hyp-prob}(b).
\end{enumerate}

\subsubsection{Ladder Operators for the Hyperbolic Subgroup}
\label{sec:hiperbolic-subgroup}%
\index{ladder operator|(}%
\index{operator!ladder|(}

Consider the case of the Hamiltonian \(H=2B\), which is a repulsive
(hyperbolic) harmonic oscillator%
\index{harmonic!oscillator!repulsive (hyperbolic)}%
\index{repulsive!harmonic oscillator}%
\index{hyperbolic!harmonic oscillator}%
\index{oscillator!harmonic!repulsive
  (hyperbolic)}~\amscite{Wulfman10a}*{\S~3.8}. The corresponding
one-dimensional subgroup of symplectomorphisms produces hyperbolic
rotations of the phase space, see Fig.~\ref{fig:rotations}.%
\index{symplectic!transformation}%
\index{transformation!symplectic} The
eigenvectors \(v_\mu\) of the operator
\begin{displaymath}
  \uir{\text{SW}}{\myhbar}(2B)v_\nu
  =-\rmi\left(\frac{\myhbar}{4\pi}\frac{d^2}{dq^2}+\frac{\pi q^2}{\myhbar}\right)v_\nu
  =\rmi\nu v_\nu, 
\end{displaymath}
are \emph{Weber--Hermite}%
\index{Weber--Hermite function}%
\index{function!Weber--Hermite} (or \emph{parabolic cylinder}%
\index{parabolic!cylinder function|see{Weber--Hermite function}}%
\index{function!parabolic cylinder|see{Weber--Hermite function}}) functions
\(v_{\nu}=D_{\nu-\frac{1}{2}}\left(\pm2e^{\rmi \frac{\pi}{4}}\sqrt{\frac{\pi}{\myhbar}} q\right)\),
see~\citelist{\amscite{ErdelyiMagnusII}*{\S~8.2}
  \cite{SrivastavaTuanYakubovich00a}} for fundamentals of
Weber--Hermite functions and~\cite{ATorre08a} for further
illustrations and applications in optics.

The corresponding one-parameter group is not compact and the
eigenvalues of the operator \(2\tilde{B}\) are not restricted by any
integrality condition, but the raising/lowering operators are
still important~\citelist{\amscite{HoweTan92}*{\S~II.1}
  \amscite{Mazorchuk09a}*{\S~1.1}}. We again seek solutions in two
subalgebras \(\algebra{h}_1\) and \(\algebra{sp}_2\) separately.
However the additional options will be provided by a choice of the
number system: either complex or double.

\begin{example}[Complex Ladder Operators]
\label{sec:compl-ladd-oper}
  
Assuming
\(\ladder[h]{+}=a\tilde{X}+b\tilde{Y}\) from the
commutators~(\ref{eq:cross-comm}--\ref{eq:cross-comm1}) we obtain
the linear equations:
\begin{equation}
  \label{eq:hyp-ladder-compatib}
  -a=\lambda_+ b, \qquad -b=\lambda_+ a.
\end{equation}
The equations have a solution if and only if \(\lambda_+^2-1=0\).
Taking the real roots \(\lambda=\pm1\) we obtain that the
raising/lowering operators are
\(\ladder[h]{\pm}=\tilde{X}\mp\tilde{Y}\).  In the Schr\"odinger
representation%
\index{Schr\"odinger!representation}%
\index{representation!Schr\"odinger}~\eqref{eq:schroedinger-rep-conf-der}
the ladder operators are
\begin{equation}
  \label{eq:ell-ladder-heisen-rep1}
  \ladder[h]{\pm}= 2\pi\rmi q\pm \myhbar \frac{d}{dq}.
\end{equation}
The null solutions \(v_{\pm\frac{1}{2}}(q)=e^{\pm\frac{\pi\rmi}{\myhbar}
  q^2}\) to operators \(\uir{}{\myhbar}(\ladder{\pm})\) are also
eigenvectors of the Hamiltonian \(\uir{\text{SW}}{\myhbar}(2B)\) with the
eigenvalue \(\pm\frac{1}{2}\).  However the important distinction from
the elliptic case is, that they are not square-integrable on the real line
anymore.

We can also look for ladder operators within the \(\algebra{sp}_2\),
that is in the form \(\ladder[2h]{+}=a\tilde{A}+b\tilde{B}+c\tilde{Z}\)
for the commutator \([2\tilde{B},\ladder[h]{+}]=\lambda
\ladder[h]{+}\), see \S~\ref{sec:hyperb-ladd-oper}. Within complex
numbers we get only the values \(\lambda=\pm 2\) with the ladder
operators \(\ladder[2h]{\pm}=\pm2\tilde{A}+\tilde{Z}/2\),
see~\citelist{\amscite{HoweTan92}*{\S~II.1}
  \amscite{Mazorchuk09a}*{\S~1.1}}. Each indecomposable \(\algebra{h}_1\)-
or \(\algebra{sp}_2\)-module is formed by a one-dimensional chain of
eigenvalues with a transitive action of ladder operators \(\ladder[h]{\pm}\)
or \(\ladder[2h]{\pm}\) respectively. And we again have a
quadratic relation between the ladder operators:
\begin{displaymath}
  \ladder[2h]{\pm}=\frac{\rmi}{4\pi\myhbar}(\ladder[h]{\pm})^2.
\end{displaymath}
\end{example}

\subsubsection{Double Ladder Operators}
\label{sec:double-ladd-oper}
  
There are extra possibilities in in the context of hyperbolic quantum
mechanics~\citelist{\cite{Khrennikov03a} \cite{Khrennikov05a}
  \cite{Khrennikov08a}}.  Here we use the representation of
\(\Space{H}{1}\) induced by a hyperbolic character \(e^{\rmh \myh
  t}=\cosh (\myh t)+\rmh\sinh(\myh t)\), see
\amscite{Kisil10a}*{(4.5)}, and obtain
the hyperbolic representation of \(\Space{H}{1}\),
cf.~\eqref{eq:schroedinger-rep-conf}:  
\begin{equation}
  \label{eq:schroedinger-rep-conf-hyp}
    [\uir{\rmh}{\myh}(s',x',y') \hat{f}\,](q)=e^{\rmh\myh (s'-x'y'/2)
    +\rmh x' q}\,\hat{f}(q-\myh y').  
\end{equation}
The corresponding derived representation is
\begin{equation}
  \label{eq:schroedinger-rep-conf-der-hyp}
  \uir{\rmh}{\myh}(X)=\rmh q,\qquad \uir{\rmh}{\myh}(Y)=-\myh \frac{d}{dq},
  \qquad
  \uir{\rmh}{\myh}(S)=\rmh\myh I.
\end{equation}
Then the associated Shale--Weil derived  representation of \(\algebra {sp}_2\) in
the Schwartz space \(\FSpace{S}{}(\Space{R}{})\) is, cf.~\eqref{eq:shale-weil-der}:
\begin{equation}
  \label{eq:shale-weil-der-double}
  \uir{\text{SW}}{\myh}(A) =-\frac{q}{2}\frac{d}{dq}-\frac{1}{4},\quad
  \uir{\text{SW}}{\myh}(B)=\frac{\rmh\myh}{4}\frac{d^2}{dq^2}-\frac{\rmh q^2}{4\myh},\quad
  \uir{\text{SW}}{\myh}(Z)=-\frac{\rmh\myh}{2}\frac{d^2}{dq^2}-\frac{\rmh q^2}{2\myh}.
\end{equation}
Note that \(\uir{\text{SW}}{\myh}(B)\) now generates a usual harmonic
oscillator, not the repulsive one like 
\(\uir{\text{SW}}{\myhbar}(B)\) in \eqref{eq:shale-weil-der}. 
However the expressions in the quadratic algebra are still the same (up to a factor),
cf.~(\ref{eq:quadratic-A}--\ref{eq:quadratic-Z}):
\begin{eqnarray}
  \label{eq:quadratic-A-hyp}
  \qquad\uir{\text{SW}}{\myh}(A) &=&
  -\frac{\rmh}{2\myh}(\uir{\rmh}{\myh}(X)\uir{\rmh}{\myh}(Y)
  -{\textstyle\frac{1}{2}}\uir{\rmh}{\myh}(S))\\
  &=&-\frac{\rmh}{4\myh}(\uir{\rmh}{\myh}(X)\uir{\rmh}{\myh}(Y)
  +\uir{\rmh}{\myh}(Y)\uir{\rmh}{\myh}(X)),\nonumber \\ 
  \label{eq:quadratic-B-hyp}
  \uir{\text{SW}}{\myh}(B) &=&
  \frac{\rmh}{4\myh}(\uir{\rmh}{\myh}(X)^2-\uir{\rmh}{\myh}(Y)^2), \\
  \label{eq:quadratic-Z-hyp}
  \uir{\text{SW}}{\myh}(Z)
  &=&-\frac{\rmh}{2\myh}(\uir{\rmh}{\myh}(X)^2+\uir{\rmh}{\myh}(Y)^2). 
\end{eqnarray}
This is due to the Principle~\ref{pr:simil-corr-principle} of
similarity and correspondence%
\index{principle!similarity and correspondence}: we can swap operators \(Z\) and \(B\) with
simultaneous replacement of hypercomplex units \(\rmi\) and \(\rmh\).

The eigenspace of the operator \(2\uir{\text{SW}}{\myh}(B)\) with an
eigenvalue \(\rmh \nu\) are spanned by the Weber--Hermite
functions%
\index{Weber--Hermite function}%
\index{function!Weber--Hermite} \(D_{-\nu-\frac{1}{2}}\left(\pm\sqrt{\frac{2}{\myh}}x\right)\),
see~\amscite{ErdelyiMagnusII}*{\S~8.2}.  Functions \(D_\nu\) are
generalisations of the Hermit functions%
\index{Hermite!polynomial}%
\index{polynomial!Hermite}~\eqref{eq:hermit-poly}.

The compatibility condition for a ladder operator within the Lie algebra
\(\algebra{h}_1\) will be~\eqref{eq:hyp-ladder-compatib} as before,
since it depends only on the
commutators~(\ref{eq:cross-comm}--\ref{eq:cross-comm1}). Thus we still
have the set of ladder operators corresponding to values
\(\lambda=\pm1\):
\begin{displaymath}
  \ladder[h]{\pm}=\tilde{X}\mp\tilde{Y}=\rmh q\pm\myh \frac{d}{dq}.
\end{displaymath}
Admitting double numbers%
\index{number!double}%
\index{double!number} we have an extra way to satisfy
\(\lambda^2=1\) in~\eqref{eq:hyp-ladder-compatib} with values
\(\lambda=\pm\rmh\).  Then there is an additional pair of hyperbolic
ladder operators, which are identical (up to factors)
to~\eqref{eq:ell-ladder-heisen-rep}:
\begin{displaymath}
  \ladder[\rmh]{\pm}=\tilde{X}\mp\rmh\tilde{Y}=\rmh q\pm\rmh\myh \frac{d}{dq}.
\end{displaymath}
Pairs \(\ladder[h]{\pm}\) and \(\ladder[\rmh]{\pm}\) shift
eigenvectors in the ``orthogonal'' directions changing their
eigenvalues by \(\pm1\) and \(\pm\rmh\).  Therefore an indecomposable
\(\algebra{sp}_2\)-module can be para\-metrised by a two-dimensional
lattice of eigenvalues in double numbers, see
Fig.~\ref{fig:2D-lattice}.

The following functions 
\begin{eqnarray*}
  v_{\frac{1}{2}}^{\pm\myh}(q)&=&e^{\mp\rmh
    q^2/(2\myh)}=\cosh\frac{q^2}{2\myh}\mp \rmh\sinh \frac{q^2}{2\myh},\\
  v_{\frac{1}{2}}^{\pm\rmh}(q)&=&e^{\mp  q^2/(2\myh)}
\end{eqnarray*}
are null solutions to the operators \(\ladder[h]{\pm}\) and
\(\ladder[\rmh]{\pm}\) respectively. They are also eigenvectors of
\(2\uir{\text{SW}}{\myh}(B)\) with eigenvalues \(\mp\frac{\rmh}{2}\)
and \(\mp\frac{1}{2}\) respectively. If these functions are used as
mother wavelets for the wavelet transforms generated by the Heisenberg
group, then the image space will consist of the null-solutions of the
following differential operators, see Cor.~\ref{co:cauchy-riemann}:
\begin{displaymath}\textstyle
  D_{h}=\overline{X^{r} - Y^{r}}=(\partial_{ x} -\partial_{y})+\frac{\myh}{2}(x+y),
\qquad
  D_{\rmh}=\overline{X^{r} - \rmh Y^{r}}=(\partial_{ x} +\rmh\partial_{y})-\frac{\myh}{2}(x-\rmh
y),
\end{displaymath}
for \(v_{\frac{1}{2}}^{\pm\myh}\) and \(v_{\frac{1}{2}}^{\pm\rmh}\)
respectively. This is again in line with the classical
result~\eqref{eq:CR-Bargmann}. However annihilation of the eigenvector
by a ladder operator does not mean that the part of the 2D-lattice
becomes void since it can be reached via alternative routes on this
lattice. Instead of multiplication by a zero, as it happens in the
elliptic case, a half-plane of eigenvalues will be multiplied by the
divisors of zero%
\index{divisor!zero}%
\index{zero!divisor} \(1\pm\rmh\).

We can also search ladder operators within the algebra
\(\algebra{sp}_2\) and admitting double numbers we will again find two
sets of them, cf.~\S~\ref{sec:hyperb-ladd-oper}:
\begin{eqnarray*}
  \ladder[2h]{\pm} &=&\pm\tilde{A}+\tilde{Z}/2 =
   \mp\frac{q}{2}\frac{d}{dq}\mp\frac{1}{4}- \frac{\rmh\myh}{4}\frac{d^2}{dq^2}-\frac{\rmh q^2}{4\myh}=-\frac{\rmh}{4\myh}(\ladder[h]{\pm})^2,\\
  \ladder[2\rmh]{\pm}&=&\pm\rmh\tilde{A}+\tilde{Z}/2=  
  \mp\frac{\rmh q}{2}\frac{d}{dq}\mp\frac{\rmh}{4}-\frac{\rmh\myh}{4}\frac{d^2}{dq^2}-\frac{\rmh q^2}{4\myh}=-\frac{\rmh}{4\myh}(\ladder[\rmh]{\pm})^2.
\end{eqnarray*}
Again the operators \(\ladder[2h]{\pm}\) and \(\ladder[2h]{\pm}\) produce
double shifts in the orthogonal directions on the same two-dimensional
lattice in Fig.~\ref{fig:2D-lattice}.%
\index{ladder operator|)}%
\index{operator!ladder|)}%

\subsection{Parabolic (Classical) Representations on the Phase Space}
\label{sec:class-repr-phase}
After the previous two cases it is natural to link classical mechanics
with dual numbers%
\index{dual!number|indef}%
\index{number!dual|indef}%
\index{dual!number|(}%
\index{number!dual|(} generated by the parabolic unit \(\rmp^2=0\).
Connection of the parabolic unit \(\rmp\) with the Galilean group of
symmetries of classical mechanics is around for a
while~\cite{Yaglom79}*{App.~C}. 

However the nilpotency of the parabolic unit \(\rmp\) make
it difficult if we will work with dual number valued functions only.
To overcome this issue we consider a commutative real algebra
\(\algebra{C}\) spanned by \(1\), \(\rmi\), \(\rmp\) and \(\rmi\rmp\)
with identities \(\rmi^2=-1\) and \(\rmp^2=0\). A seminorm on
\(\algebra{C}\) is defined as follows:
\begin{displaymath}
  \modulus{a+b\rmi+c\rmp+d\rmi\rmp}^2=a^2+b^2.
\end{displaymath}

\subsubsection{Classical Non-Commutative Representations}
\label{sec:class-non-comm}
We wish to build a representation of the Heisenberg group which will
be a classical analog of the Fock--Segal--Barg\-mann
representation~\eqref{eq:stone-inf}.  To this end we introduce the
space \(\FSpace[\rmp]{F}{\myh}(\Space{H}{n})\) of
\(\algebra{C}\)-valued functions on \(\Space{H}{n}\) with the
property:
\begin{equation}
  \label{eq:induced-prop-p}
  f(s+s',h,y)=e^{\rmp \myh s'} f(s,x,y), \qquad \text{ for all }
  (s,x,y)\in \Space{H}{n},\ s'\in \Space{R}{} ,
\end{equation}
and the square integrability condition~\eqref{eq:L2-condition}. It is
invariant under the left shifts and we restrict the left group action to
\(\FSpace[\rmp]{F}{\myh}(\Space{H}{n})\). 

There is an unimodular
\(\algebra{C}\)-valued function on the Heisenberg group parametrised
by a point \((\myh, q, p)\in\Space{R}{2n+1}\):
\begin{displaymath}
  E_{(\myh,q,p)}(s,x,y)= e^{2\pi(\rmp s\myhbar+\rmi xq + \rmi yp)}=e^{2\pi\rmi (xq + yp)}(1+\rmp s\myh).
\end{displaymath}
This function, if used instead of the ordinary exponent, produces a modification
\(\oper{F}_c\) of the Fourier transform~\eqref{eq:fourier-transform}.
The transform intertwines the left regular representation with the following
action on \(\algebra{C}\)-valued functions on the phase space:%
\index{Fock--Segal--Bargmann!representation!classic (parabolic)}%
\index{representation!Fock--Segal--Bargmann!classic (parabolic)}%
\index{classic!Fock--Segal--Bargmann representation}%
\index{parabolic!Fock--Segal--Bargmann representation}%
\index{Heisenberg!group!Fock--Segal--Bargmann representation!classic (parabolic)}%
\index{group!Heisenberg!Fock--Segal--Bargmann representation!classic (parabolic)}
\begin{equation}
  \label{eq:dual-repres}
  \uir{\rmp}{\myh}(s,x,y): f(q,p) \mapsto e^{-2\pi\rmi(xq+yp)}(f(q,p)
  +\rmp\myh(s f(q,p) +\frac{y}{2\pi\rmi}f'_q(q,p)-\frac{x}{2\pi\rmi}f'_p(q,p))).
\end{equation}
\begin{rem}
  \label{re:classic-rep}
  Comparing the traditional
  infinite-dimensional~\eqref{eq:stone-inf} and
  one-dimen\-sional~\eqref{eq:commut-repres} representations of
  \(\Space{H}{n}\) we can note that the properties of the
  representation~\eqref{eq:dual-repres} are a non-trivial mixture of
  the former:  
  \begin{enumerate}
  \item \label{it:class-non-commut} 
    The action~\eqref{eq:dual-repres}
    is non-commutative, similarly to the quantum
    representation~\eqref{eq:stone-inf} and unlike the classical
    one~\eqref{eq:commut-repres}. This non-commutativity will produce
    the Hamilton equations below in a way very similar to Heisenberg
    equation, see Rem.~\ref{re:hamilton-from-nc}.
  \item \label{it:class-locality} The
    representation~\eqref{eq:dual-repres} does not change the support
    of a function \(f\) on the phase space, similarly to the
    classical representation~\eqref{eq:commut-repres} and unlike the
    quantum one~\eqref{eq:stone-inf}. Such a localised action will be
    responsible later for an absence of an interference in classical
    probabilities.
  \item The parabolic representation~\eqref{eq:dual-repres} can not be
    derived from either the elliptic~\eqref{eq:stone-inf} or
    hyperbolic~\eqref{eq:segal-bargmann-hyp} by the plain substitution
    \(\myh=0\).
  \end{enumerate}
\end{rem}
We may also write a classical Schr\"odinger type representation.%
\index{representation!Heisenberg group!classic (parabolic)}%
\index{Schr\"odinger!representation!classic (parabolic)}%
\index{parabolic!Schr\"odinger representation!}%
\index{classic!Schr\"odinger representation!}%
\index{representation!Schr\"odinger!classic (parabolic)}
According to \S~\ref{sec:concl-induc-repr} we get a representation formally
very similar to the elliptic~\eqref{eq:schroedinger-rep} and
hyperbolic versions~\eqref{eq:schroedinger-rep-hyp}:
\begin{eqnarray}
    \label{eq:schroedinger-rep-par}
    [\uir{\rmp}{\chi}(s',x',y') f](x)&=&e^{-\rmp\myh
      (s'+xy'-x'y'/2)}f(x-x')\\
    &=&(1-\rmp\myh (s'+xy'-\textstyle\frac{1}{2}x'y')) f(x-x').\nonumber 
\end{eqnarray}
However due to nilpotency of \(\rmp\) the (complex) Fourier transform
\(x\mapsto q\) produces a different formula for parabolic
Schr\"odinger type representation in the configuration space%
\index{configuration!space}%
\index{space!configuration},
cf.~\eqref{eq:schroedinger-rep-conf}
and~\eqref{eq:schroedinger-rep-conf-hyp}:
\begin{displaymath}
    [\uir{\rmp}{\chi}(s',x',y') \hat{f}](q)= e^{2\pi\rmi x' q}\left(
    \left(1-\rmp\myh (s'-{\textstyle\frac{1}{2}}x'y')\right)    \hat{f}(q)
    +\frac{\rmp\myh y'}{2\pi\rmi}\hat{f}'(q)\right).
\end{displaymath}
This representation shares all properties mentioned in
Rem.~\ref{re:classic-rep} as well.

\subsubsection{Hamilton Equation}
\label{sec:hamilton-equation}


The identity \(e^{ \rmp t}-e^{ -\rmp t}= 2\rmp t\) can be interpreted
as a parabolic version of the sine function, while the parabolic
cosine is identically equal to
one, cf. \S~\ref{sec:hyperc-char} and~\cites{HerranzOrtegaSantander99a,Kisil07a}.  From this we obtain
the parabolic version of the commutator~\eqref{eq:repres-commutator}:
\begin{eqnarray*}
  [k',k]\hat{_s}(\rmp \myh, x,y) 
  &=& 
  \rmp \myh\int_{\Space{R}{2n}}
 (xy'-yx') \\
 &&{}\times\, \hat{k}'_s(\rmp \myh,x',y')  \,
 \hat{k}_s(\rmp \myh,x-x',y-y')\,dx'dy', \nonumber 
\end{eqnarray*}
for the partial parabolic Fourier-type transform \(\hat{k}_s\) of the
kernels.  Thus the parabolic representation of the dynamical
equation~\eqref{eq:universal} becomes:
\begin{equation}
  \label{eq:dynamics-par}
  \rmp\myh \frac{d\hat{f}_s}{dt}(\rmp\myh,x,y;t)=
  \rmp \myh \int_{\Space{R}{2n}}
 (xy'-yx')\, 
\hat{H}_s(\rmp \myh,x',y')  \,
 \hat{f}_s(\rmp \myh,x-x',y-y';t)\,dx'dy', 
\end{equation}
Although there is no possibility to divide by \(\rmp\) (since it is a
zero divisor) we can obviously eliminate \(\rmp \myh \) from the both
sides if the rest of the expressions are real.  Moreover this can be
done ``in advance'' through a kind of the antiderivative operator
considered in~\cite{Kisil02e}*{(4.1)}. This will prevent ``imaginary
parts'' of the remaining expressions (which contain the factor
\(\rmp\)) from vanishing.
\begin{rem}
  It is noteworthy that the Planck constant%
  \index{Planck!constant}%
  \index{constant!Planck} completely disappeared
  from the dynamical equation. Thus the only prediction about it
  following from our construction is \(\myh\neq 0\), which was
  confirmed by experiments, of course. 
\end{rem}
Using the duality between the Lie algebra of \(\Space{H}{n}\) and the
phase space we can find an adjoint equation for observables on the
phase space. To this end we apply the usual Fourier transform
\((x,y)\mapsto(q,p)\). It turn to be the Hamilton
equation%
\index{Hamilton!equation}%
\index{equation!Hamilton}~\cite{Kisil02e}*{(4.7)}.  However the transition to the phase
space is more a custom rather than a necessity and in many cases we
can efficiently work on the Heisenberg group itself.

\begin{rem}
  \label{re:hamilton-from-nc}
  It is noteworthy, that the non-commutative
  representation~\eqref{eq:dual-repres} allows to obtain the Hamilton
  equation directly from the commutator
  \([\uir{\rmp}{\myh}(k_1),\uir{\rmp}{\myh}(k_2)]\). Indeed its
  straightforward  evaluation will produce exactly the above expression. On
  the contrast such a commutator for the commutative
  representation~\eqref{eq:commut-repres} is zero and to obtain the
  Hamilton equation we have to work with an additional tools, e.g. an
  anti-derivative~\cite{Kisil02e}*{(4.1)}. 
\end{rem}

\begin{example}
  \begin{enumerate}
  \item For the harmonic oscillator%
    \index{harmonic!oscillator}%
    \index{oscillator!harmonic} in Example~\ref{ex:p-harmonic} the
    equation~\eqref{eq:dynamics-par} again reduces to the
    form~\eqref{eq:p-harm-osc-dyn} with the solution given
    by~\eqref{eq:p-harm-sol}. The adjoint equation of the harmonic
    oscillator on the phase space is not different from the quantum
    written in
    Example~\ref{ex:quntum-oscillators}(\ref{it:q-harmonic}). This is
    true for any Hamiltonian of at most quadratic order.
  \item 
    For non-quadratic Hamiltonians classical and quantum dynamics
    are different, of course. For example, 
    the cubic term of \(\partial_s\) in the
    equation~\eqref{eq:p-unharm-osc-dyn} will generate the factor
    \(\rmp^3=0\) and thus vanish. Thus the
    equation~\eqref{eq:dynamics-par} of the unharmonic oscillator%
    \index{unharmonic!oscillator}%
    \index{oscillator!unharmonic} on
    \(\Space{H}{n}\) becomes:
    \begin{displaymath}
      \dot{f}=     \left(m  k^2 y\frac{\partial}{\partial x}
        +\frac{\lambda y}{2}\frac{\partial^2}{\partial x^2} 
          -\frac{1}{m} x
        \frac{\partial}{\partial y} \right) f. 
   \end{displaymath}
   The adjoint equation on the phase space is:
    \begin{displaymath}
      \dot{f}=     \left(\left(m  k^2
          q+\frac{\lambda}{2}q^2\right)
        \frac{\partial}{\partial p} -\frac{1}{m}  p \frac{\partial}{\partial q} \right) f. 
    \end{displaymath}
    The last equation is the classical
    Hamilton equation generated by the
    cubic potential~\eqref{eq:unharmonic-hamiltonian}. Qualitative
    analysis of its dynamics can be found in many textbooks
    \citelist{\cite{Arnold91}*{\S~4.C, Pic.~12} \cite{PercivalRichards82}*{\S~4.4}}. 
  \end{enumerate}
\end{example}

\begin{rem}
  We have obtained the \emph{Poisson bracket}%
  \index{Poisson!bracket}%
  \index{bracket!Poisson} from the commutator of
  convolutions on \(\Space{H}{n}\) without any quasiclassical limit
  \(\myh\rightarrow 0\). This has a common source with the deduction
  of main calculus theorems in~\cite{CatoniCannataNichelatti04} based
  on dual numbers. As explained in~\cite{Kisil05a}*{Rem.~6.9} this is
  due to the similarity between the parabolic unit \(\rmp\) and the
  infinitesimal number used in non-standard analysis~\cite{Devis77}.
  In other words, we never need to take care about terms of order
  \(O(\myh^2)\) because they will be wiped out by \(\rmp^2=0\).
\end{rem}
An alternative derivation of classical dynamics from the Heisenberg
group is given in the recent paper~\cite{Low09a}.

\subsubsection{Classical probabilities}
\label{sec:class-prob}
It is worth to notice that dual numbers are not only helpful in
reproducing classical Hamiltonian dynamics, they also provide the
classic rule for addition of probabilities.%
\index{probability!classic (parabolic)}%
\index{classic!probability}%
\index{parabolic!probability|see{classic probability}}
We use the same formula~\eqref{eq:kernel-state} to calculate kernels of
the states. The important difference now that the
representation~\eqref{eq:dual-repres} does not change the support of
functions. Thus if we calculate the correlation term
\(\scalar{v_1}{\uir{}{}(g)v_2}\) in~\eqref{eq:kernel-add}, then it
will be zero for every two vectors \(v_1\) and \(v_2\) which have 
disjoint supports in the phase space. Thus no interference%
\index{interference} similar to
quantum or hyperbolic cases (Subsection~\ref{sec:quantum-probabilities})
is possible.

\subsubsection{Ladder Operator for the Nilpotent Subgroup}
\label{sec:nilpotent-subgroup}%
\index{ladder operator|(}%
\index{operator!ladder|(}

Finally we look for ladder operators for the Hamiltonian
\(\tilde{B}+\tilde{Z}/2\) or, equivalently,
\(-\tilde{B}+\tilde{Z}/2\). It can be identified with a free
particle~\amscite{Wulfman10a}*{\S~3.8}. 

We can look for ladder operators in the
representation~(\ref{eq:schroedinger-rep-conf-der}--\ref{eq:shale-weil-der})
within the Lie algebra \(\algebra{h}_1\) in the form
\(\ladder[\rmp]{\pm}=a\tilde{X}+b\tilde{Y}\). This is possible if and only if
\begin{equation}
  \label{eq:compatib-parab}
  -b=\lambda a,\quad 0=\lambda b.
\end{equation}
The compatibility condition \(\lambda^2=0\) implies \(\lambda=0\)
within complex numbers. However such a ``ladder'' operator produces
only the zero shift on the eigenvectors, cf.~\eqref{eq:ladder-action}.

Another possibility appears if we consider the representation of the
Heisenberg group induced by dual-valued characters. On the
configuration space%
\index{configuration!space}%
\index{space!configuration} such a representation
is~\amscite{Kisil10a}*{(4.11)}: 
\begin{equation}
  \label{eq:schroedinger-rep-conf-par}
    [\uir{\rmp}{\chi}(s,x,y) f](q)= e^{2\pi\rmi x q}\left(
    \left(1-\rmp\myh (s-{\textstyle\frac{1}{2}}xy)\right) f(q)
    +\frac{\rmp\myh y}{2\pi\rmi} f'(q)\right).
\end{equation}
The corresponding derived representation of \(\algebra{h}_1\) is 
\begin{equation}
  \label{eq:schroedinger-rep-conf-der-par}
  \uir{p}{\myh}(X)=2\pi\rmi q,\qquad
  \uir{p}{\myh}(Y)=\frac{\rmp\myh}{2\pi \rmi} \frac{d}{dq},
  \qquad
  \uir{p}{\myh}(S)=-\rmp\myh I.
\end{equation}
However the Shale--Weil extension generated by this representation is
inconvenient.  It is better to consider the FSB--type parabolic
representation~\eqref{eq:dual-repres}%
\index{Fock--Segal--Bargmann!representation}%
\index{representation!Fock--Segal--Bargmann} on the phase space
induced by the same dual-valued character.
Then the derived representation of \(\algebra{h}_1\) is:
\begin{equation}
  \label{eq:schroedinger-rep-conf-der-par1}
  \uir{p}{\myh}(X)=-2\pi\rmi q-\frac{\rmp\myh}{4\pi\rmi}\partial_{p},\qquad
  \uir{p}{\myh}(Y)=-2\pi\rmi p+\frac{\rmp\myh}{4\pi\rmi}\partial_{q},
  \qquad
  \uir{p}{\myh}(S)=\rmp\myh I.
\end{equation}
An advantage of the FSB representation%
  \index{Fock--Segal--Bargmann!representation}%
  \index{representation!Fock--Segal--Bargmann} is that the
derived form of the parabolic Shale--Weil representation coincides
with the elliptic one~\eqref{eq:shale-weil-der-ell}.

Eigenfunctions with the eigenvalue \(\mu\) of the parabolic
Hamiltonian \(\tilde{B}+\tilde{Z}/2=q\partial_p\) have the form
\begin{equation}
  \label{eq:par-eigenfunctions}
  v_\mu (q,p)=e^{\mu p/q} f(q), \text{ with an arbitrary function }f(q).
\end{equation}

The linear equations defining the corresponding ladder operator
\(\ladder[\rmp]{\pm}=a\tilde{X}+b\tilde{Y}\) in the algebra
\(\algebra{h}_1\) are~\eqref{eq:compatib-parab}.  The compatibility
condition \(\lambda^2=0\) implies \(\lambda=0\) within complex numbers
again. Admitting dual numbers we have additional values
\(\lambda=\pm\rmp\lambda_1\) with \(\lambda_1\in\Space{C}{}\) with the
corresponding ladder operators
\begin{displaymath}
  \ladder[\rmp]{\pm}=\tilde{X}\mp\rmp\lambda_1\tilde{Y}=
  -2\pi\rmi q-\frac{\rmp\myh}{4\pi\rmi}\partial_{p}\pm 2\pi\rmp\lambda_1\rmi p= 
  -2\pi\rmi q+   \rmp\rmi( \pm 2\pi\lambda_1 p+\frac{\myh}{4\pi}\partial_{p}).
\end{displaymath}
For the eigenvalue \(\mu=\mu_0+\rmp\mu_1\) with \(\mu_0\),
\(\mu_1\in\Space{C}{}\) the
eigenfunction~\eqref{eq:par-eigenfunctions} can be rewritten as:
\begin{equation}
  \label{eq:par-eigenfunctions-1}
  v_\mu (q,p)=e^{\mu  p/q} f(q)= e^{\mu_0  p/q}\left(1+\rmp\mu_1
    \frac{p}{q}\right) f(q)
\end{equation}
due to the nilpotency of \(\rmp\).  Then the ladder action of
\(\ladder[\rmp]{\pm}\) is \(\mu_0+\rmp\mu_1\mapsto \mu_0+\rmp(\mu_1\pm
\lambda_1)\).  Therefore these operators are suitable for building
\(\algebra{sp}_2\)-modules with a one-dimensional chain of
eigenvalues.

Finally, consider the ladder operator for the same element \(B+Z/2\)
within the Lie algebra \(\algebra{sp}_2\), cf.
\S~\ref{sec:parab-ladd-oper}.  There is the only operator
\(\ladder[p]{\pm}=\tilde{B}+\tilde{Z}/2\) corresponding to complex
coefficients, which does not affect the eigenvalues.  However the dual
numbers lead to the operators
\begin{displaymath}
  \ladder[\rmp]{\pm}=\pm \rmp\lambda_2\tilde{A}+\tilde{B}+\tilde{Z}/2
  =
  \pm\frac{\rmp\lambda_2}{2}\left(q\partial_{q}-p\partial_{p}\right)+q\partial_{p}, 
  \qquad \lambda_2\in\Space{C}{}. 
\end{displaymath}
These operator act on eigenvalues in a non-trivial way.%
\index{dual!number|)}%
\index{number!dual|)}

\subsubsection{Similarity and Correspondence}
\label{sec:concl-simil-corr}

We wish to summarise our findings. Firstly, the appearance of
hypercomplex numbers%
\index{number!hypercomplex}%
\index{hypercomplex!number} in ladder operators for \(\algebra{h}_1\) follows
exactly the same pattern as was already noted for
\(\algebra{sp}_2\), see Rem.~\ref{re:hyper-number-necessity}:
\begin{itemize}
\item the introduction of complex numbers is a necessity for the
  \emph{existence} of ladder operators in the elliptic
  case;
\item in the parabolic case we need dual numbers to make
  ladder operators \emph{useful};
\item in the hyperbolic case double numbers are not required 
  neither for the existence or for the usability of ladder operators, but
  they do provide an enhancement. 
\end{itemize}
In the spirit of the Similarity and Correspondence
Principle~\ref{pr:simil-corr-principle}%
\index{principle!similarity and correspondence} we have the following
extension of Prop.~\ref{pr:ladder-sim-eq}: 
\begin{prop}
  \label{pr:ladder-sim-eq1}
  Let a vector \(H\in\algebra{sp}_2\) generates the subgroup \(K\),
  \(N'\) or \(\Aprime\), that is \(H=Z\), \(B+Z/2\), or
  \(2B\) respectively. Let \(\alli\) be the respective hypercomplex
  unit. Then the ladder operators \(\ladder[]{\pm}\)  satisfying to the
    commutation relation:
  \begin{displaymath}
    [H,\ladder[2]{\pm}]=\pm\alli \ladder{\pm}
  \end{displaymath}
  are given by:
  \begin{enumerate}
  \item Within the Lie algebra \(\algebra{h}_1\): \(\ladder{\pm}=\tilde{X}\mp\alli \tilde{Y}.\)
  \item Within the Lie algebra \(\algebra{sp}_2\): \(
    \ladder[2]{\pm}=\pm\alli \tilde{A} +\tilde{E}\).
  Here \(E\in\algebra{sp}_2\) is a linear combination of  \(B\) and
  \(Z\) with the properties:
  \begin{itemize}
  \item \(E=[A,H]\).
  \item \(H=[A,E]\).
  \item Killings form \(K(H,E)\)~\amscite{Kirillov76}*{\S~6.2} vanishes.
  \end{itemize}
  Any of the above properties defines the vector \(E\in\loglike{span}\{B,Z\}\)
  up to a real constant factor.
  \end{enumerate}
\end{prop}
It is worth continuing  this investigation and describing in details
hyperbolic and parabolic versions of FSB spaces.%
\index{Fock--Segal--Bargmann!space}%
\index{space!Fock--Segal--Bargmann}
\index{Heisenberg!group|)}%
\index{group!Heisenberg|)}
\index{ladder operator|)}%
\index{operator!ladder|)}%

\section{Open Problems}
\label{sec:open-problems}
A reader may already note numerous objects and results, which deserve
a further consideration. It may also worth to state some open problems
explicitly.  In this section we indicate several directions for
further work, which go through four main areas described in the paper.
\makeatletter \def\p@enumi{\thesubsection.}  \makeatother

\subsection{Geometry}
\label{sec:geometry-problems}
Geometry\index{geometry} is most elaborated area so far, yet many
directions are waiting for further exploration.
\begin{enumerate}
\item M\"obius transformations~\eqref{eq:moebius} with three types
  of hypercomplex units appear from the action of the group \(\SL\) on
  the homogeneous space \(\SL/H\)~\cite{Kisil09c}, where \(H\) is any
  subgroup \(A\), \(N\), \(K\) from the Iwasawa
  decomposition~\eqref{eq:iwasawa-decomp}%
  \index{Iwasawa decomposition}%
  \index{decomposition!Iwasawa}. Which other actions and
  hypercomplex numbers can be obtained from other Lie groups and
  their subgroups?
\item Lobachevsky geometry%
  \index{Lobachevsky!geometry}%
  \index{geometry!Lobachevsky} of the upper half-plane is extremely
  beautiful and well-developed subject~\citelist{\cite{Beardon05a}
    \cite{CoxeterGreitzer}}. However the traditional study is limited
  to one subtype out of nine possible: with the complex numbers for
  M\"obius transformation and the complex imaginary unit used in
  FSCc~\eqref{eq:FSCc-matrix}%
  \index{Fillmore--Springer--Cnops construction}%
  \index{construction!Fillmore--Springer--Cnops}.  The remaining eight
  cases shall be explored in various directions, notably in the
  context of discrete subgroups~\cite{Beardon95}.
\item The Fillmore-Springer-Cnops construction, see
  subsection~\ref{sec:cycles-as-invariant}, is closely related to the
  \emph{orbit method}%
  \index{orbit!method}%
  \index{method!orbits, of}~\cite{Kirillov99} applied to \(\SL\). An extension of
  the orbit method from the Lie algebra dual to matrices representing
  cycles may be fruitful for semisimple Lie groups.
\item \label{it:discrete-geom}
  A development of a discrete version of the geometrical notions
  can be derived from suitable discrete groups. A natural first
  example is the group
  \(\mathrm{SL}_2(\Space{F}{})\), where \(\Space{F}{}\) is a finite field,
  e.g. \(\Space[p]{Z}{}\) the field of integers modulo a prime \(p\).%
  \index{geometry!discrete}%
  \index{discrete!geometry}
\end{enumerate}

\subsection{Analytic Functions}
\label{sec:analytic-functions-problems}
It is known that in several dimensions there are different notions of
analyticity, e.g. several complex variables and Clifford analysis.%
\index{Clifford!algebra}%
\index{algebra!Clifford}
However, analytic functions of a complex variable are
usually thought to be the only options in a plane domain. The
following seems to be promising:
\begin{enumerate}
\item \label{it:hyp-functions}
  Development of the basic components of analytic function theory
  (the Cauchy integral%
  \index{integral!Cauchy}%
  \index{Cauchy!integral}, the Taylor 
  expansion, the Cauchy-Riemann%
\index{Cauchy-Riemann operator}%
\index{operator!Cauchy-Riemann} and Laplace equations%
  \index{Laplacian}, etc.) from the same
  construction and principles in the elliptic, parabolic and hyperbolic
  cases and respective subcases.
\item \label{it:bergman}
  Identification of Hilbert spaces of analytic functions of Hardy%
  \index{space!Hardy}%
  \index{Hardy!space} and
  Bergman types%
  \index{space!Bergman}%
  \index{Bergman!space}, investigation of their properties. Consideration of the
  corresponding Toeplitz operators%
  \index{Toeplitz!operator}%
  \index{operator!Toeplitz} and algebras generated by them.
\item Application of analytic methods to elliptic, parabolic and hyperbolic
  equations and corresponding boundary and initial values problems. 
\item \label{it:mult-dim-funct}
  Generalisation of the results obtained to higher dimensional
  spaces. Detailed investigation of physically significant cases of three
  and four dimensions.
\item There is a current interest in construction of analytic function
  theory on discrete sets. Our approach is ready for application
  to an analytic functions in discrete geometric set-up outlined in 
  item~\ref{it:discrete-geom} above.%
  \index{analytic!function on discrete sets}%
  \index{discrete!analytic function}
\end{enumerate}

\subsection{Functional Calculus}
\label{sec:functional-calculus-problems}
The functional calculus%
\index{functional!calculus}%
\index{calculus!functional} of a finite dimensional operator considered in
Section~\ref{sec:functional-calculus} is elementary but provides a
coherent and comprehensive treatment. It shall be extended to further
cases where other approaches seems to be rather limited.
\begin{enumerate}
\item Nilpotent and quasinilpotent operators have the most trivial
  spectrum possible (the single point \(\{0\}\)) while their structure
  can be highly non-trivial. Thus the standard spectrum is
  insufficient for this class of operators. In contract, the covariant
  calculus and the spectrum give complete description of nilpotent
  operators---the basic prototypes of quasinilpotent ones.  For
  quasinilpotent operators%
  \index{quasinilpotent!operator}%
  \index{operator!quasinilpotent} the construction will be more complicated
  and shall use analytic functions mentioned in \ref{it:hyp-functions}.
  
\item The version of covariant calculus described above is based on the
  \emph{discrete series}%
  \index{representation!discrete series} representations of \(\SL\) group and is particularly
  suitable for the description of the \emph{discrete spectrum}%
  \index{spectrum!discrete}%
  \index{discrete!spectrum} (note the
  remarkable coincidence in the names). 
  
  It is interesting to develop similar covariant calculi based on the
  two other representation series of \(\SL\): \emph{principal}%
  \index{representation!principal series} and \emph{complementary}%
  \index{representation!complementary series}~\cite{Lang85}. The
  corresponding versions of analytic function theories for
  principal~\cite{Kisil97c} and complementary series~\cite{Kisil05a}
  were initiated within a unifying framework. The classification of
  analytic function theories into elliptic, parabolic,
  hyperbolic~\cite{Kisil05a,Kisil06a} hints the following associative
  chains:
  \begin{center}
    \begin{tabular}{c@{---}c@{---}c}
      \textbf{Representations} &  \textbf{ Function Theory } & 
      \textbf{ Type of Spectrum }\\[1mm] \hline\hline
      discrete series%
      \index{representation!discrete series} &  elliptic   & discrete spectrum\\
      principal series%
      \index{representation!principal series}& hyperbolic & continuous spectrum\\
      complementary series%
      \index{representation!complementary series}& parabolic & residual spectrum
    \end{tabular}
  \end{center}

\item Let \(a\) be an operator with \(\spec a\in\bar{\Space{D}{}}\)
  and \(\norm{a^k}< C k^p\). It is typical to consider instead of
  \(a\) the \emph{power bounded} operator%
  \index{power bounded operator}%
  \index{operator!power bounded} \(ra\), where \(0<r< 1\),
  and consequently develop its \(\FSpace{H}{\infty}\) calculus.
  However such a regularisation is very rough and hides the nature of
  extreme points of \(\spec{a}\). To restore full information a
  subsequent limit transition \(r\rightarrow 1\) of the regularisation
  parameter \(r\) is required. This make the entire technique rather
  cumbersome and many results have an indirect nature.

  The regularisation \(a^k\rightarrow a^k/k^p\) is more natural and
  accurate for polynomially bounded operators. However it cannot be
  achieved within the homomorphic calculus Defn.~\ref{de:calculus-old}
  because it is not compatible with any algebra homomorphism. Albeit
  this may be achieved within the covariant
  calculus~Defn.~\ref{de:functional-calculus-new} and Bergman type space
  from~\ref{it:bergman}.

\item Several non-commuting operators are especially difficult to
  treat with functional calculus Defn.~\ref{de:calculus-old} or a joint
  spectrum. For example, deep insights on joint spectrum of commuting
  tuples~\cite{JTaylor72} refused to be generalised to non-commuting
  case so far.  The covariant calculus was initiated~\cite{Kisil95i}
  as a new approach to this hard problem and was later found useful
  elsewhere as well.  Multidimensional covariant
  calculus~\cite{Kisil04d} shall use analytic functions described
  in~\ref{it:mult-dim-funct}.

\item As we noted above there is a duality between the co- and
  contravariant calculi%
  \index{contravariant!calculus}%
  \index{calculus!contravariant}%
  \index{covariant!calculus}%
  \index{calculus!covariant} from Defins.~\ref{de:covariant-calculus}
  and~\ref{de:conravariant-calculus}. We also seen in
  Section~\ref{sec:functional-calculus} that functional calculus%
  \index{functional!calculus}%
  \index{calculus!functional} is an
  example of contravariant calculus and the functional model%
  \index{functional!model}%
  \index{model!functional} is a case of
  a covariant one. It is interesting to explore the duality between
  them further.
\end{enumerate}

\subsection{Quantum Mechanics}
\label{sec:quantum-mechanics}

Due to the space restrictions we only touched quantum mechanics,
further details can be found in~\citelist{\cite{Kisil96a}
  \cite{Kisil02e} \cite{Kisil05c} \cite{Kisil04a} \cite{Kisil09a}
  \cite{Kisil10a}}. In general, Erlangen approach is much more popular
among physicists rather than mathematicians. Nevertheless its
potential is not exhausted even there.

\begin{enumerate}
\item There is a possibility to build representation of the Heisenberg
  group%
  \index{Heisenberg!group}%
  \index{group!Heisenberg} using characters of its centre with values
  in dual and double numbers rather than in complex ones. This will
  naturally unifies classical mechanics, traditional QM and hyperbolic
  QM~\cite{Khrennikov08a}. In particular, a full construction of the
  corresponding Fock--Segal--Bargmann%
  \index{Fock--Segal--Bargmann!space}%
  \index{space!Fock--Segal--Bargmann} spaces would be of interest.
\item  Representations of nilpotent Lie groups with multidimensional
  centres in Clifford algebras%
  \index{Clifford!algebra}%
  \index{algebra!Clifford} as a framework for consistent quantum
  field%
  \index{field!quantum}%
  \index{quantum!field} theories based on De Donder--Weyl formalism%
  \index{De Donder--Weyl formalism}~\cite{Kisil04a}.
\end{enumerate}
\begin{rem}
  This work is performed within the ``Erlangen programme at large''
  framework~\cites{Kisil06a,Kisil05a}, thus it would be suitable to
  explain the numbering of various papers. Since the logical order may be
  different from chronological one the following numbering  scheme
  is used:
  \begin{center}
  \begin{tabular}{||c|p{.7\textwidth}||}
    \hline\hline
    Prefix&Branch description\\
    \hline\hline
    ``0'' or no prefix & Mainly geometrical works, within the classical
    field of Erlangen programme by F.~Klein, see~\citelist{\cite{Kisil05a} \cite{Kisil09c}}\\
    \hline 
    ``1'' & Papers on analytical functions theories and wavelets, e.g.~\cite{Kisil97c}\\
    \hline
    ``2'' & Papers on operator theory, functional calculi and
    spectra, e.g.~\cite{Kisil02a}\\ 
    \hline 
    ``3'' & Papers on mathematical physics, e.g.~\cite{Kisil10a}\\
    \hline\hline
  \end{tabular}    
  \end{center}
  For example, \cite{Kisil10a} is the first paper in the mathematical
  physics area. The present paper~\cite{Kisil11c} outlines the whole
  framework and thus does not carry a subdivision number. The on-line
  version of this paper may be updated in due course to reflect the
  achieved progress.
\end{rem}

\section*{Acknowledgement}
\label{sec:aknowledgement}

Material of these notes was lectured at various summer/winter schools
and advanced courses. Those presentations helped me to clarify ideas
and improve my understanding of the subject. In particular, I am grateful to
Prof.~S.V.~Rogosin and Dr.~A.A.~Koroleva for a kind invitation to the
Minsk Winter School in 2010, which were an exciting event. I would
like also to acknowledge support of
\href{http://maths.dept.shef.ac.uk/magic/index.php}{MAGIC group}
during my work on those notes.

{\small

\bibliography{arare,abbrevmr,akisil,ageometry,analyse,aphysics,algebra,aclifford,acombin}
}
\printindex

\end{document}